\documentclass[sn-aps]{sn-jnl}% American Physical Society (APS) Reference Style
\usepackage[utf8]{inputenc}
\usepackage{amsmath}
\usepackage{amsfonts}
\usepackage{amssymb}
\usepackage{mathtools}
\usepackage{lipsum}

\newcommand{\calO}{{\mathcal{O}}}

\newcommand{\calQ}{{\mathcal{Q}}}

\newcommand{\bbE}{{\mathbb{E}}}

\newcommand{\bbV}{{\mathbb{V}}}

\newcommand\E[1]{{\mathbb{E}\left[#1\right]}}
\newcommand\V[1]{{\mathbb{V}\left[#1\right]}}
\newcommand\MSE[1]{{\text{MSE}\left(#1\right)}}

\newcommand\cost[1]{{\text{cost}\left(#1\right)}}

\everymath{\displaystyle}

\usepackage{graphicx}
\usepackage{subcaption}

\jyear{2023}%

\theoremstyle{thmstyleone}%
\theoremstyle{thmstyletwo}%

\theoremstyle{thmstylethree}%

\begin{document}

\title[MF UQ for homogenization in CPFEM]{Multi-fidelity uncertainty quantification for homogenization problems in structure-property relationships from crystal plasticity finite elements}

%%=============================================================%%
%% Prefix   -> \pfx{Dr}
%% GivenName   -> \fnm{Joergen W.}
%% Particle -> \spfx{van der} -> surname prefix
%% FamilyName  -> \sur{Ploeg}
%% Suffix   -> \sfx{IV}
%% NatureName  -> \tanm{Poet Laureate} -> Title after name
%% Degrees  -> \dgr{MSc, PhD}

\author*[1]{\fnm{Anh} \sur{Tran}}\email{anhtran@sandia.gov}

\author[2]{\fnm{Pieterjan} \sur{Robbe}}\email{pmrobbe@sandia.gov}

\author[1]{\fnm{Theron} \sur{Rodgers}}\email{trodgers@sandia.gov}

\author[1]{\fnm{Hojun} \sur{Lim}}\email{hnlim@sandia.gov}

\affil[1]{\orgname{Sandia National Laboratories}, \orgaddress{\city{Albuquerque}, \postcode{87123}, \state{NM}, \country{USA}}}

\affil[2]{\orgname{Sandia National Laboratories}, \orgaddress{\street{Street}, \city{Livermore}, \postcode{94550}, \state{CA}, \country{USA}}}


\abstract{
Crystal plasticity finite element method (CPFEM) has been an integrated computational materials engineering (ICME) workhorse to study materials behaviors and structure-property relationships for the last few decades. These relations are mappings from the microstructure space to the materials properties space. 
Due to the stochastic and random nature of microstructures, there is always some uncertainty associated with materials properties, for example, in homogenized stress-strain curves. 
For critical applications with strong reliability needs, it is often desirable to quantify the microstructure-induced uncertainty in the context of structure-property relationships. 
However, this uncertainty quantification (UQ) problem often incurs a large computational cost because many statistically equivalent representative volume elements (SERVEs) are needed. 
In this paper, we apply a multi-level Monte Carlo (MLMC) method to CPFEM to study the uncertainty in stress-strain curves, given an ensemble of SERVEs at multiple mesh resolutions. 
By using the information at coarse meshes, we show that it is possible to approximate the response at fine meshes with a much reduced computational cost. 
We focus on problems where the model output is multi-dimensional, which requires us to track multiple quantities of interest (QoIs) at the same time. 
Our numerical results show that MLMC can accelerate UQ tasks around 2.23$\times$, compared to the classical Monte Carlo (MC) method, which is widely known as ensemble average in the CPFEM literature.
}




\maketitle

\section{Introduction}
\label{sec:Intro}

Designing materials with tailored properties requires a comprehensive understanding of the composition--process--structure--property relationship~\cite{olson1997computational,arroyave2019systems}, which has been studied extensively with four paradigms of materials science, namely experimental, theoretical, computational, and (scientific) machine learning~\cite{hey2009fourth,agrawal2016perspective}.
For polycrystalline metals and alloys, the materials design problem can often be formulated as optimization under uncertainty, where uncertainty comes from many sources. 
Some common epistemic uncertainty source includes observation noise, numerical solver tolerance, approximation error, besides the aleatory uncertainty source, which often manifests in the naturally random microstructure. 
Unfortunately, the materials design cycle is one of the major computational bottlenecks for transformative technologies~\cite{arroyave2019systems} due to their resource-intensive nature. 
To accelerate the materials design cycle, computational materials paradigms were developed and implemented in that context, across multiple length-scales and time-scales, with the overarching goal of modeling and simulating experimental and theoretical physics~\cite{olson1997computational}. 
Over the last few decades, multiple integrated computational materials engineering (ICME)~\cite{national2008integrated} models and simulations have been developed. ICME has become the third paradigm in materials science, the so-called computational paradigm~\cite{agrawal2016perspective}. 
The Materials Genome Initiative (MGI)~\cite{national2011materials} was created in 2011 in that scientific computing context, with the hope that ICME models can significantly reduce the research and development time and cost by leveraging the computational resource of high-performance computers. 
With the emerging field of machine learning, the revised MGI~\cite{de2019new} was updated in 2019 to include scientific machine learning (SciML) as the fourth paradigm of materials science~\cite{hey2009fourth,agrawal2016perspective}. 

Uncertainty quantification (UQ) plays several important roles in the materials design problem\cite{oden2010computer1,oden2010computer2}. 
First of all, many inverse UQ and optimization tools are used to perform deterministic or statistical ICME model calibration. 
Second, forward UQ tools are applied to quantify uncertainty associated with calibrated ICME models, in order to establish a reliable, robust, and predictive computing capability for ICME. 
Third, many UQ tools are currently developed and applied to quantify uncertainty for ICME-based SciML. 
Last but not least, multi-fidelity approaches also play a crucial role in ICME paradigm with multiple fidelity parameters, such as meshes, numerical integrators, iterative solvers, order of element. 
Therefore, from this perspective, UQ is naturally immersed in all four paradigms of materials science, from experiments to SciML. 

Within the composition-process-structure-property relationship, the microstructure is perhaps most uncertain due to its naturally inherent randomness, which can be associated with the aleatory uncertainty. 
It is well-known that the final microstructure of an alloy is generated ``from a very complex, process-specific, history-dependent sequence of transformation''~\cite{arroyave2019systems}. 
% Just a few years ago, inverting process-structure linkages are challenging and remains largely ``an unsolved problem due to the non-equilibrium, path-dependent nature of process-structure relations''~\cite{arroyave2019systems}. 
% Epistemic uncertainty, on the other hand, are often related to experimental noise, modeling assumptions, and numerical approximations. 
In this paper, we are concerned with quantifying the microstructure-induced uncertainty from the microstructure to the property space, where microstructures are represented by an ensemble of statistically equivalent representative volume elements (SERVEs) in a hierarchical multi-fidelity manner using crystal plasticity finite element model (CPFEM). 

Numerous works have been done for UQ in ICME in the last decades. 
Some notable works are summarized as follows. 
Nguyen et al.~\cite{nguyen2021bayesian} employed a Markov chain Monte Carlo to calibrate constitutive models for CPFEM in a Bayesian context. 
Hasan and Acar~\cite{hasan2022microstructure} developed a microstructure-sensitive design for performance optimization in titanium, aluminum, and galfenol. 
Tran et al.~\cite{tran2019quantifying,tran2022microstructure} applied stochastic collocation method, which is composed of generalized polynomial chaos expansion and sparse grid for phase-field simulation~\cite{tran2019quantifying} and CPFEM~\cite{tran2022microstructure}, respectively. 
Venkatraman et al.~\cite{venkatraman2023new} developed a three-step Bayesian protocol for model calibration and model-form UQ for CPFEM constitutive models in $\alpha+\beta$ titanium alloys. 
Rixner and Koutsourelakis~\cite{rixner2022self} formulated and developed a probabilistic, data-driven convolutional neural network to actively solve an inverse problem in structure-property linkage. 
Tran and Wildey~\cite{tran2020solving} applied a data-consistent inversion method~\cite{butler2018combining,butler2018convergence} to infer a distribution of microstructure features from a distribution of yield stress, where the push-forward density map is consistent with a target yield stress density through a heteroscedastic Gaussian process regression. 
Tran et al.~\cite{tran2023multi} employed multi-level and multi-index Monte Carlo to estimate effective yield strength and effective Young's modulus in a hierarchical multi-fidelity manner for a single quantity of interest (QoI). 
Rodgers et al.~\cite{rodgers2020three} simulated the process-structure with kinetic Monte Carlo~\cite{mitchell2023parallel} and the structure-property with CPFEM to study the effects of common microstructure in additive manufacturing. 
Tran et al.~\cite{tran2023high} employed an asynchronous parallel Bayesian optimization~\cite{tran2022aphbo} to calibrate phenomenological constitutive models for several materials systems. 
Ricciardi et al.~\cite{ricciardi2020uncertainty} applied the model discrepancy approach~\cite{kennedy2001bayesian} with Gaussian process using a visco-plastic self-consistent model to quantify uncertainty of homogenized stress-strain curve. 
Khalil et al.~\cite{khalil2021modeling} also applied an adaptive Metropolis-Hastings Markov chain Monte Carlo to calibrate constitutive models accounting for uncertainty in homogenized materials responses. 
Ghoreishi et al.~\cite{ghoreishi2018multi} proposed a weighted linear average approaches to combine multiple prediction from Gaussian process regression for homogenized stress-strain response for dual-phase microstructures.

% \textcolor{red}{
% Earlier work in Zhao et al.~\cite{zhao2022quantifying} incorporated measurement and parametric uncertainty to quantify the uncertainty of critical resolved share stress for hexagonal close-packed titanium alloys from nano-indentation. 
% Lim et al.~\cite{lim2019investigating} investigated the mesh sensitivity and polycrystalline representative volume element (RVE), where initial textures, hardening models, and boundary conditions are uncertain, and showed that an adequate polycrystalline RVE is obtained by capturing 1000 or more grains. 
% Fernadez et al.~\cite{fernandez2018estimating} utilized Bayesian inference to quantify the uncertainty in stress-strain curves, where model parameters are treated as random variables. 
% Tallman et al.~\cite{tallman2019gaussian,tallman2020uncertainty} applied Gaussian process regression and the Materials Knowledge System framework to predict a set of homogenized materials properties with uncertainty from a distribution function for crystallography orientations and textures.
% The inductive design exploration method (IDEM) \cite{ellis2017application,mcdowell2009integrated,choi2008inductive} has been introduced as a materials design methodology to identify feasible and robust design for microstructure features, which has been broadly applied to many practical problems. 
% }

In this paper, we extend our previous work in~\cite{tran2023multi} from a single QoI to multiple QoIs, employing the multi-level Monte Carlo (MLMC) method to quantifying uncertainty in the stress-strain curves for a magnesium alloy. 
The rest of the paper is organized as follows. 
In section~\ref{sec:MLMC}, we describe the MLMC method in the context of CPFEM applications. 
In section~\ref{sec:CPFEM}, we describe the CPFEM case study for magnesium. Our main numerical results, including a comparison between the MLMC method and the brute-force MC method, are discussed in Section~\ref{sec:results}. 
In sections~\ref{sec:Discussion} and~\ref{sec:Conclusion}, we discuss our results and formulate a conclusion.

% \section{Quantifying microstructure-induced variability of materials properties}
\section{Multi-level Monte Carlo for microstructure-induced uncertainty quantification}
\label{sec:MLMC}

\subsection{Microstructure from a statistical perspective}

% \textcolor{red}{
% We denote the uncertain microstructure of the material under consideration by $m$, where $\bsx \in \calD \subseteq \bbR^d$ denotes a spatial location, with $\calD$ the computational domain where the $d$-dimensional microstructure is defined, and where $m \in \mathcal{M}$ denotes an outcome, with $\mathcal{M}$ the sample space of the underlying probability space $(\mathcal{M}, \calF, p)$. Thus, the microstructure $m$ is spatially varying with the location $\bsx$ and depends on a random parameter $m$ representing a possible outcome of the microstructure. Note that we made the assumption that all microstructure SERVEs follow the same underlying distribution according to the probability density function (PDF) $p(m)$. The rationale behind this assumption is that microstructure reconstruction should theoretically produce statistically equivalent microstructure SERVEs.
% }

% \pr{I would really argue to revert the changes in notation here... The way it is currently written is mathematically not correct imho... Since we are stating in the introduction that this is a ``mathematically rigorous approach to CPFEM'', we should make sure that this is correct. I think a good way is to consider the microstructure $m(\bsx, m)$ as a random field / random process, i.e., it is spatially varying and thus depends on an $\bsx$-coordinate, as well as a hidden stochastic parameter $m$. For fixed location $\bsx$, the SERVE $m(\bsx, \cdot)$ is a random variable, while for fixed $m$, $m(\cdot, m)$ is a sample of the SERVE...}

Given a probability space $(\Omega, \mathcal{A}, \mathbb{P})$ and $\omega \in \Omega$, we consider an uncertain microstructure $m(\boldsymbol{x}, \omega)$ defined on a bounded domain $D \subset \mathbb{R}^d$, with $d=2$ or $d=3$. The notation $m(\boldsymbol{x}, \omega)$ entails that the microstructure depends on both the spatial variable $\boldsymbol{x} \in D$ as well as the stochastic sample $\omega$ from the sample space of the corresponding probability space. For a fixed $\omega$, the microstructure $m(\boldsymbol{x}, \cdot)$ is a deterministic but space-dependent function. For a fixed location $\boldsymbol{x}$, the microstructure $m(\cdot, \omega)$ is a random variable. We are interested in computing statistics of a certain QoI $Q$ that is defined by a mapping $f$ from the microstructure space to the homogenized material property, i.e., $Q = f(m(\cdot, \omega))$. $m(\boldsymbol{x}, \omega)$ is often known as the statistically equivalent representative volume element (SERVE) in the materials science literature, sampled independently and identically (i.i.d.) from statistical microstructure measure $\mathbb{P}(\omega)$ that describes all statistical microstructure descriptors, including shape, size, morphology, neighboring, chord length, and crystal orientation distribution functions. 

For the remainder of this paper, we will be interested in computing the first-order moment or expected value of the QoI $Q(m)$, defined as
\begin{equation}
\E{Q} \coloneqq \int_\Omega Q(m(\cdot, \omega)) \mathrm{d} \mathbb{P}(\omega).
\label{eq:integral}
\end{equation}
Since the microstructure $m(\cdot, \omega)$ is not directly accessible, the QoI is often approximated by the quantity $Q_L = f(m_L(\cdot, \omega))$, where $m_L$ is a finite-dimensional approximation to the microstructure $m$ at a particular mesh resolution level $L$, typically the solution of a microstructure reconstruction problem. 
For readers interested in statistical microstructure descriptors and microstructure reconstruction problems, we refer to some notable works in the literature, for example, \cite{groeber2008framework1,groeber2008framework2,torquato2002statistical,bostanabad2018computational}.

\subsection{Monte Carlo method}

Given an ensemble of microstructure of $N$ SERVEs, denoted as $\{m(\cdot, \omega^{(n)})\}_{n=1}^N$ with corresponding predictions for the material property of interest $\{Q_L(m(\cdot, \omega^{(n)}))\}_{n=1}^N$ where $Q_L$ is the map from microstructure to properties (i.e. a CPFEM simulation), we can approximate Equation~\eqref{eq:integral} by the average
\begin{equation}
\calQ_{\text{MC}} \coloneqq \frac{1}{N} \sum_{n=1}^{N} Q_L(m(\cdot, \omega^{(n)})),
\label{eq:mc_method}
\end{equation}
where the SERVE $m(\cdot, \omega)$ is sampled from the measure $\mathbb{P}(\omega)$, i.e. $\omega \sim \mathbb{P}(\omega)$. 
The ensemble average approach in Equation~\eqref{eq:mc_method}, also known as the Monte Carlo (MC) method, is widely used in the CPFEM literature, see, e.g., \cite{paulson2017reduced,teferra2018random,paulson2018data}. 
For each microstructure SERVE $m(\cdot, \omega^{(n)})$, CPFEM is deployed to evaluate the corresponding material property $Q_L(m(\cdot, \omega^{(n)}))$.

It is natural to propose the average of an ensemble of material properties extracted from the microstructures $\{m(\cdot, \omega^{(n)})\}_{n=1}^N$ to approximate the expected value in Equation~\eqref{eq:integral}. Since the sequence of microstructure SERVEs are drawn from the same underlying sample space $\Omega$ with probability $\mathbb{P}(\omega)$, we have that the expected value $\bbE[Q_L(m^{(1)})] = \dots = \bbE[Q_L(m^{(N)})] = \bbE[Q_L]$, and the strong law of large numbers guarantees that $\calQ_{\text{MC}} \to \E{Q_L}$ almost surely as the number of microstructure realizations (i.e. the number of SERVE) $N$ goes  to infinity, see \cite{robert2013monte}.

There are two sources of error in the MC estimator in Equation~\eqref{eq:mc_method}: a stochastic error, present because we approximate the expected value by an average, and a bias, present because samples of $Q(m)$ are approximated by samples of $Q_L(m)$. 
% These two contributions become apparent in the expression for the mean square error (MSE) of the MC estimator. 
The mean square error (MSE) of the estimator can be decomposed into to the sum of these two terms:
\begin{equation}
\begin{array}{ll}
    \MSE{\calQ_{\text{MC}}} &\coloneqq \E{ (\calQ_{\text{MC}} - \E{Q})^2 } \\
    &= \E{((\calQ_{\text{MC}} - \E{\calQ_{\text{MC}}}) + (\E{\calQ_{\text{MC}}} - \E{Q}))^2} \\
    &= \E{(\calQ_{\text{MC}} - \E{\calQ_{\text{MC}}})^2} + (\E{Q_L - Q})^2 \\
    &= \V{\calQ_{\text{MC}}} + (\E{Q_L - Q})^2 ,
    \label{eq:mc_mse}
\end{array}
\end{equation}
where the cross-product term vanishes because the MC estimator is an unbiased estimator for $Q_L$, i.e., $\E{\calQ_{\text{MC}}} = \E{Q_L}$. The first term in Equation~\eqref{eq:mc_mse} is the variance of the estimator and represents the stochastic error. Because we assume the ensemble is uncorrelated, the variance can be written as
\begin{equation}
\label{eq:mc_variance_decays}
\V{\calQ_{\text{MC}}} = \frac{1}{N^2} \sum_{n=1}^{N} \V{Q_L} = \frac{\V{Q_L}}{N}.
\end{equation}
The variance for the MC estimator $\V{\calQ_{\text{MC}}}$ decays as $\calO(N^{-1})$ and can be reduced by considering more SERVEs, i.e. increasing $N$. The second term in Equation~\eqref{eq:mc_mse} is the square of the bias. It can be reduced by increasing the level of resolution $L$, i.e., by decreasing the mesh size, which leads to the notion of $h$-convergence in FEM~\cite{babuvska1990p}. 
% In term of $hp$ convergence for FEM (including CPFEM), decreasing the mesh size 

% If we require an MSE smaller than or equal to $\epsilon^2$, a sufficient condition is
If we restrict the MSE to be less than or equal to $\epsilon^2$, then a sufficient condition could be
\begin{equation}
\frac{\V{Q_L}}{N} \leq \frac{\epsilon^2}{2} \quad \text{and} \quad \vert \E{Q_L - Q} \vert \leq \frac{\epsilon}{\sqrt{2}}.
\end{equation}
Hence, the number of microstructure SERVEs $N$ scales as $\calO(\epsilon^{-2})$. 
Let $C_L$ be the cost of a single model evaluation, the total computational cost of the MC estimator in Equation~\eqref{eq:mc_method} can be expressed as
\begin{equation}
    \mathrm{cost}(\calQ_{\text{MC}}) = N C_L.
\end{equation}
Since $C_L$ is considered as a constant, the cost complexity of the MC estimator $\mathrm{cost}(\calQ_{\text{MC}})$ is the same as the number of SERVEs, i.e. $N$ $\calO(\epsilon^{-2})$.

\subsection{Multi-level Monte Carlo}

% \hl{Need to simplify this section.}

% The central idea in MLMC sampling is that we do not sample from a single approximation $Q_L$ for the QoI, but instead compute samples on a hierarchy of approximations $\{Q_\ell\}_{\ell=0}^L$ for the quantity of interest $Q$. In the context of CPFEM, this hierarchy corresponds to an approximation for the material parameter computed from a microstructure defined on a sequence of meshes with increasing resolution levels, where level $\ell=0$ corresponds to the cheapest approximation with the coarsest mesh size, and level $\ell = L$ corresponds to the most expensive approximation with the finest mesh size.  
The main idea of MLMC sampling~\cite{giles2015multilevel} is that instead of sampling from only one approximation $Q_L$, we sample from a hierarchy of approximation $\{Q_\ell\}_{\ell=0}^L$ for the QoI, denoted as $Q$. 
$\ell=0$ indicates the lowest level of fidelity and $\ell = L$ indicates the highest level of fidelity. 
From the multi-fidelity perspective, the intuition is to approximate the high-fidelity data using low-fidelity data by their correlation and hopefully, reducing the computational cost while retaining the approximation accuracy.  
In the context of this paper, $\ell = 0$ corresponds to the coarsest mesh size and $\ell = L$ corresponds to the finest mesh size, respectively, and $\{Q_\ell\}_{\ell=0}^L$ is a sequence of meshes.

% Because the expected value is a linear operator, we have that
Thanks to the linear property of the expectation operator, we can decompose the expectation at the highest level of fidelity into a telescoping sum as
\begin{equation}
\label{eq:tel_sum}
    \E{Q_L} = \E{Q_0} + \sum_{\ell=1}^L \E{Q_\ell - Q_{\ell-1}} = \sum_{\ell=0}^L \E{\Delta Q_\ell}
\end{equation}
where the backward difference is
\begin{equation}
\label{eq:mlmc_diff}
    \Delta Q_\ell \coloneqq \begin{cases}
    Q_{\ell} - Q_{\ell-1}, & \text{for } \ell > 0, \\
    Q_\ell, & \text{for } \ell = 0. \\
    \end{cases}
\end{equation}
% Using an independent MC estimator for each of the $L+1$ terms in the right-hand side of~\eqref{eq:tel_sum}, we obtain the MLMC estimator
Using the MC estimator to estimate $\E{\Delta Q_\ell}$ for $0 \leq \ell \leq L$, the MLMC estimator can be written as
\begin{equation}
\label{eq:mlmc}
    \calQ_{\text{MLMC}} \coloneqq \sum_{\ell=0}^L \frac{1}{N_\ell} \sum_{n=1}^{N_\ell} \Delta Q_\ell(m(\cdot, \omega^{(n)})).
\end{equation}
% In effect, this means that we use an ensemble of microstructure SERVEs $\{m^{(n)}\}_{n=1}^{N_\ell}$ on each level $\ell=0, 1, \ldots, L$ to estimate the expected values on the right-hand side of~\eqref{eq:tel_sum}, where we assume that the microstructure instances on each level are mutually independent.
Theoretically, this means an ensemble of microstructure SERVEs $\{m^{(n)}\}_{n=1}^{N_\ell}$ is used to estimate $\Delta Q_\ell$ from the telescoping sum in Equation~\eqref{eq:tel_sum}, for $0 \leq \ell \leq L$. 
% on each level $\ell=0, 1, \ldots, L$ to estimate the expected values on the right-hand side of~\eqref{eq:tel_sum}

The MLMC estimator in Equation~\eqref{eq:mlmc} is still an unbiased estimator for $\E{Q_L}$, i.e., $\E{\calQ_{\text{MLMC}}} = \E{\Delta Q_L}$
% \begin{equation}
%     \E{\calQ_{\text{MLMC}}} = \sum_{\ell=0}^L \frac{1}{N_\ell} \sum_{n=1}^{N_\ell} \E{\Delta Q_\ell} = \sum_{\ell=0}^L \E{\Delta Q_\ell} = \E{\Delta Q_L}
% \end{equation}
with variance
\begin{equation}
    % \V{\calQ_{\text{MLMC}}} = \sum_{\ell=0}^L \frac{1}{N_\ell^2} \sum_{n=1}^{N_\ell} \V{\Delta Q_\ell} = \sum_{\ell=0}^L \frac{\V{\Delta Q_\ell}}{N_\ell}.
    \V{\calQ_{\text{MLMC}}} = \sum_{\ell=0}^L \frac{\V{\Delta Q_\ell}}{N_\ell}.
\end{equation}

Decomposing the MSE into the bias and variance as in~\eqref{eq:mc_mse}, we have 
\begin{equation}
\label{eq:mlmc_mse}
    \MSE{\calQ_{\text{MLMC}}} = \V{\calQ_{\text{MLMC}}} + (\E{Q_L - Q})^2 = \sum_{\ell=0}^L \frac{\V{\Delta Q_\ell}}{N_\ell} + (\E{Q_L - Q})^2.
\end{equation}
% Again, the MSE consists of two terms: the variance of the estimator and the square of the bias. Note that the bias of the MLMC estimator is the same as the bias of the MC estimator.
It is worthy to point out that the bias of the MLMC estimator is the same as the bias of the MC estimator.

% A crucial observation is that, instead of estimating the expected value $\E{Q_\ell}$ directly on level $\ell$, it is much cheaper to estimate the expected value of the backward difference $\E{\Delta Q_\ell}$, if the random variables $Q_\ell$ and $Q_{\ell-1}$ are strongly positively correlated, i.e.,
% \begin{equation}
% \begin{array}{lll}
% % \begin{align}
%     \V{\Delta Q_\ell} &= \bbV[Q_\ell - Q_{\ell-1}] \\
%     &= \V{Q_\ell} + \V{Q_{\ell-1}} - 2 \text{cov}(Q_\ell, Q_{\ell-1}) \\
%     &\ll \V{Q_\ell} + \V{Q_{\ell-1}},
% % \end{align}
% \end{array}
% \label{eq:StrongCorr}
% \end{equation}
% where $\text{cov}(Q_\ell, Q_{\ell-1}) = \rho_{\ell, \ell-1} \sqrt{\V{Q_\ell} \V{Q_{\ell-1}}}$ is the covariance between $Q_\ell$ and $Q_{\ell-1}$ and $\rho_{\ell, \ell-1}$ is the Pearson correlation coefficient. 

If the $Q_\ell$ and $Q_{\ell-1}$ are strongly positively correlated, then
\begin{equation}
\begin{array}{lll}
% \begin{align}
    \V{\Delta Q_\ell} &= \bbV[Q_\ell - Q_{\ell-1}] \\
    &= \V{Q_\ell} + \V{Q_{\ell-1}} - 2 \text{cov}(Q_\ell, Q_{\ell-1}) \\
    &\ll \V{Q_\ell} + \V{Q_{\ell-1}},
% \end{align}
\end{array}
\label{eq:StrongCorr}
\end{equation}
where $\text{cov}(Q_\ell, Q_{\ell-1}) = \rho_{\ell, \ell-1} \sqrt{\V{Q_\ell} \V{Q_{\ell-1}}}$ is the covariance between $Q_\ell$ and $Q_{\ell-1}$ and $\rho_{\ell, \ell-1}$ is the Pearson correlation coefficient. 
Notice that the backward difference $\E{\Delta Q_\ell}$ is computed on the \textit{same} input microstructure $m(\cdot, \omega^{(n)})$. 
% In order to ensure this strong correlation, it is important to note that the difference $\Delta Q_\ell$ in~\eqref{eq:mlmc} is evaluated for the \textit{same} input microstructure $m(\cdot, \omega^{(n)})$. In the context of CPFEM, this means that the difference is computed from the material parameter prediction for the \emph{same} underlying SERVE $m$, but on two different mesh resolution levels $\ell$ and $\ell-1$ in the hierarchy.
As $\ell \rightarrow \infty$, we recover the $h$-convergence from the asymptotic analysis of FEM when the mesh size approaches zero, i.e. $h \to 0$. 
% As the level parameter $\ell \rightarrow \infty$, we expect the approximations $Q_\ell$ to converge towards the true QoI $Q$ in mean square sense, i.e., $\V{\Delta Q_\ell} \rightarrow 0$ as $\ell \rightarrow \infty$. In effect, this means that fewer model evaluations are required in the successive MC estimators for the difference $\Delta Q_\ell$ with increasing $\ell$. Under this assumption, we find that most samples will be taken on level $\ell=0$, where model evaluations are cheap, and fewer samples are required on the higher levels, where model evaluations are increasingly more expensive. Often, only a handful of samples with the highest resolution level are required. 
% Compare this to the MC method outlined in \cref{sec:mc}, where all samples are taken on the same high-resolution level.

If we require an MSE smaller than or equal to $\epsilon^2$, a sufficient condition is
\begin{equation}
\label{eq:mlmc_constraint}
    \sum_{\ell=0}^L \frac{\V{\Delta Q_\ell}}{N_\ell} \leq \frac{\epsilon^2}{2} \quad \text{and} \quad \vert \E{Q_L - Q} \vert \leq \frac{\epsilon}{\sqrt{2}}.
\end{equation}
% An expression for the required number of SERVEs $N_\ell$ on each level $\ell = 0, 1, \ldots, L$ can be obtained by minimizing the total cost of the MLMC estimator, taking into account the above constraint on the variance of the estimator. 
The total cost of the MLMC estimator is
\begin{equation}
    \mathrm{cost}(\calQ_{\text{MLMC}}) = \sum_{\ell=0}^L N_\ell \Delta C_\ell,
\end{equation}
where $\Delta C_\ell$ is the sampling cost of the backward difference $\Delta Q_\ell$. 
The optimal number of samples is~\cite{giles2015multilevel}
\begin{equation}
\label{eq:optimalSamples_mlmc}
N_{\ell} = \frac{2}{\epsilon^2}\sqrt{\frac{\V{\Delta Q_\ell}}{\Delta C_{\ell}}} \left( \sum_{\ell=0}^L \sqrt{\V{\Delta Q_{\ell}} \Delta C_{\ell}} \right),
\end{equation}
% see~\cite{giles2015multilevel} for details on the derivation. 
Numerically, $N_\ell$ in~\eqref{eq:optimalSamples_mlmc} is rounded up to the nearest integer, which in turn, may increase the cost of the estimator by at most one sample per level.

% In~\cite{cliffe2011multilevel}, a theoretical bound for the asymptotic cost complexity of the MLMC estimator is provided. 
Assuming that the expectation, variance, and sampling cost are bounded as
\begin{align}
\vert \E{\Delta \calQ_{\ell}} \vert &\leq c_1 \; 2^{-\alpha \ell}, \label{eq:C1}\tag{C1} \\
\V{\Delta \calQ_{\ell}} &\leq c_2 \; 2^{-\beta \ell} \text{ and} \label{eq:C2}\tag{C2} \\
\Delta C_\ell &\leq c_3 \; 2^{\gamma \ell} \label{eq:C3}\tag{C3}
\end{align}
with $2\alpha \geq \min(\beta,\gamma)$, the asymptotic cost complexity of the MLMC estimator is \cite{cliffe2011multilevel}
% we have that
\begin{equation}
\label{eq:mlmc_theorem}
\cost{\calQ_{\text{MLMC}}} \leq
\begin{cases}
c_4 \;\epsilon^{-2} & \text{if } \beta > \gamma, \\
c_4 \;\epsilon^{-2} (\log \epsilon)^2 & \text{if } \beta = \gamma, \\
c_4 \;\epsilon^{-2 - (\gamma - \beta) / \alpha} & \text{if } \beta < \gamma.
\end{cases}
\end{equation}
The extension of MLMC methods to multi-output MLMC is simply done by imposing the convergence criteria on the $L_\infty$ norm of the vector-valued output. 

% \subsection{Multi-output multi-level Monte Carlo}

% \hl{To be added: descriptions / algorithms}

\subsection{Integrated workflow: \texttt{MultilevelEstimators.jl} + \texttt{DREAM.3D} + \texttt{DAMASK}}

\begin{figure}[!htbp]
\centering
\includegraphics[width=\textwidth]{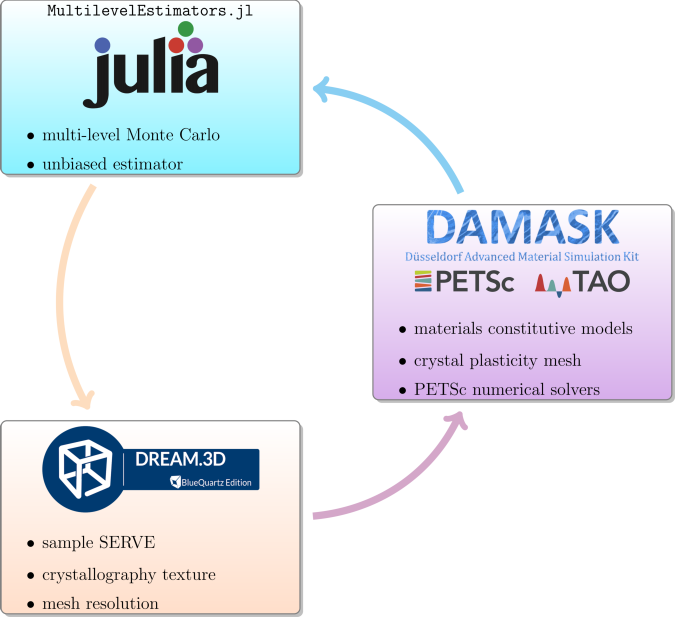}
\caption{Multi-fidelity uncertainty quantification workflow for CPFEM.}
\label{fig:uq_workflow}
\end{figure}

Figure~\ref{fig:uq_workflow} describes the integrated workflow used in this paper. 
In this workflow, three different software packages are integrated to create a UQ framework for CPFEM problems. 
At each iteration, \texttt{MultilevelEstimators.jl}~\cite{robbe2017multi,robbe2019recycling} requests an evaluation of the CPFEM model at fidelity levels $\ell$ and $\ell-1$. \texttt{DREAM.3D}~\cite{groeber2014dream} is then employed to generate a microstructure SERVE on these mesh resolution levels. \texttt{DAMASK}~\cite{roters2019damask} uses the generated microstructure geometries and computes the QoI using a \texttt{PETSc}~\cite{balay1997efficient,petsc-web-page} backend. The workflow is coupled together using a combination of Python and shell scripts.

% \pjr{\sout{For example,}} \texttt{DAMASK}~\cite{roters2019damask} is mainly programmed in Fortran 90 with pre-/post-processing pipeline in Python. 
% Among many numerical solvers, we chose to build \texttt{DAMASK} with \texttt{PETSc}~\cite{balay1997efficient,petsc-web-page}, which is programmed in C. 
% \texttt{MultilevelEstimators.jl}~\cite{robbe2017multi,robbe2019recycling} is programmed in Julia, while \texttt{DREAM.3D}~\cite{groeber2014dream} is programmed in C++, JavaScript, and Python. 
% Bash/Shell scripts are used to bind results from one package to another. 

\section{Crystal plasticity finite element for magnesium}
\label{sec:CPFEM}

\subsection{Phenomenonlogical constitutive model}

% For small deformations, the elasto-plastic decomposition can be computed additively, whereas for large deformations, a multiplicative decomposition of deformation gradient is more appropriate, i.e.,
In this section, we followed ~\cite{roters2010overview,roters2011crystal,roters2019damask} to summarize the basic of CPFEM with phenomenological constitutive model. 
The multiplicative elasto-plastic decomposition of deformation gradient for large deformations reads as
\begin{equation}
\mathbf{F} = \mathbf{F}_\text{e} \cdot \mathbf{F}_\text{p},
\label{eq:strain_decomp}
\end{equation}
and the elasto-plastic decomposition of the velocity gradient is
\begin{equation}
\mathbf{L} = \dot{\mathbf{F}} \cdot \mathbf{F}^{-1} = \dot{\mathbf{F}}_\text{e} \cdot \mathbf{F}_\text{e}^{-1} + \mathbf{F}_\text{e} \cdot \dot{\mathbf{F}}_\text{p} \cdot \mathbf{F}_\text{p} \cdot \mathbf{F}_\text{e}^{-1}
= \mathbf{L}_\text{e} + \mathbf{F}_\text{e} \cdot \mathbf{L}_\text{p} \cdot \mathbf{F}_\text{e}^{-1},
\end{equation}
$\mathbf{L}_\text{p}$ and $\mathbf{L}_\text{e}$ are the plastic and elastic velocity gradient, respectively. 
% The second Piola-Kirchhoff stress tensor $\mathbf{S}$, which is a symmetric second-order tensor defined in the intermediate configuration, is given by
% \begin{equation}
% \mathbf{S}
% = \frac{\mathbb{C}}{2} : (\mathbf{F}_\text{e}^T \mathbf{F}_\text{e} - \mathbf{I})
% = \mathbb{C} : \bse_\text{e}
% = J \mathbf{F}^{-1} \cdot \mathbf{\sigma} \cdot \mathbf{F}^{-T},
% \end{equation}
% where $\mathbb{C}$ is the elasticity fourth-order tensor, $\mathbf{F}_\text{e}$ is the elastic deformation gradient, $\mathbf{F}_\text{p}$ is the plastic deformation gradient \cite{roters2010overview}, $\bse_\text{e} = \frac{1}{2}\left( \mathbf{F}^T_e \mathbf{F}_\text{e} - \mathbf{I} \right)$ is the elastic Green's Lagrangian strain and $\sigma$ is the Cauchy stress tensor (cf. \cite{roters2011crystal}, Section 3.3).
% The evolution of the inelastic deformation gradient $\mathbf{F}_\text{p}$ is given in terms of their respective velocity gradients $\mathbf{L}_\text{p}$ by the flow rules
The flow rules models the evolution of the inelastic deformation gradient $\mathbf{F}_\text{p}$ as
\begin{align}
\dot{\mathbf{F}}_\text{p} = \mathbf{L}_\text{p} \mathbf{F}_\text{p},\\
\end{align}
where the plasticity velocity gradient in the intermediate configuration $\mathbf{L}_\text{p}$ is determined by
\begin{equation}
\mathbf{L}_\text{p}
= \dot{\mathbf{F}}_\text{p} \cdot \mathbf{F}_\text{p}^{-1}
= \sum_{\alpha} \dot{\gamma}^{\alpha} \left( \mathbf{s}_{\text{s}}^{\alpha} \otimes \mathbf{n}_{\text{s}}^{\alpha} \right).
\end{equation}
$\mathbf{s}_{\text{s}}^{\alpha}$ is the unit vector along the slip direction and $\mathbf{n}_{\text{s}}^{\alpha}$ is the unit vectors normal to the slip plane. 
% The driving force $\tau^{\alpha}$ for $\dot{\gamma}^{\alpha}$ is given by the Schmid law as
% \begin{equation}
% \tau^{\alpha} = \mathbf{M}_\text{p} \cdot \left( \mathbf{s}_{\text{s}}^{\alpha} \otimes \mathbf{n}_{\text{s}}^{\alpha} \right),
% \end{equation}
% where $\mathbf{M}_\text{p}$ is the Mandel stress in the plastic configuration, calculated from the second Piola-Kirchhoff stress $\mathbf{S}$. 

% A phenomenological crystal plasticity constitutive model used for face-centered cubic (FCC) crystals was first proposed by Hutchinson~\cite{hutchinson1976bounds} and extended for deformation twinning by Kalidindi~\cite{kalidindi1998incorporation}.
% The plastic component is parameterized in terms of resistance $\xi$ on $N_{\text{s}}$ slip and $N_{\text{tw}}$ twin systems. 
Plasticity is considered in terms of resistance $\xi$ on $N_{\text{s}}$ slip and $N_{\text{tw}}$ twin systems. 
The resistances on $\alpha=1,\dots,N_{\text{s}}$ slip systems evolve as
% to a system-dependent saturation value and depend on shear on slip and twin systems according to
\begin{equation}
\begin{array}{lll}
\dot{\xi}^{\alpha} &=& h_0^{\text{s-s}} \left( 1 + c_1 \left(f^{\text{tot}}_{\text{tw}}\right)^{c_2} \right) (1 + h^{\alpha}_{\text{int}}) \left[ \sum_{\alpha' = 1}^{N_{\text{s}}} \vert \dot{\gamma}^{\alpha'} \vert \left\vert 1 - \frac{\xi^{\alpha'}}{\xi^{\alpha'}_{\infty}} \right\vert^a \text{sgn}\left( 1 - \frac{\xi^{\alpha'}}{\xi^{\alpha'}_{\infty}} \right) h^{\alpha \alpha'} \right] \\
&& + \sum_{\beta'=1}^{N_{\text{tw}}} \dot{\gamma}^{\beta'} h^{\alpha \beta'} ,
\label{eq:slipping}
\end{array}
\end{equation}
where
$f^{\text{tot}}_{\text{tw}}$ is the total twin volume fraction,
$h$ is the matrices of the slip-slip and slip-twin interactions,
$h_0^{\text{s-s}}$, $h_{\text{int}}$, $c_1$, $c_2$ are fitting parameters, $\xi_{\infty}$ is the saturated resistance.

The resistances on the $\beta = 1, \dots, N_{\text{tw}}$ twin systems evolve similarly:
\begin{equation}
\dot{\xi}^{\beta} = h_0^{\text{tw-s}} \left( \sum_{\alpha=1}^{N_{\text{s}}} \vert\gamma_{\alpha}\vert \right)^{c_3} \left( \sum_{\alpha'=1}^{N_{\text{s}}} \vert\dot{\gamma}^{\alpha'}\vert h^{\beta \alpha'} \right)
+ h_0^{\text{tw-tw}} \left(f^{\text{tot}}_{\text{tw}} \right)^{c_4} \left( \sum_{\beta'=1}^{N_{\text{tw}}} \dot{\gamma}^{\beta'} h^{\beta \beta'} \right),
\label{eq:twinning}
\end{equation}
where $h_0^{\text{tw-s}}$, $h_0^{\text{tw-tw}}$, $c_3$, and $c_4$ are fitting parameters.
Shear on each slip system evolves as
\begin{equation}
\dot{\gamma}^{\alpha} = (1 - f^{\text{tot}}_{\text{tw}}) \dot{\gamma_0}^{\alpha} \left\vert \frac{\tau^{\alpha}}{\xi^{\alpha}} \right\vert^n \text{sgn}(\tau^{\alpha}).
\end{equation}
where slip due to mechanical twinning is 
% accounting for the unidirectional character of twin formation is computed slightly differently,
\begin{equation}
\dot{\gamma} = (1 - f^{\text{tot}}_{\text{tw}}) \dot{\gamma_0} \left\vert \frac{\tau}{\xi} \right\vert^n \mathcal{H}(\tau).
\end{equation}
$\mathcal{H}$ is the Heaviside step function. 
The total twin volume is
\begin{equation}
f^{\text{tot}}_{\text{tw}} = \max\left(1.0, \sum_{\beta=1}^{N_{\text{tw}}} \frac{\gamma^{\beta}}{\gamma^{\beta}_{\text{char}}} \right),
\end{equation}
where $\gamma_{\text{char}}$ is the characteristic shear due to mechanical twinning and depends on the twin system. Interested readers are referred to the work of Roters et al.~\cite{roters2010overview} for a complete picture of CPFEM model in general and for using \texttt{DAMASK}~\cite{roters2019damask} in particular. 
% For spectral solver implementation, readers are referred to Eisenlohr et al.~\cite{eisenlohr2013spectral} and Shanthraj et al.~\cite{shanthraj2015numerically,shanthraj2019spectral}.
Table~\ref{tab:MgConstitutiveParameters} lists the constitutive parameter used in this example from the literature ~\cite{sedighiani2020efficient,sedighiani2022determination,wang2014situ,tromans2011elastic,agnew2006validating}. 

\begin{table}[!htbp]
% \tiny
\caption{Parameters for Mg used in this case study (cf. Tables 7 and 8~\cite{sedighiani2020efficient,sedighiani2022determination}, ~\cite{wang2014situ,tromans2011elastic,agnew2006validating}).}
\label{tab:MgConstitutiveParameters}
\begin{tabular*}{\textwidth}{c @{\extracolsep{\fill}} cccc} \hline
variable                    & description                       & units     &  reference value  \\ \hline
$c/a$                       & lattice parameter ratio                 &   --    &  1.635        \\   
$C_{11}$                    & elastic constant                    &   GPa     &   59.3        \\   
$C_{12}$                    & elastic constant                    &   GPa     &   61.5        \\
$C_{44}$                    & elastic constant                    &   GPa     &   16.4        \\
$C_{44}$                    & elastic constant                    &   GPa     &   25.7        \\
$C_{44}$                    & elastic constant                    &   GPa     &   21.4        \\ 
$\dot{\gamma}_0$                & twinning reference shear rate             &  s$^{-1}$   &  0.001        \\
$\dot{\gamma}_0$                & slip reference shear rate               &  s$^{-1}$   &  0.001        \\ \hline
$\tau_{0,\text{basal}}$             & basal slip resistance                 &   MPa     &  10         \\
$\tau_{0,\text{pris}}$              & prismatic slip resistance               &   MPa     &  55         \\
$\tau_{0,\text{pyr} \langle a \rangle}$     & pyramidal $\langle a \rangle$ slip resistance     &   MPa     &  60         \\
$\tau_{0,\text{pyr} \langle c+a \rangle}$     & pyramidal $\langle c+a \rangle$ slip resistance     &   MPa     &  60         \\
$\tau_{0,\text{T}1}$              & tensile twin resistance                 &   MPa     &  45         \\
$\tau_{0,\text{C}2}$              & compressive twin resistance               &   MPa     &  80         \\
$\tau_{\infty,\text{basal}}$          & basal saturation stress                 &   MPa     &  45         \\
$\tau_{\infty,\text{pris}}$           & prismatic saturation stress               &   MPa     &  135        \\
$\tau_{\infty,\text{pyr} \langle a \rangle}$  & pyramidal $\langle a \rangle $ saturation stress    &   MPa     &  150        \\
$\tau_{\infty,\text{pyr} \langle c+a \rangle}$  & pyramidal $\langle c+a \rangle$ saturation stress   &   MPa     &  150        \\
$h_{0}^{\text{tw}-\text{tw}}$          & twin-twin hardening parameter             &   MPa     &  50          \\
$h_{0}^{\text{s}-\text{s}}$          & slip-slip hardening parameter             &   MPa     &  500         \\
$h_{0}^{\text{tw}-\text{s}}$           & twin-slip hardening parameter             &   MPa     &  150         \\
$n_\text{s}$                  & slip strain rate sensitivity parameter        &   --    &  10         \\
$n_\text{tw}$                   & twinning strain rate sensitivity parameter      &   --    &  5          \\
$a$                       & slip hardening parameter                &   --    &  2.5        \\ \hline
\end{tabular*}
\normalsize
\end{table}

\begin{table}[]
    \centering
    \caption{Mesh resolutions and corresponding number of degrees of freedom (i.e. SERVE size) used in this study.}
    \begin{tabular}{cccc}\toprule
         $\ell$ & Degrees of freedom & \# processors/sample & cost/sample [s]\\ \midrule
         0 &  $2 \times 2 \times 2$    & 1 & 39    \\
         1 &  $4 \times 4 \times 4$    & 1 & 365   \\
         2 &  $8 \times 8 \times 8$    & 2 & 1955  \\
         3 &  $16 \times 16 \times 16$ & 4 & 3305  \\
         4 &  $32 \times 32 \times 32$ & 8 & 12487 \\\bottomrule
    \end{tabular}
    \label{tab:dofs}
\end{table}

\subsection{Crystallography and microstructure for magnesium}

% generated from $64\times 64 \times 64$ original SERVE, followed by down scaling

% \hl{Reminder: FEM approximates both geometry and numerical solutions. Two levels of approximations required.}

\begin{figure}[!htbp]
\centering

\begin{subfigure}[b]{0.175\textwidth}
\includegraphics[width=\textwidth]{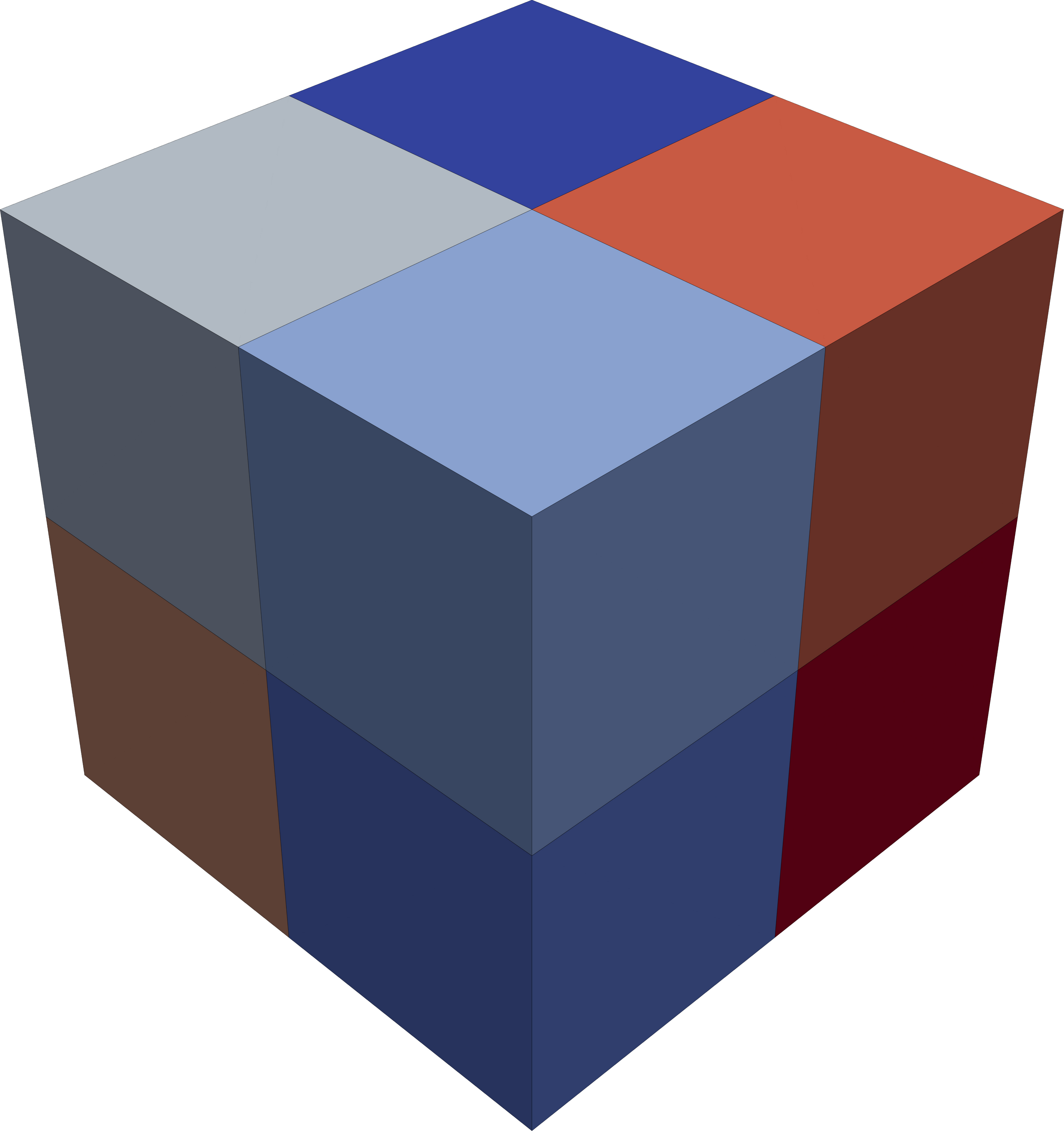}
\caption{$m_1: 2 \times 2 \times 2$}
\end{subfigure}
\begin{subfigure}[b]{0.175\textwidth}
\includegraphics[width=\textwidth]{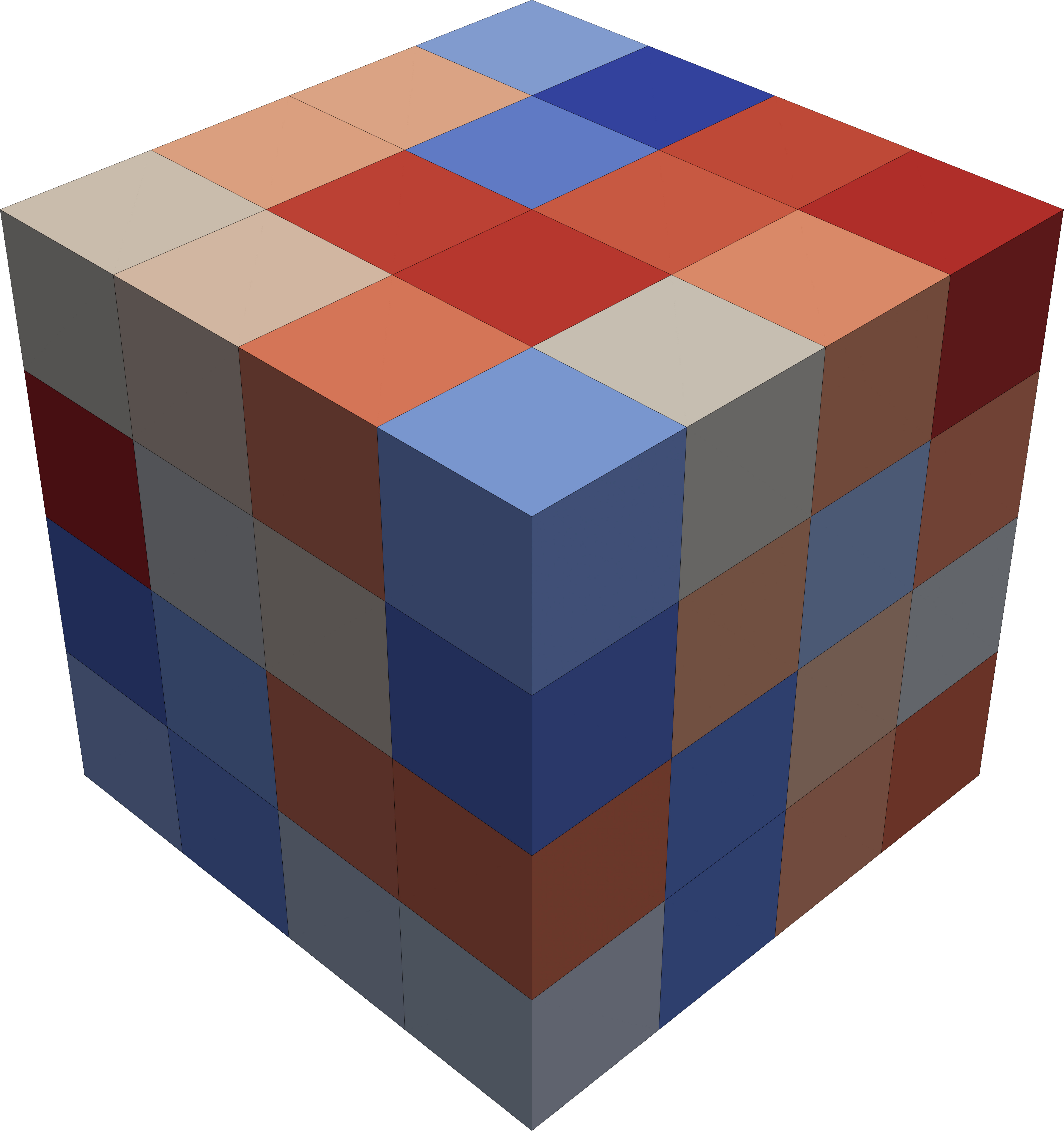}
\caption{$m_1: 4 \times 4 \times 4$}
\end{subfigure}
\begin{subfigure}[b]{0.175\textwidth}
\includegraphics[width=\textwidth]{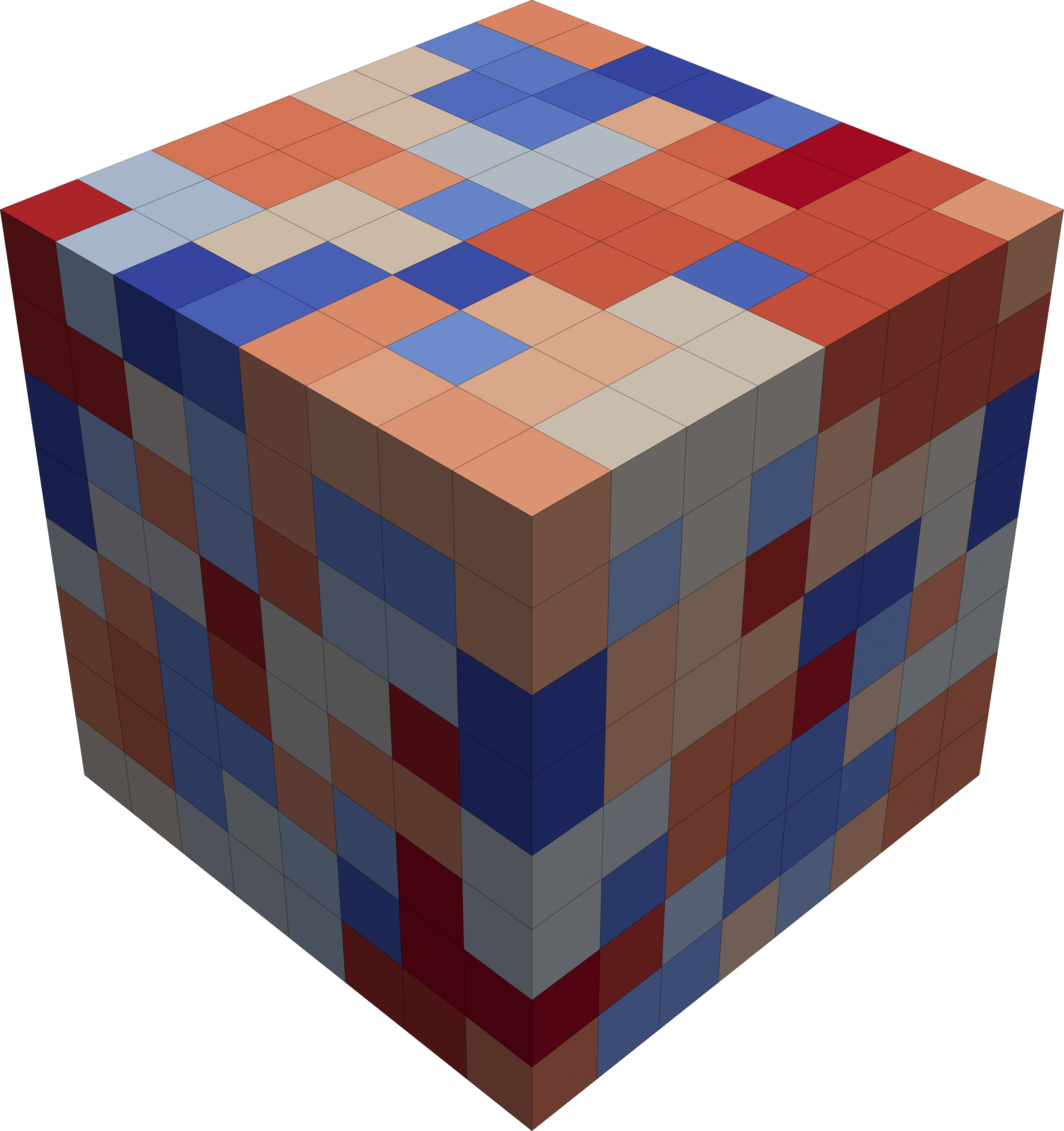}
\caption{$m_1: 8 \times 8 \times 8$}
\end{subfigure}
\begin{subfigure}[b]{0.175\textwidth}
\includegraphics[width=\textwidth]{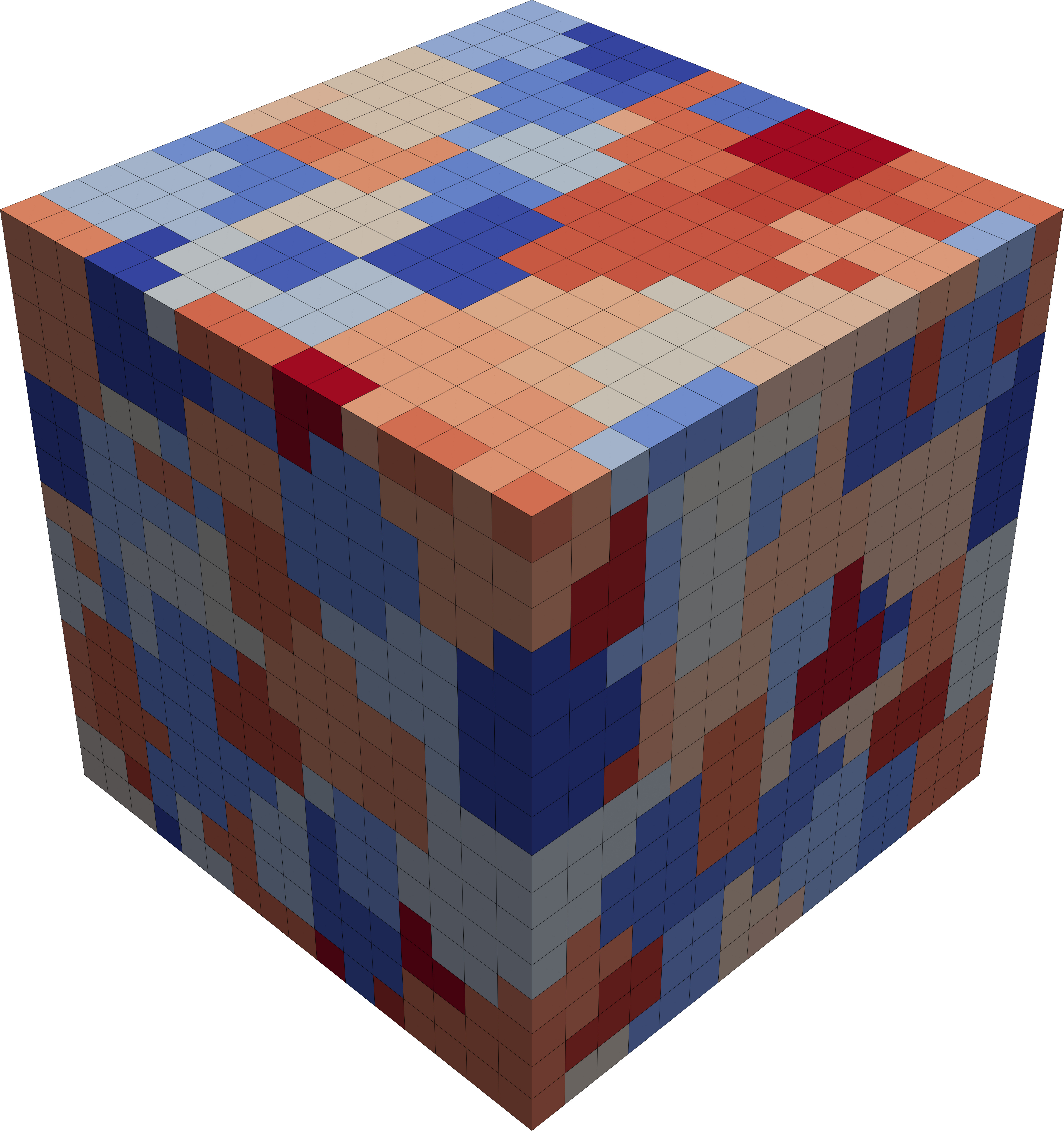}
\caption{$m_1: 16 \times 16 \times 16$}
\end{subfigure}
\begin{subfigure}[b]{0.175\textwidth}
\includegraphics[width=\textwidth]{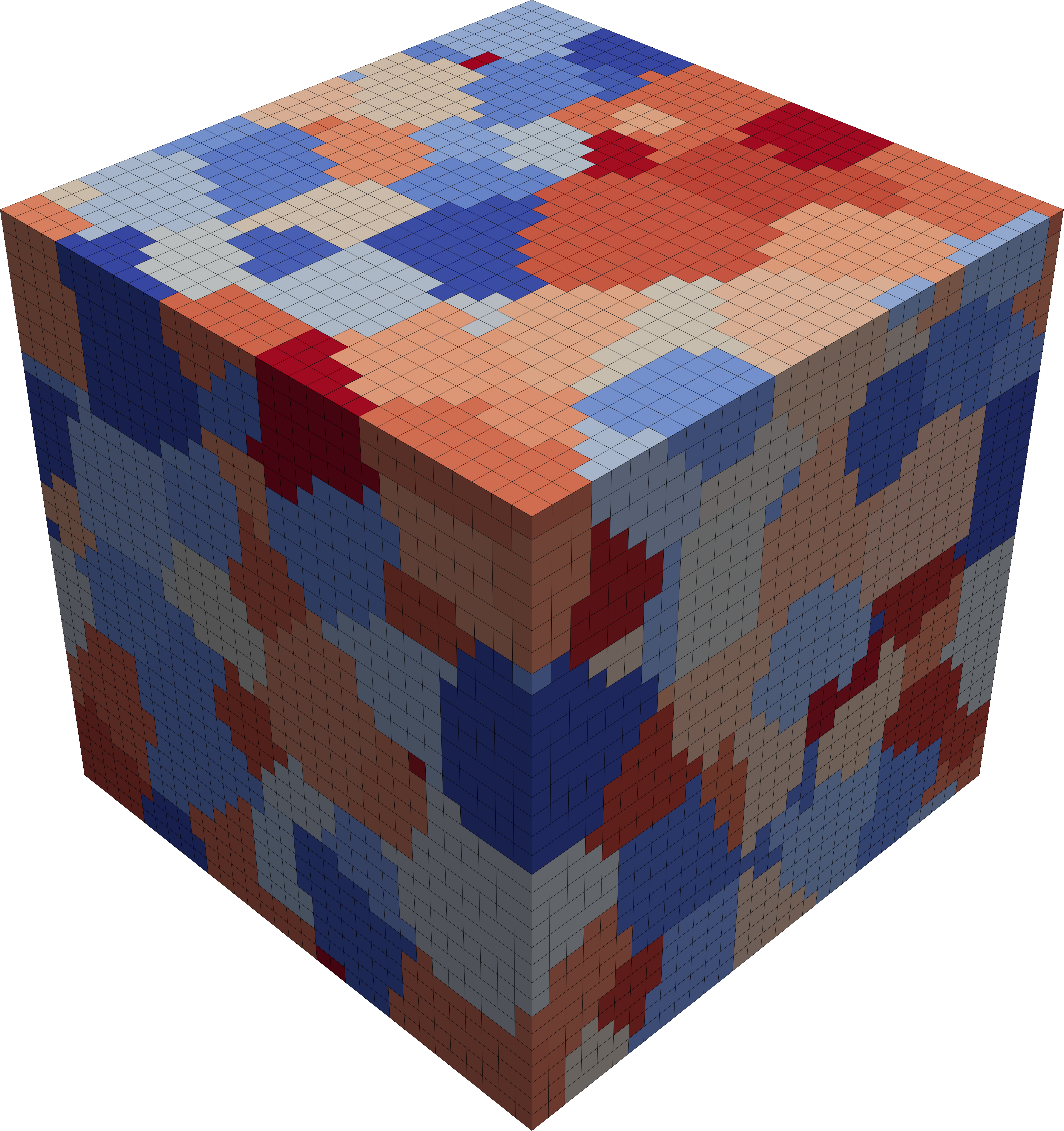}
\caption{$m_1: 32 \times 32 \times 32$}
\end{subfigure}

\begin{subfigure}[b]{0.175\textwidth}
\includegraphics[width=\textwidth]{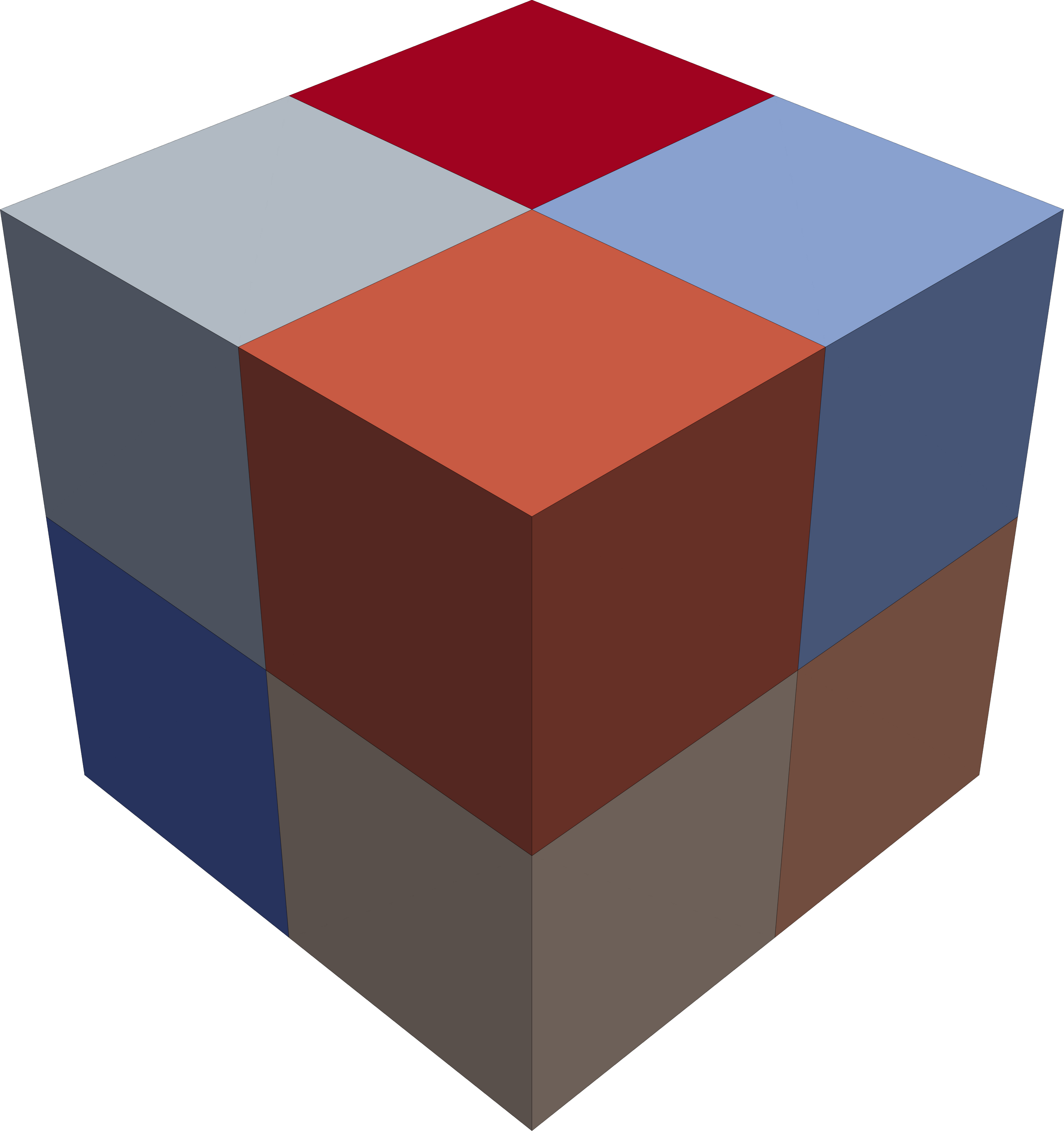}
\caption{$m_2: 2 \times 2 \times 2$}
\end{subfigure}
\begin{subfigure}[b]{0.175\textwidth}
\includegraphics[width=\textwidth]{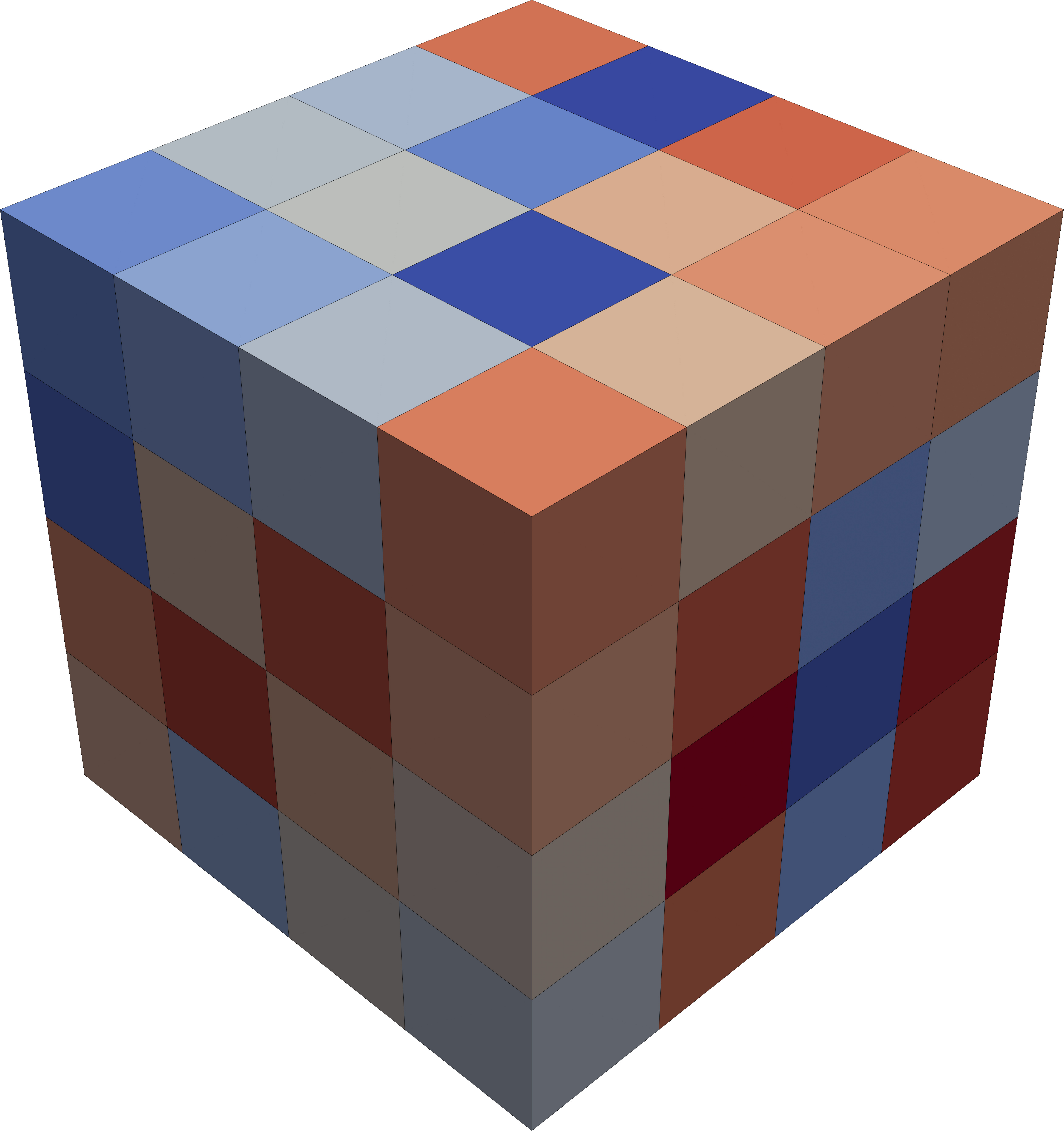}
\caption{$m_2: 4 \times 4 \times 4$}
\end{subfigure}
\begin{subfigure}[b]{0.175\textwidth}
\includegraphics[width=\textwidth]{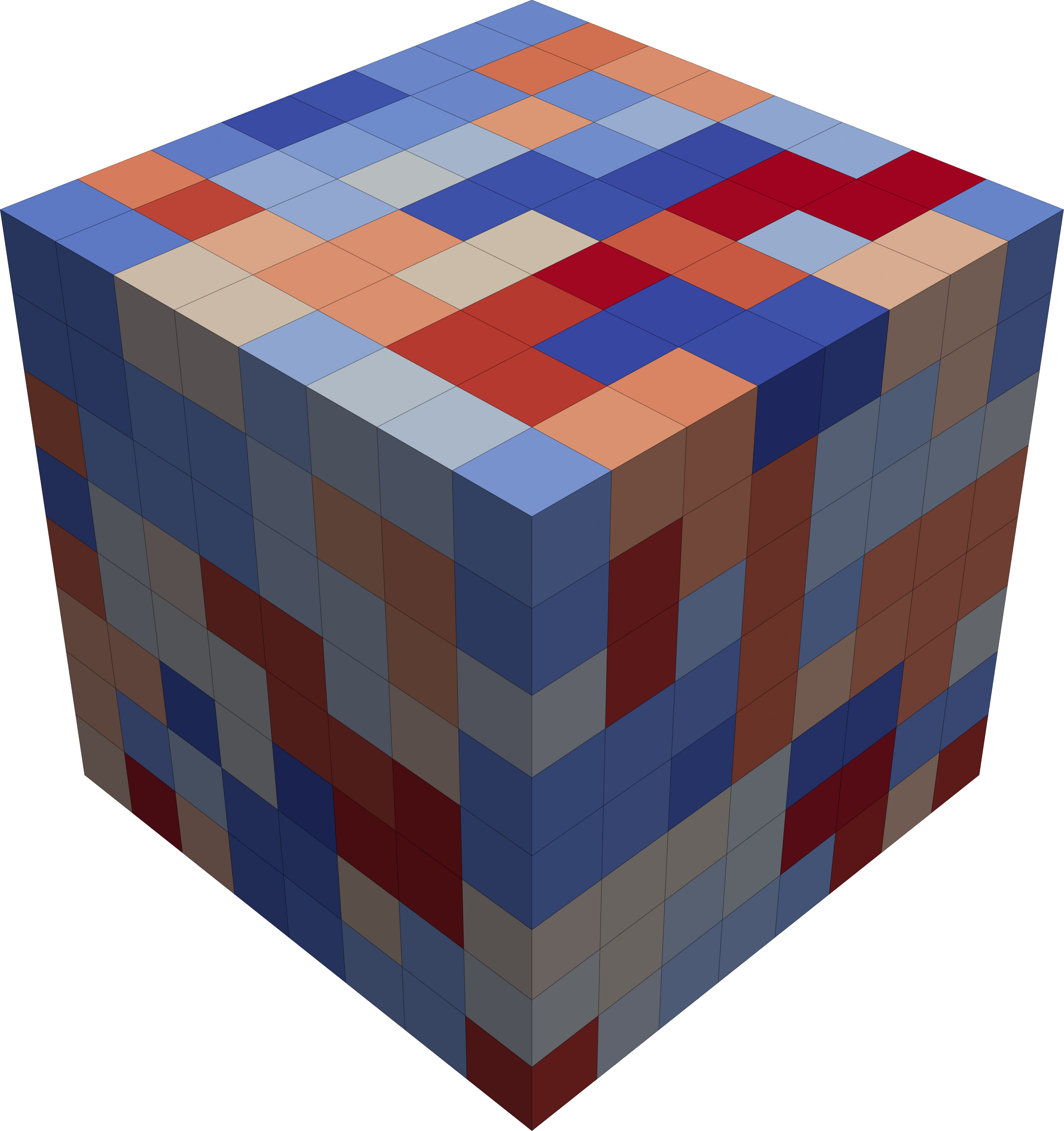}
\caption{$m_2: 8 \times 8 \times 8$}
\end{subfigure}
\begin{subfigure}[b]{0.175\textwidth}
\includegraphics[width=\textwidth]{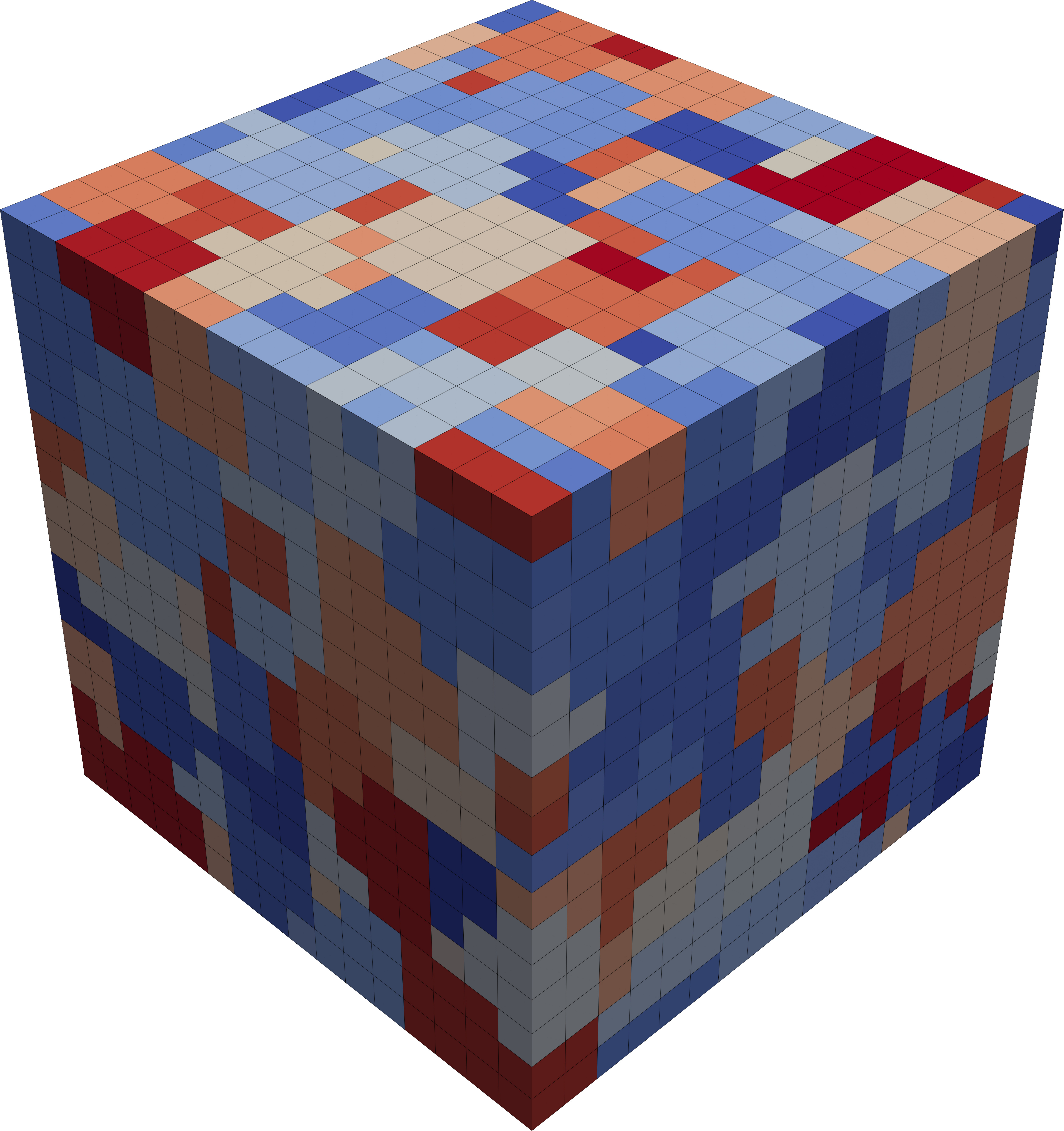}
\caption{$m_2: 16 \times 16 \times 16$}
\end{subfigure}
\begin{subfigure}[b]{0.175\textwidth}
\includegraphics[width=\textwidth]{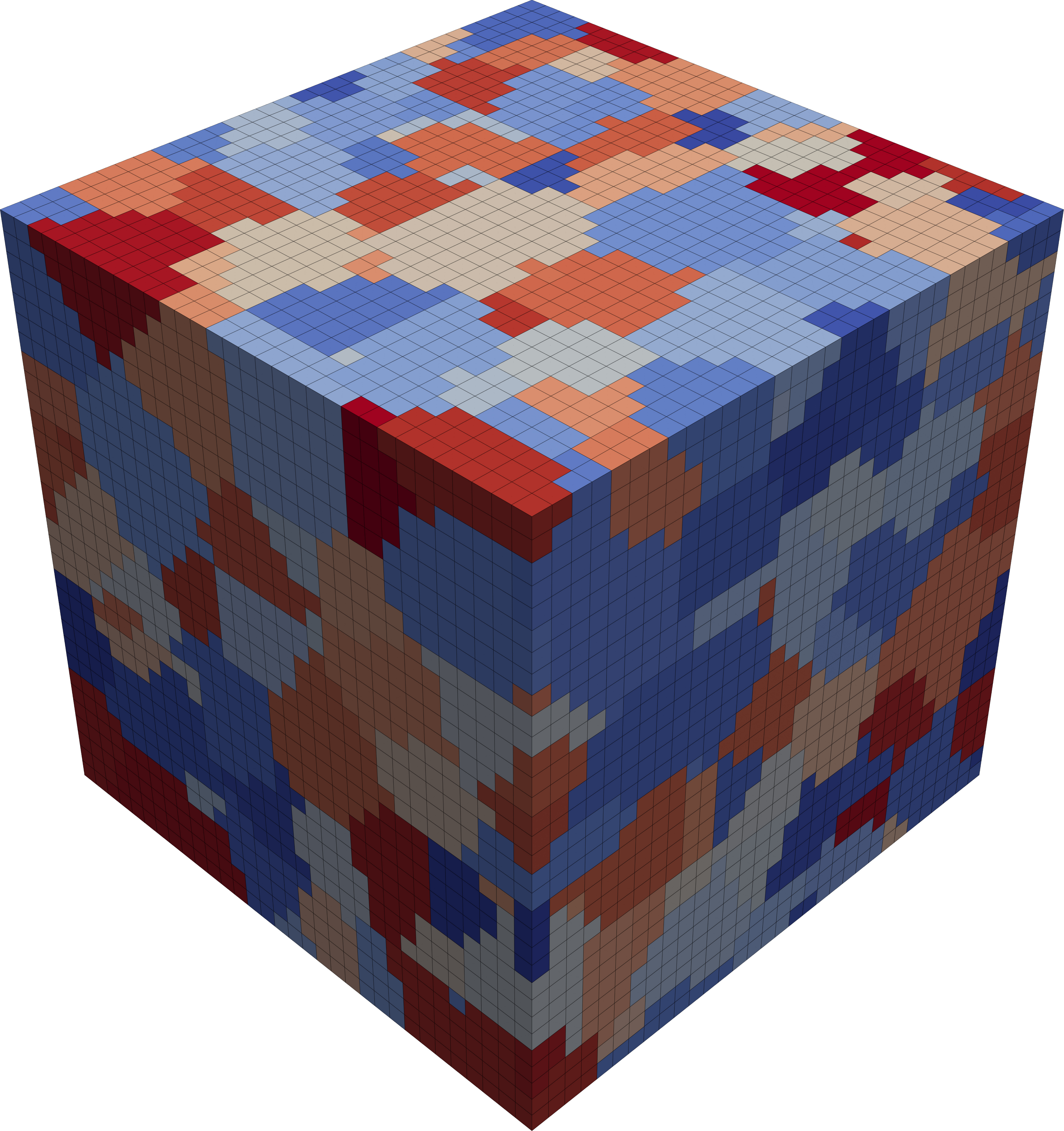}
\caption{$m_2: 32 \times 32 \times 32$}
\end{subfigure}

\begin{subfigure}[b]{0.175\textwidth}
\includegraphics[width=\textwidth]{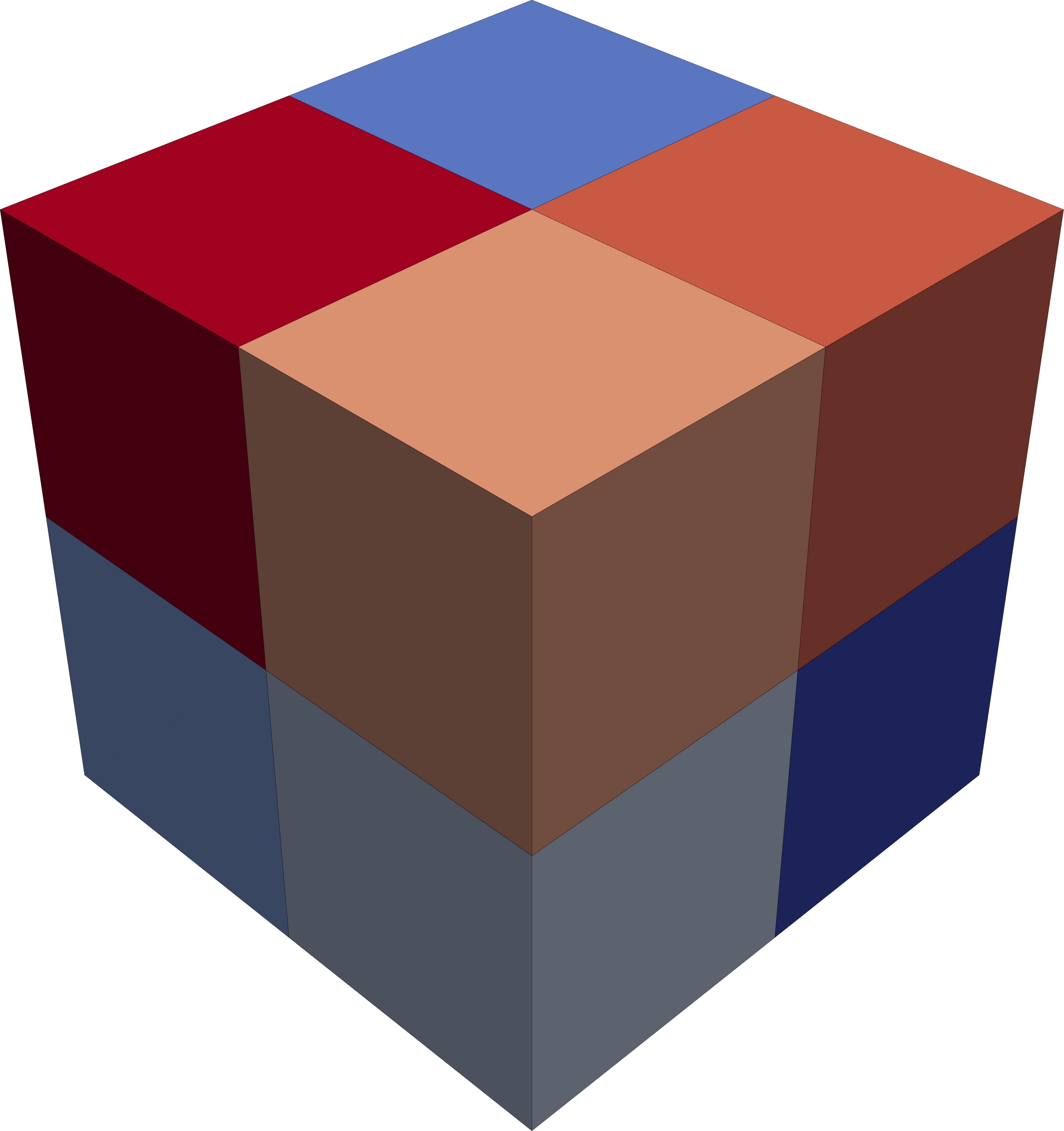}
\caption{$m_3: 2 \times 2 \times 2$}
\end{subfigure}
\begin{subfigure}[b]{0.175\textwidth}
\includegraphics[width=\textwidth]{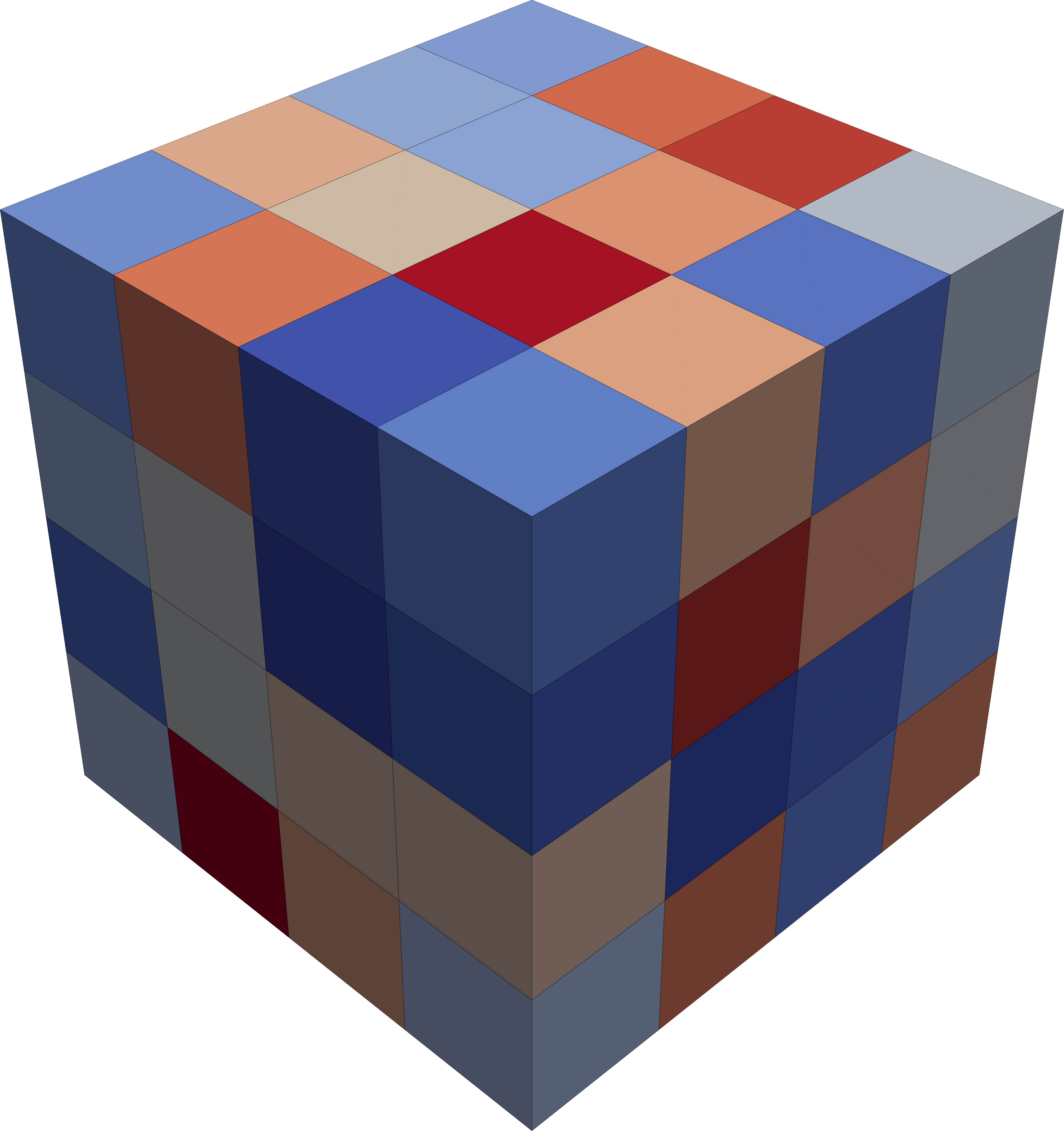}
\caption{$m_3: 4 \times 4 \times 4$}
\end{subfigure}
\begin{subfigure}[b]{0.175\textwidth}
\includegraphics[width=\textwidth]{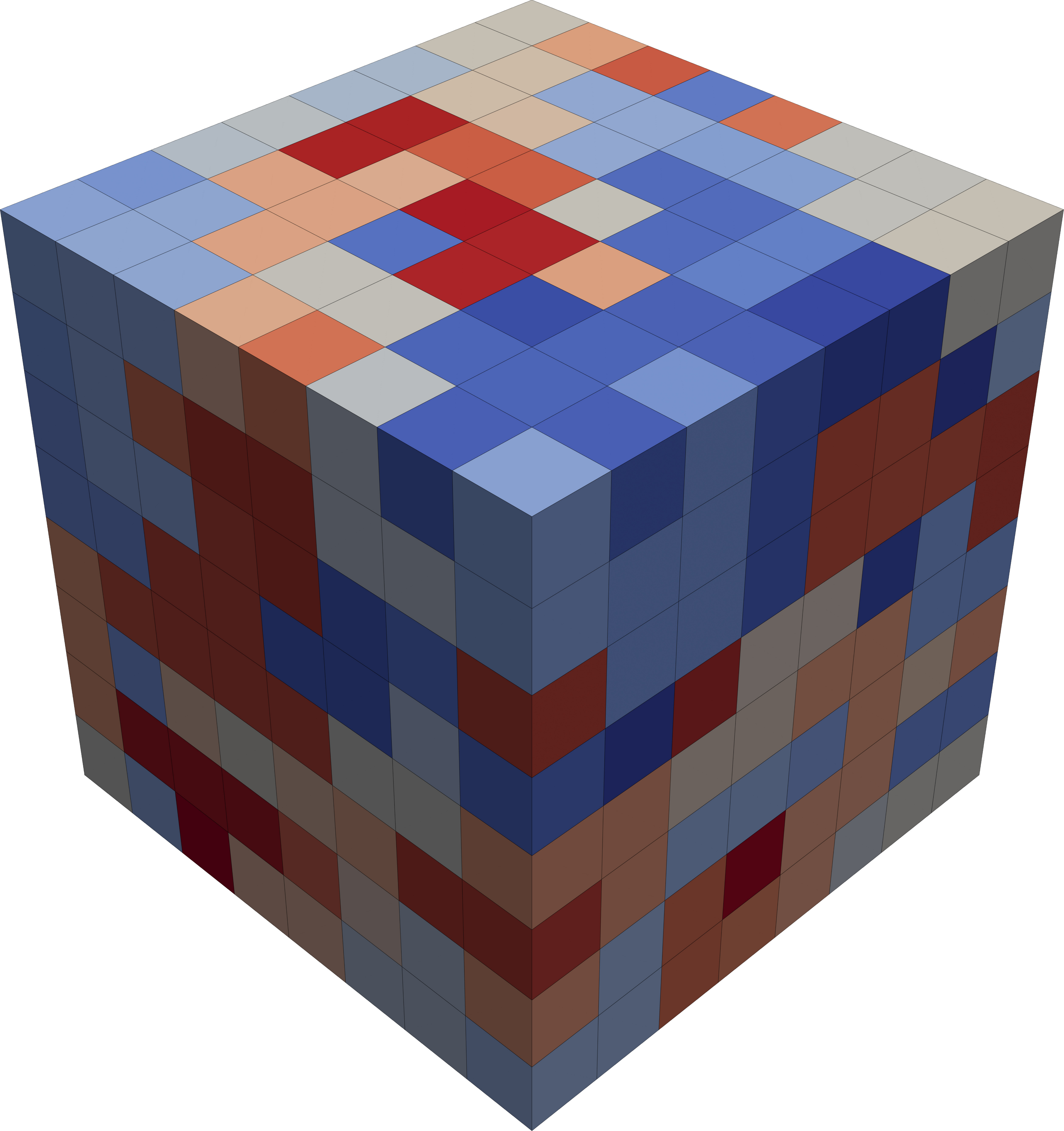}
\caption{$m_3: 8 \times 8 \times 8$}
\end{subfigure}
\begin{subfigure}[b]{0.175\textwidth}
\includegraphics[width=\textwidth]{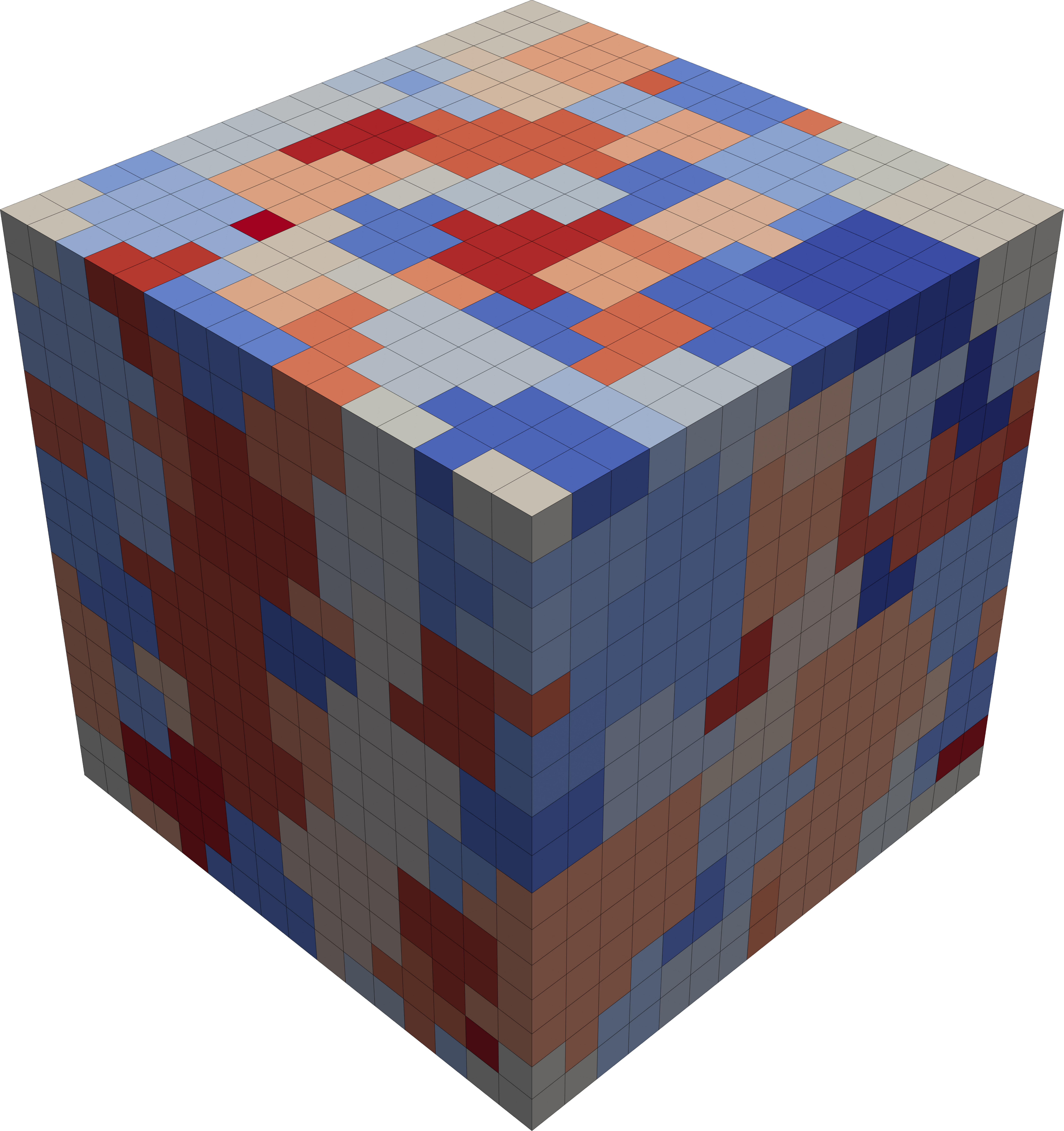}
\caption{$m_3: 16 \times 16 \times 16$}
\end{subfigure}
\begin{subfigure}[b]{0.175\textwidth}
\includegraphics[width=\textwidth]{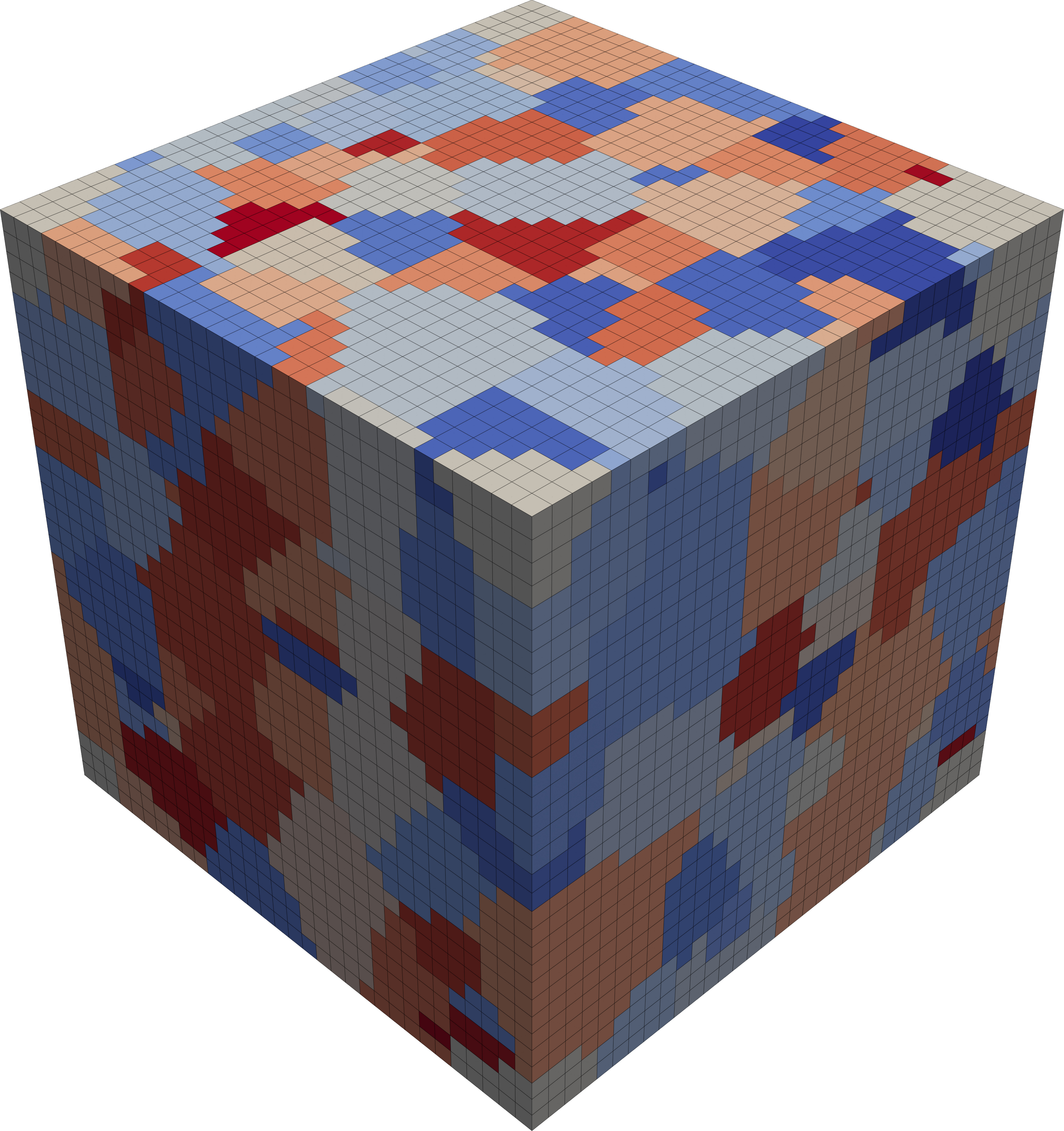}
\caption{$m_3: 32 \times 32 \times 32$}
\end{subfigure}

\caption{Three SERVE realizations, i.e. $m_1$, $m_2$, and $m_3$, discretized at 5 geometric mesh resolutions: $2 \times 2 \times 2$, $4 \times 4 \times 4$, $8 \times 8 \times 8$, $16 \times 16 \times 16$, and $32 \times 32 \times 32$.}
\label{fig:mesh_refinements}
\end{figure}

\begin{figure}
\includegraphics[width=\textwidth,keepaspectratio]{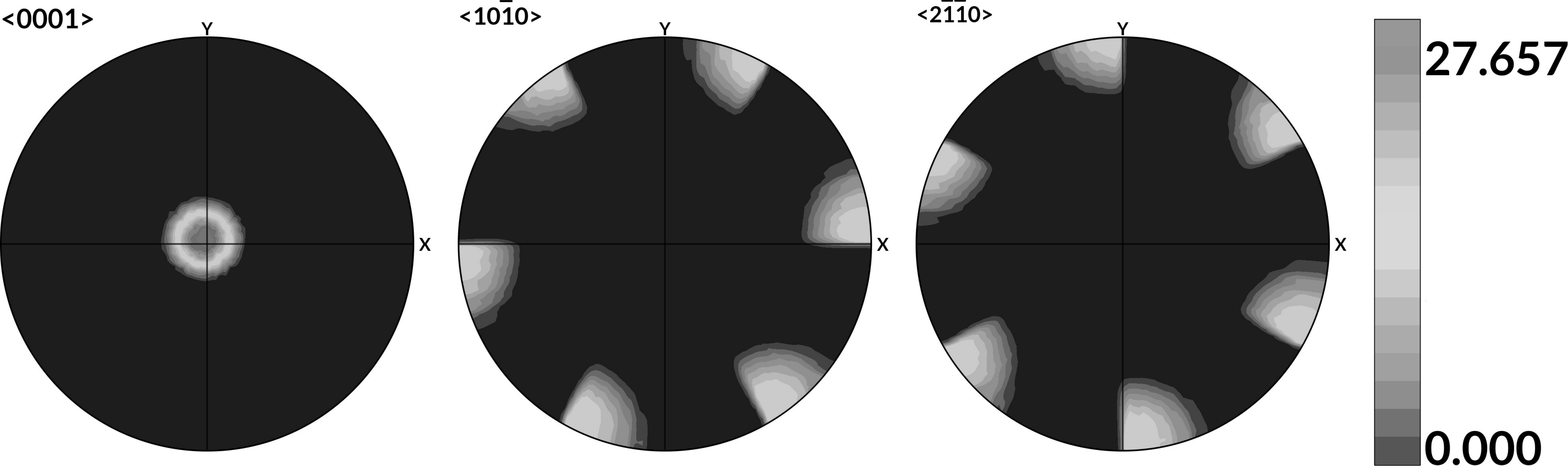}
\caption{Magnesium texture component with Euler angles $(\phi_1, \theta, \phi_2) = (90^\circ,0^\circ,0^\circ)$}.
\label{fig:Mg_odf_texture}
\end{figure}

Figure~\ref{fig:mesh_refinements} presents three different SERVEs with geometric mesh resolution levels chosen as detailed in Table~\ref{tab:dofs}. The SERVEs are reconstructed using \texttt{DREAM.3D}, where the grain sizes are log normally distributed with equiaxed grains, i.e. $d \sim \text{LogNormal}(\mu_d, \sigma_d)$, $\mu_d = 5.2983, \sigma_d = 0.2$. 
Following~\cite{mangal2018dataset,tran2022microstructure}, the crystallographic texture for magnesium is sampled from $(\phi_1, \theta, \phi_2) = (90^\circ,0^\circ,0^\circ)$ and shown in Figure~\ref{fig:Mg_odf_texture}.

\section{Numerical results}\label{sec:results}

\subsection{Problem description}

\begin{figure}
\centering
\includegraphics[width=0.5\textwidth]{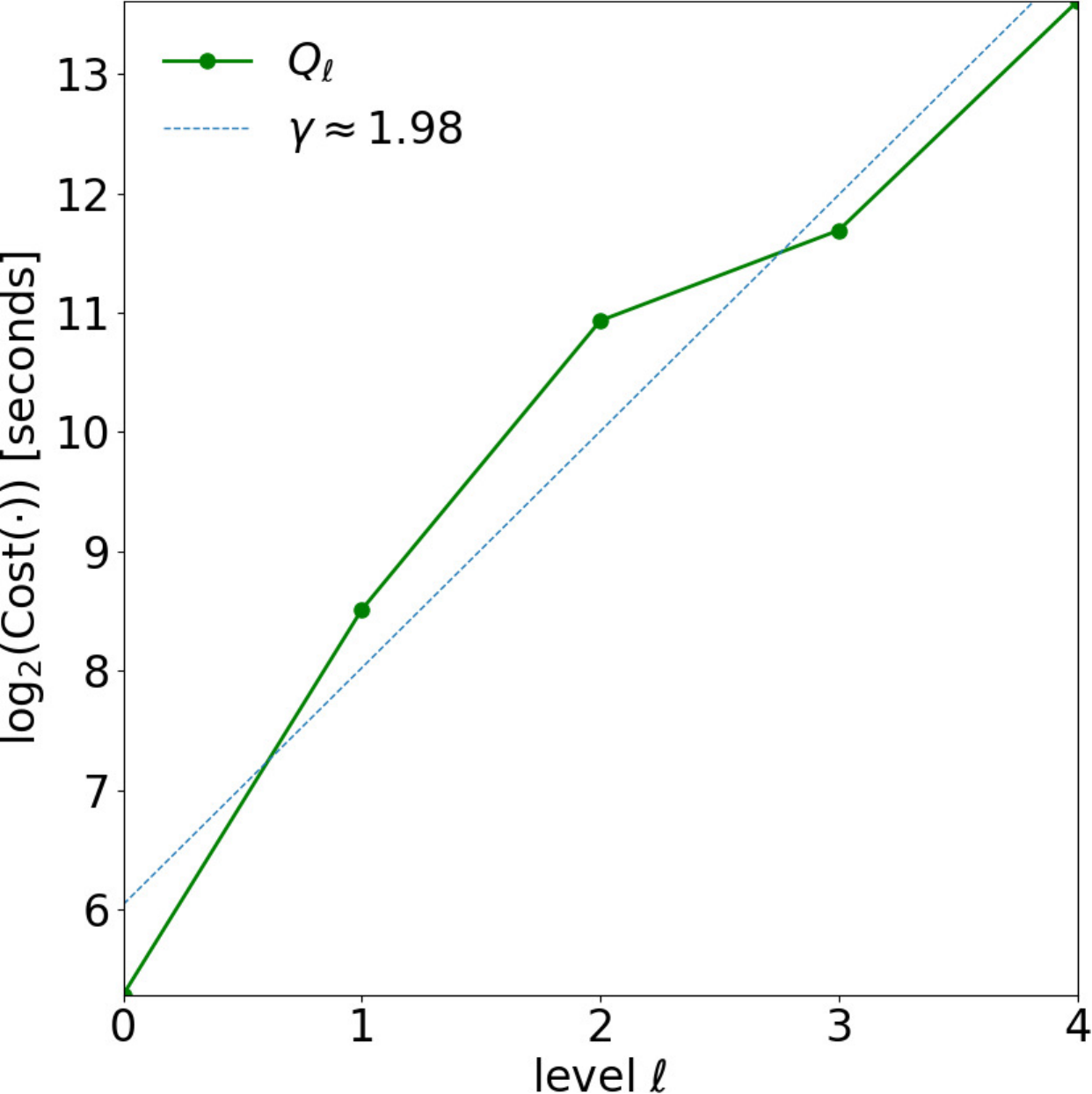}
\caption{Increase of the average computational cost per sample expressed in seconds as a function of the level parameter $\ell$ using the mesh refinement as a fidelity parameter in the MLMC experiment. The dashed line corresponds to a fit of the computational cost proportional to $2^{\gamma\ell}$ with $\gamma\approx1.98$, see condition \eqref{eq:C3}. The computational cost at each level can be found in Table~\ref{tab:dofs}).}
\label{fig:cost_per_sample_vs_level}
\end{figure}

In this section, we discuss the setup of our numerical experiments. We are interested in estimating homogenized stress curves for magnesium at different strain levels. 
We compute those curves using both the MC method and the MLMC method outlined in Section~\ref{sec:MLMC}. For the MLMC method, we use a geometric mesh resolution level hierarchy as shown in Table~\ref{tab:dofs}. 
Figure~\ref{fig:cost_per_sample_vs_level} shows the increase in computational cost associated with each level, indicating that the cost increases geometrically with the mesh resolution level number. 
Although the number of degrees of freedom increases with a factor 8 as the mesh resolution level increases, the computational cost scales only as $2^{\gamma\ell}$, where we numerically fitted $\gamma\approx1.98$. This is because we use a different number of processors at each level, following a constraint on the number of grid points per processor imposed by \texttt{DAMASK}. %In particular, for mesh resolution level 4, using $32 \times 32 \times 32$ DOFs, we use a total of 8 processors.
The number of processors used for each mesh resolution level is shown in Table~\ref{tab:dofs}. 
% \pjr{\sout{the MPI parallelization implemented within \texttt{DAMASK} imposes a constraint on the number of processors being used. For example, $32 \times 32 \times 32$ would only use a maximum of 8 cores. To exploit the MPI parallelism and reduce the wall-clock waiting time, we also double the computational resource $2 \times$, by doubling the number of processors used with the same number of nodes, due to the fact that most of the computational resources in HPC nodes are under-utilized. This mechanism results in about $4 \times$ increase in wall-clock waiting time, confirmed in Figure~\ref{fig:cost_per_sample_vs_level} with $\gamma \approx 1.98$, which is very close to 2 (i.e. $2^\gamma \approx 4$). Had the number of processors not double, $\gamma$ would have been closer to $3$.}}

\subsection{MLMC estimation of homogenized mean stress-strain responses}
\label{subsec:MLMC_Applications}

\begin{figure}[!htbp]
\centering
\includegraphics[width=1.0\textwidth]{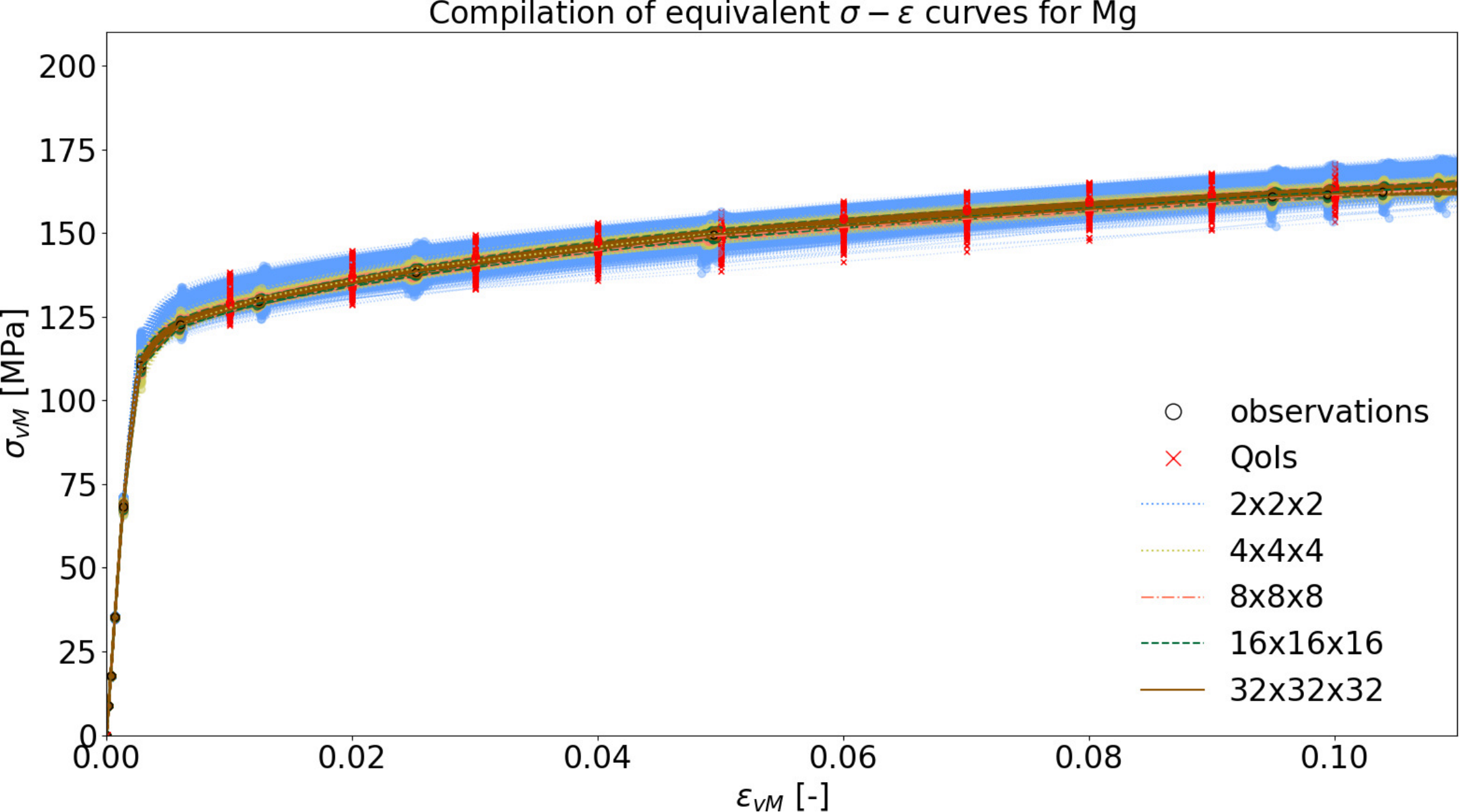}
\caption{Compilation of stress-strain curves by sampling SERVEs at the different mesh resolution levels shown in Table~\ref{tab:dofs}. The locations of the 9 QoIs are plotted as \large$\textcolor{red}{\times}$\normalsize.}
\label{fig:stress_strain_compilation_magnesium}
\end{figure}

Figure~\ref{fig:stress_strain_compilation_magnesium} shows a collection of the stress-strain curve evaluated for an ensemble of SERVEs at each mesh resolution level provided in Table~\ref{tab:dofs}. 
As the fidelity level increases, the mesh becomes finer, and the number of SERVEs decreases. 
Hence, the plot is dominated by stress-strain curves collected at $2 \times 2 \times 2$, and much less at $32 \times 32 \times 32$. 
We interpolate the obtained stress-strain curves at a set of 9 prescribed strain values $\varepsilon = \{0.1, 0.2, \ldots, 0.9\}$. We use cubic Hermite spline interpolation~\cite{virtanen2020scipy} to obtain the values of the stress-strain curves at those 9 locations. This higher-order interpolation scheme should avoid additional interpolation errors that converge at a slower rate than the bias in the predicted stress-strain curves. The value of the stress-strain curve at those 9 locations will be the set of QoIs we are interested in. Note that these QoIs are indicated by \huge$\textcolor{red}{\times}$\normalsize in Figure~\ref{fig:stress_strain_compilation_magnesium}. 
% \pjr{\sout{presents the stress-strain compilation for many SERVEs at 5 different mesh resolutions. The observations of stresses from CPFEM are plotted as \huge$\circ$\normalsize, whereas the interpolated stresses are plotted as \huge$\textcolor{red}{\times}$\normalsize. The interpolated stresses at collocated strains, which are our QoIs, are obtained via the piecewise cubic Hermite interpolating polynomial interpolator}~\cite{fritsch1984method}, \sout{also known as PCHIP interpolator, implemented in \texttt{scipy}}~\cite{virtanen2020scipy}. \sout{Up to 1000 homogenized stress-strain curves are plotted where different colors correspond to different mesh resolutions. Readers are referred to color version online.} }

\begin{figure}[!htbp]
\centering

\begin{subfigure}[b]{0.30\textwidth}
\includegraphics[width=\textwidth]{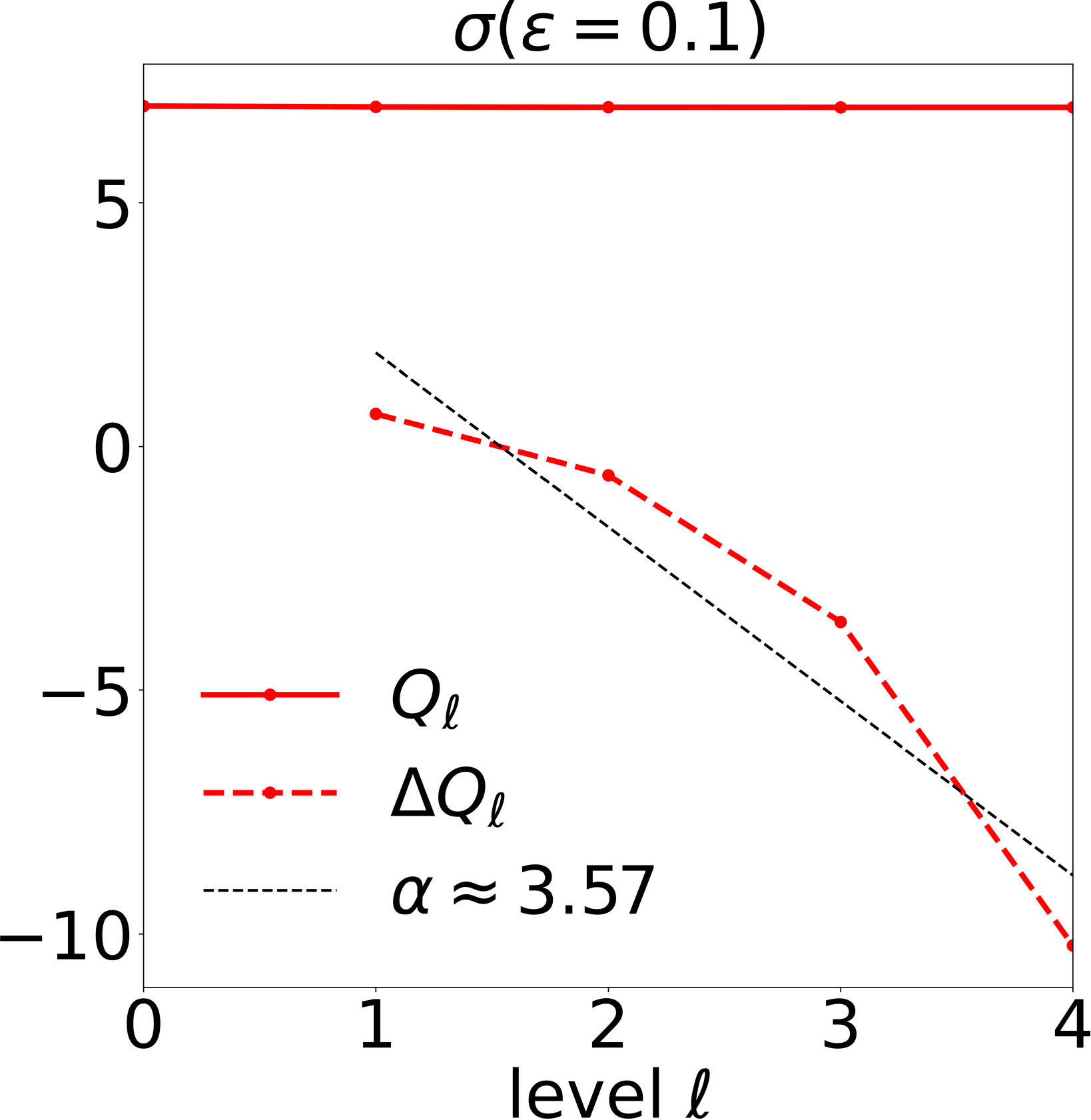}
\caption{$\varepsilon = 0.10$}
\end{subfigure}
\begin{subfigure}[b]{0.30\textwidth}
\includegraphics[width=\textwidth]{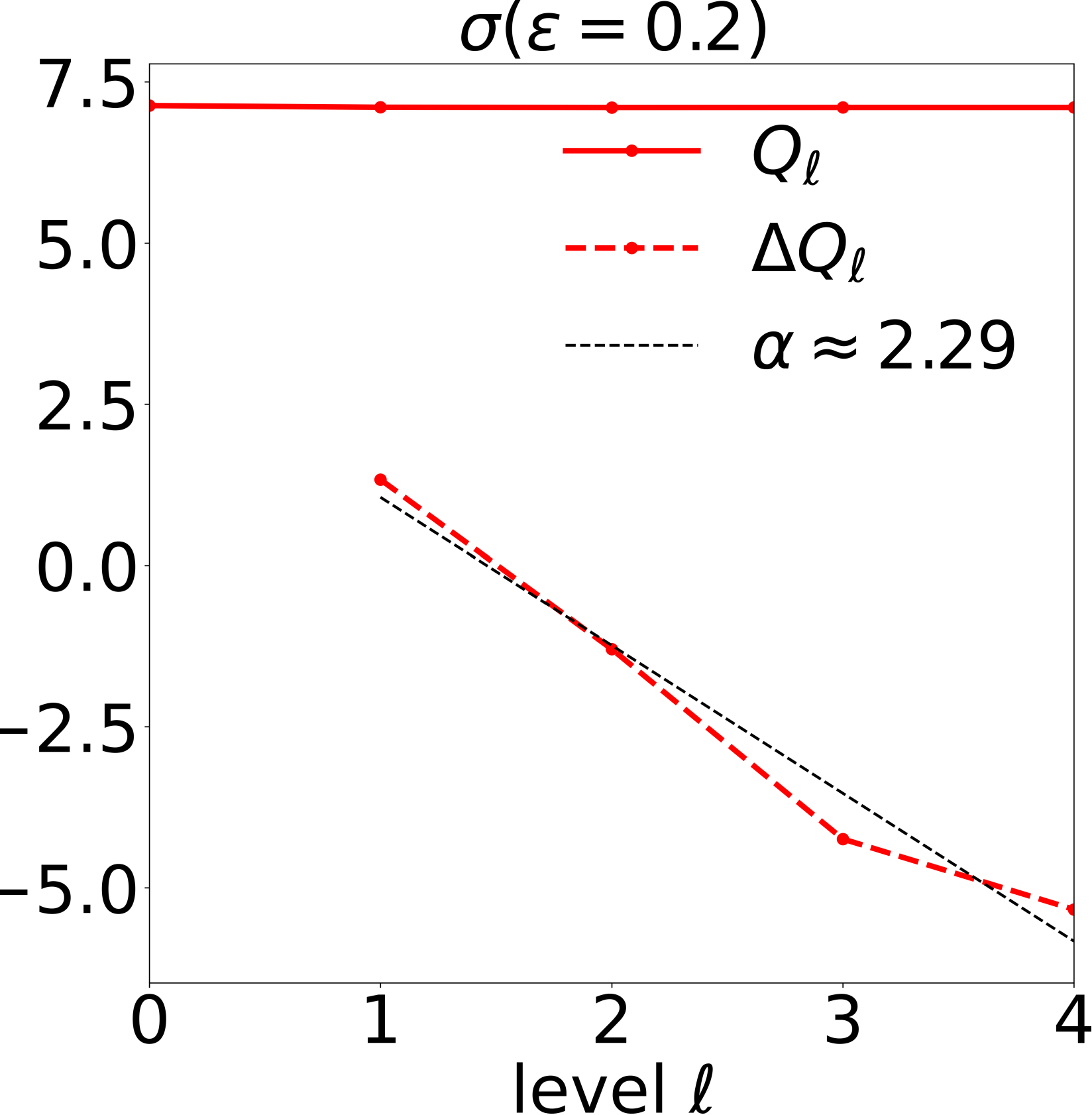}
\caption{$\varepsilon = 0.20$}
\end{subfigure}
\begin{subfigure}[b]{0.30\textwidth}
\includegraphics[width=\textwidth]{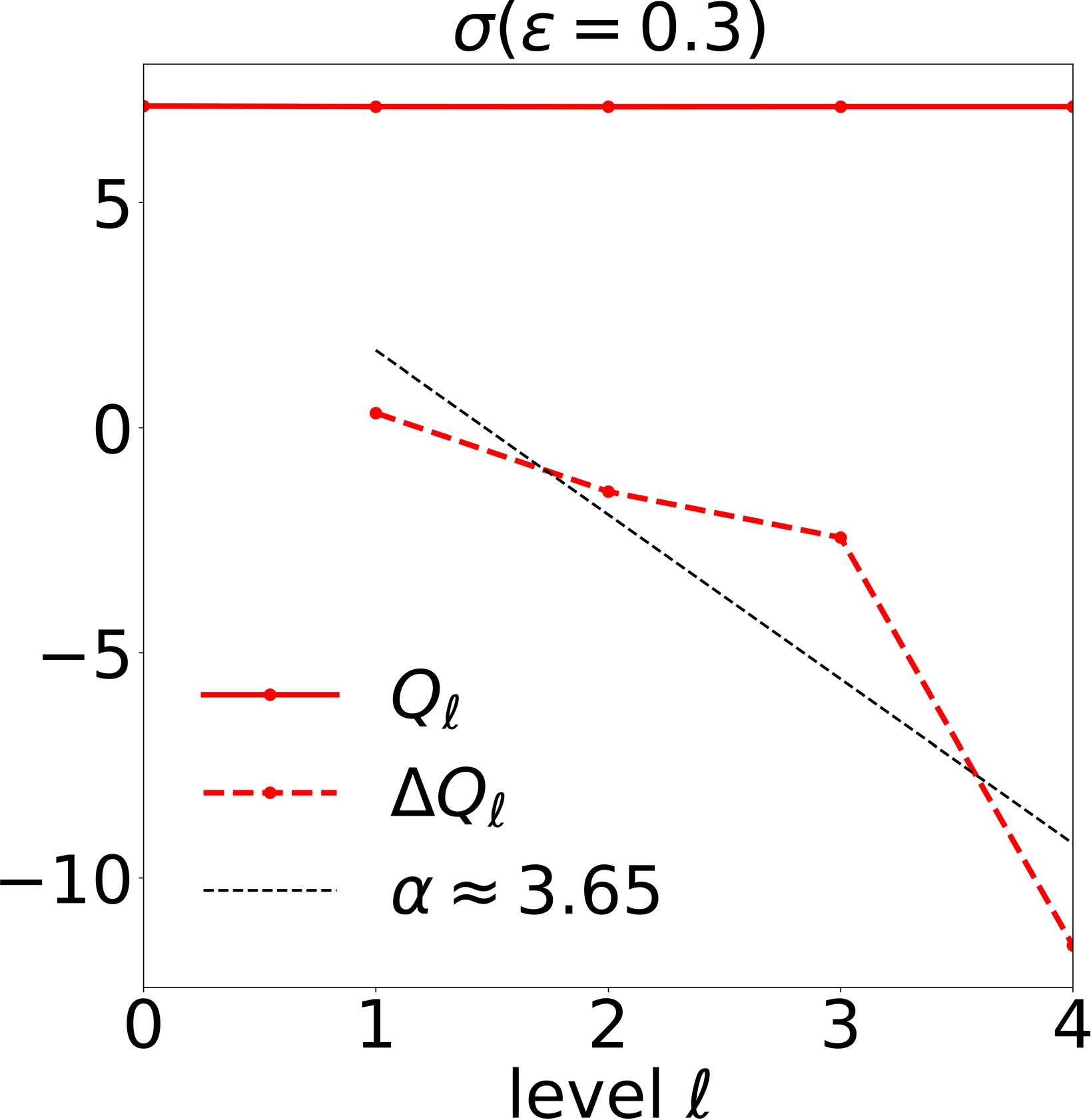}
\caption{$\varepsilon = 0.30$}
\end{subfigure}

\begin{subfigure}[b]{0.30\textwidth}
\includegraphics[width=\textwidth]{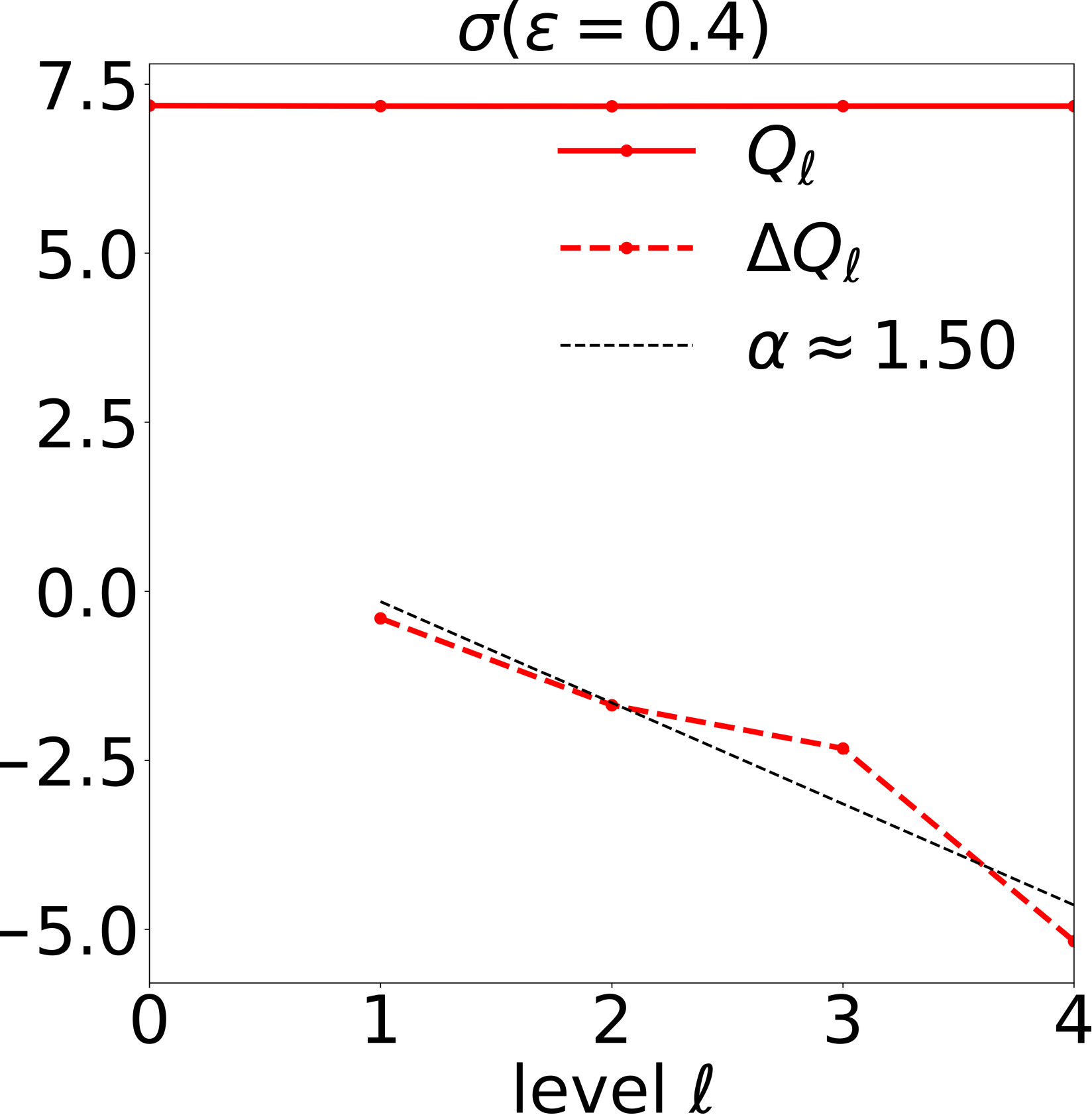}
\caption{$\varepsilon = 0.40$}
\end{subfigure}
\begin{subfigure}[b]{0.30\textwidth}
\includegraphics[width=\textwidth]{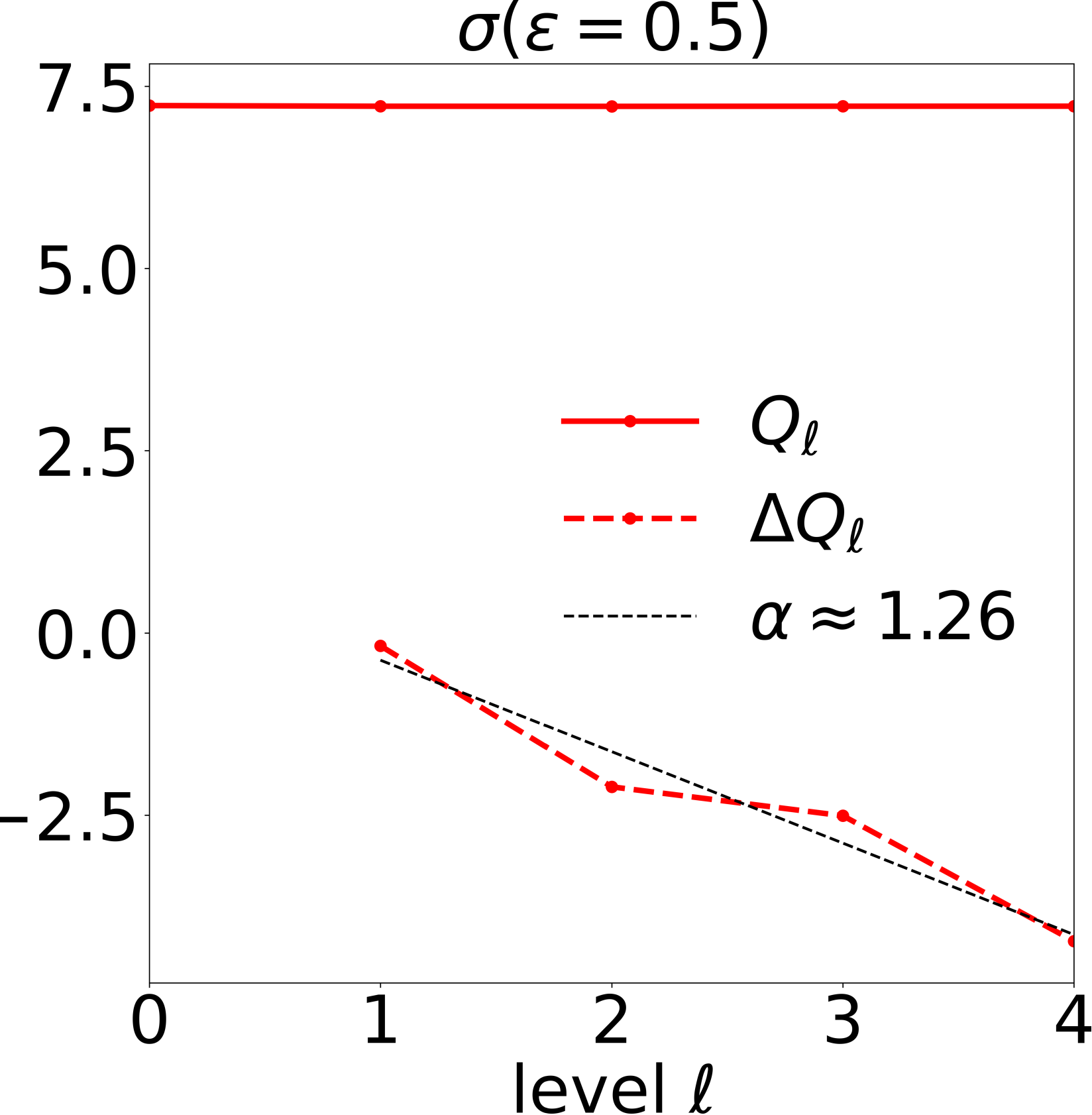}
\caption{$\varepsilon = 0.50$}
\end{subfigure}
\begin{subfigure}[b]{0.30\textwidth}
\includegraphics[width=\textwidth]{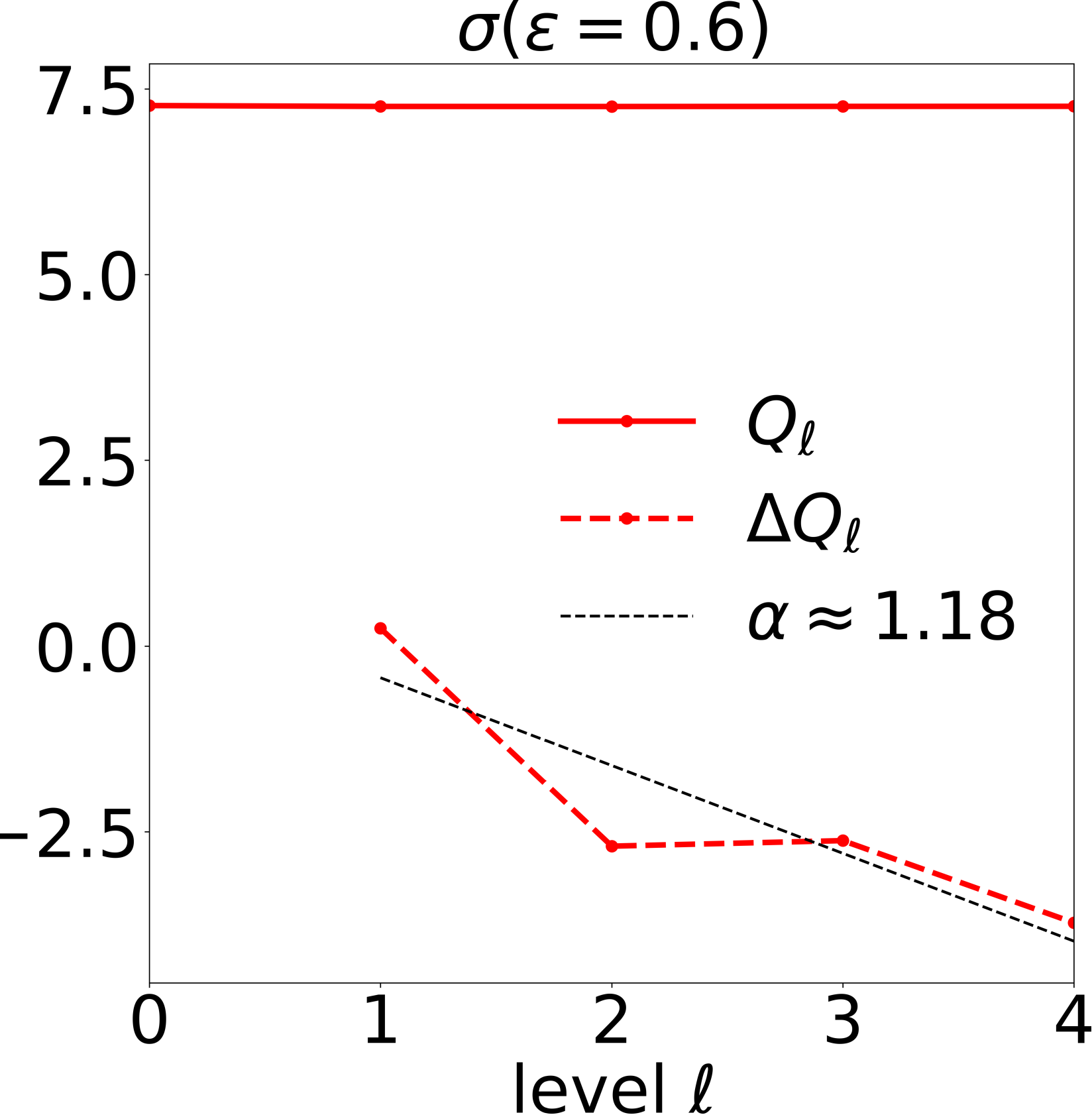}
\caption{$\varepsilon = 0.60$}
\end{subfigure}

\begin{subfigure}[b]{0.30\textwidth}
\includegraphics[width=\textwidth]{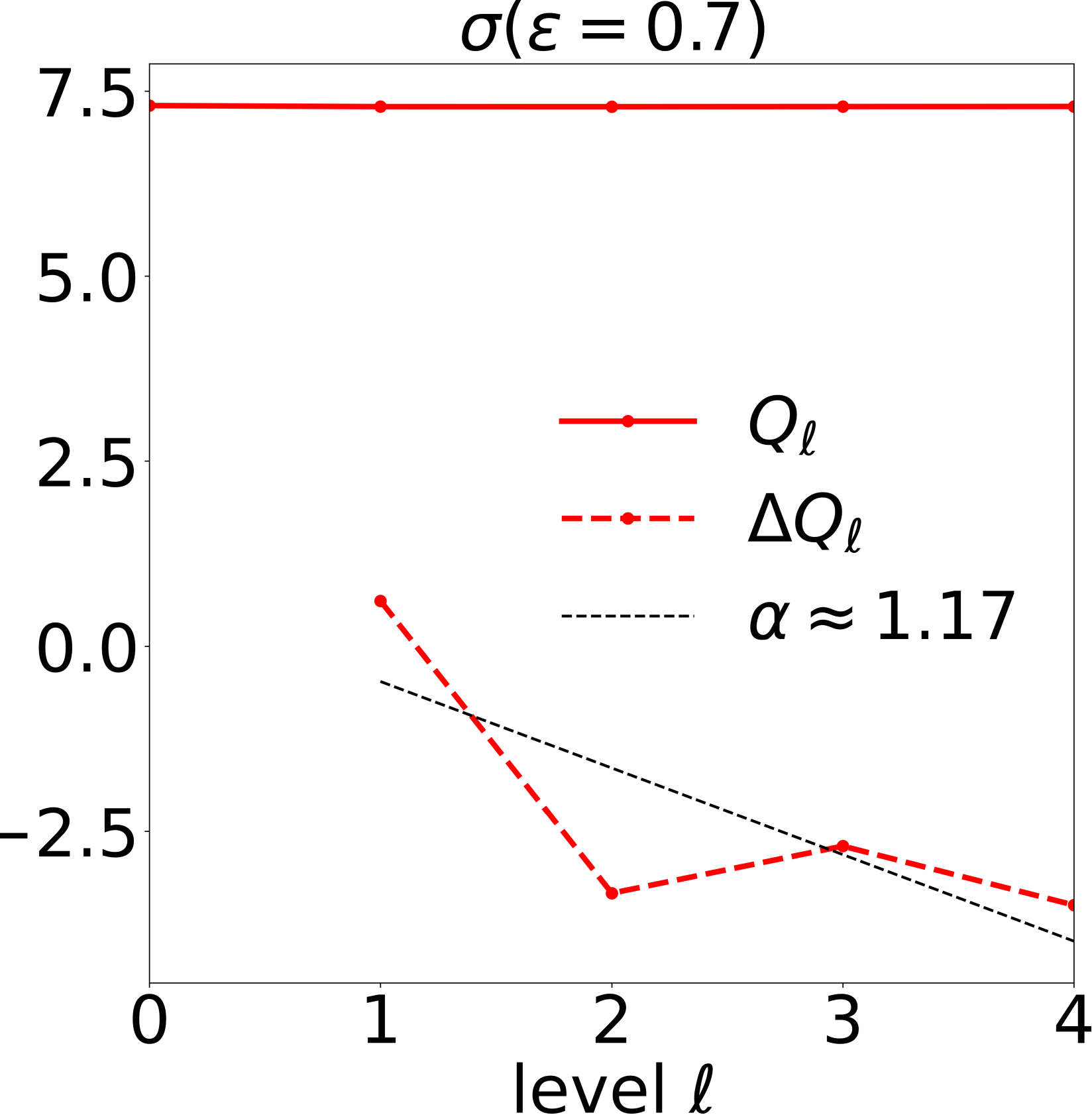}
\caption{$\varepsilon = 0.70$}
\end{subfigure}
\begin{subfigure}[b]{0.30\textwidth}
\includegraphics[width=\textwidth]{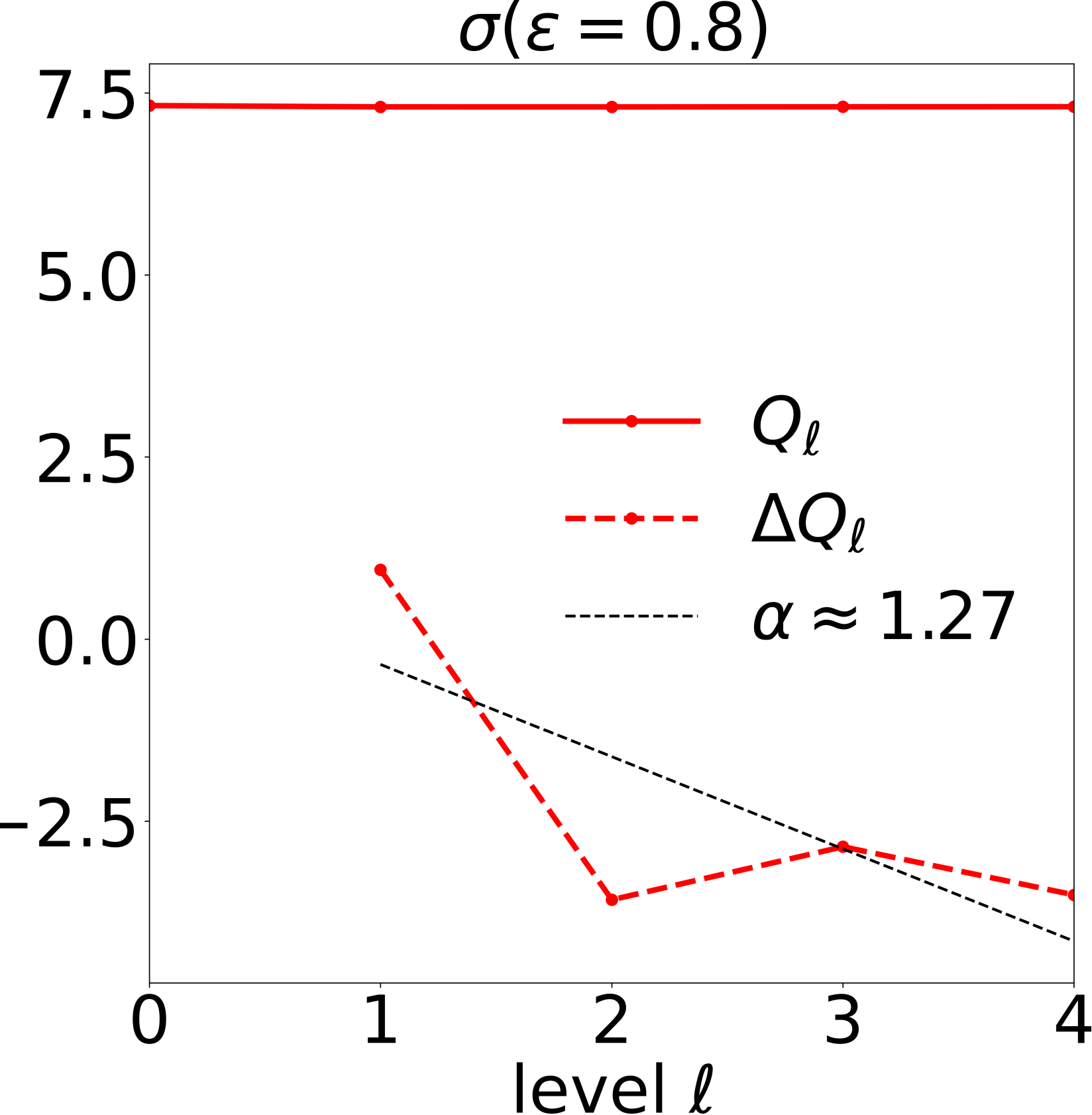}
\caption{$\varepsilon = 0.80$}
\end{subfigure}
\begin{subfigure}[b]{0.30\textwidth}
\includegraphics[width=\textwidth]{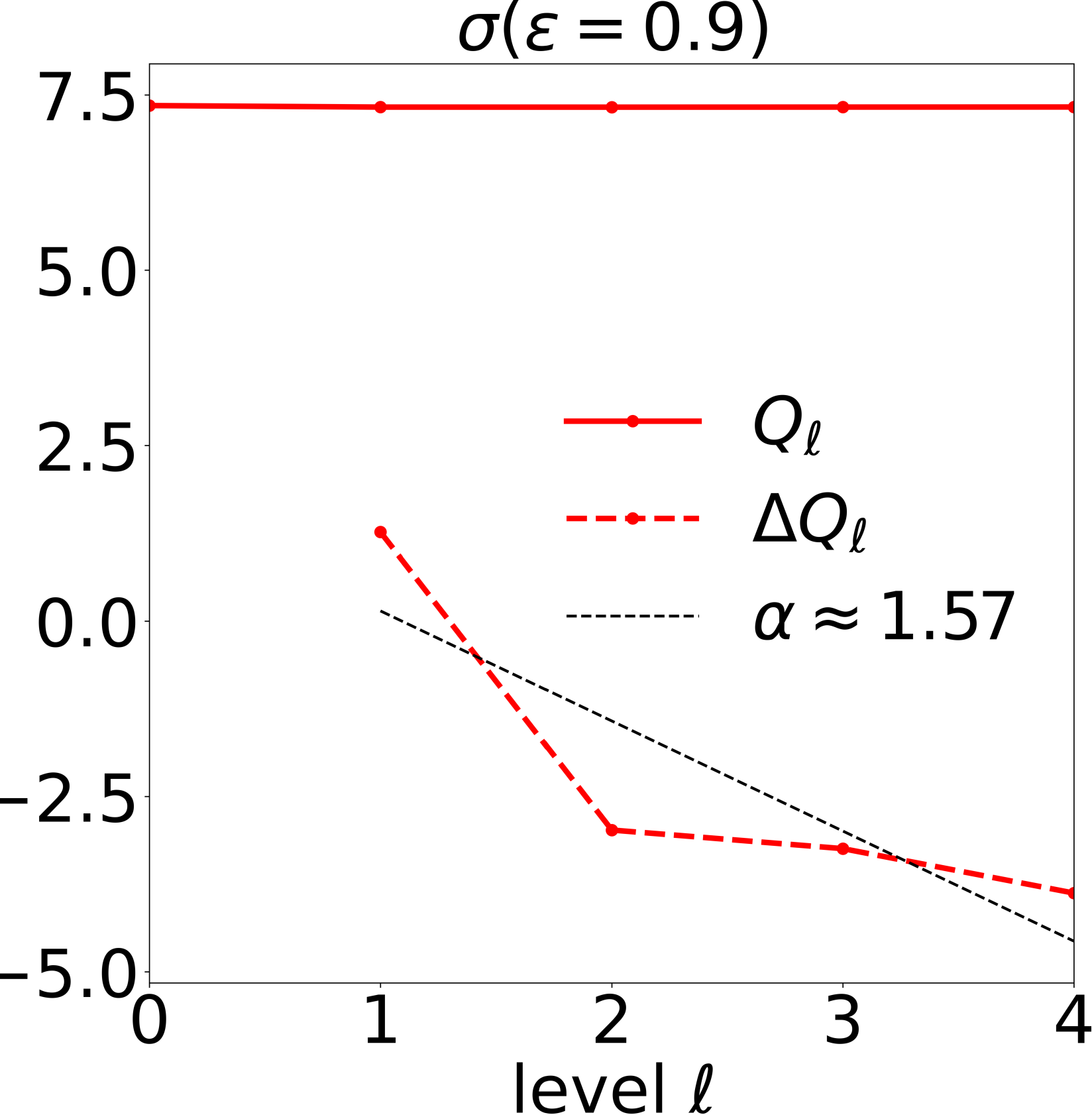}
\caption{$\varepsilon = 0.90$}
\end{subfigure}

\caption{Behavior of the expected value $\log_2 (\mathbb{E} [ \left\vert \cdot \right\vert ])$ of the quantities of interest $Q_\ell$ and the multi-level difference $\Delta Q_\ell$ as a function of the level $\ell$ using the mesh refinement as a fidelity parameter for each QoI. We numerically fitted the values $\bbE[\Delta Q_\ell] \;\propto\; 2^{-\alpha \ell}$ with $1.34 \lesssim \alpha \lesssim 2.93$, see Table~\ref{tab:alpha_beta}.}
\label{fig:mlmc_expected_value}
\end{figure}

\begin{figure}[!htbp]
\centering

\begin{subfigure}[b]{0.30\textwidth}
\includegraphics[width=\textwidth]{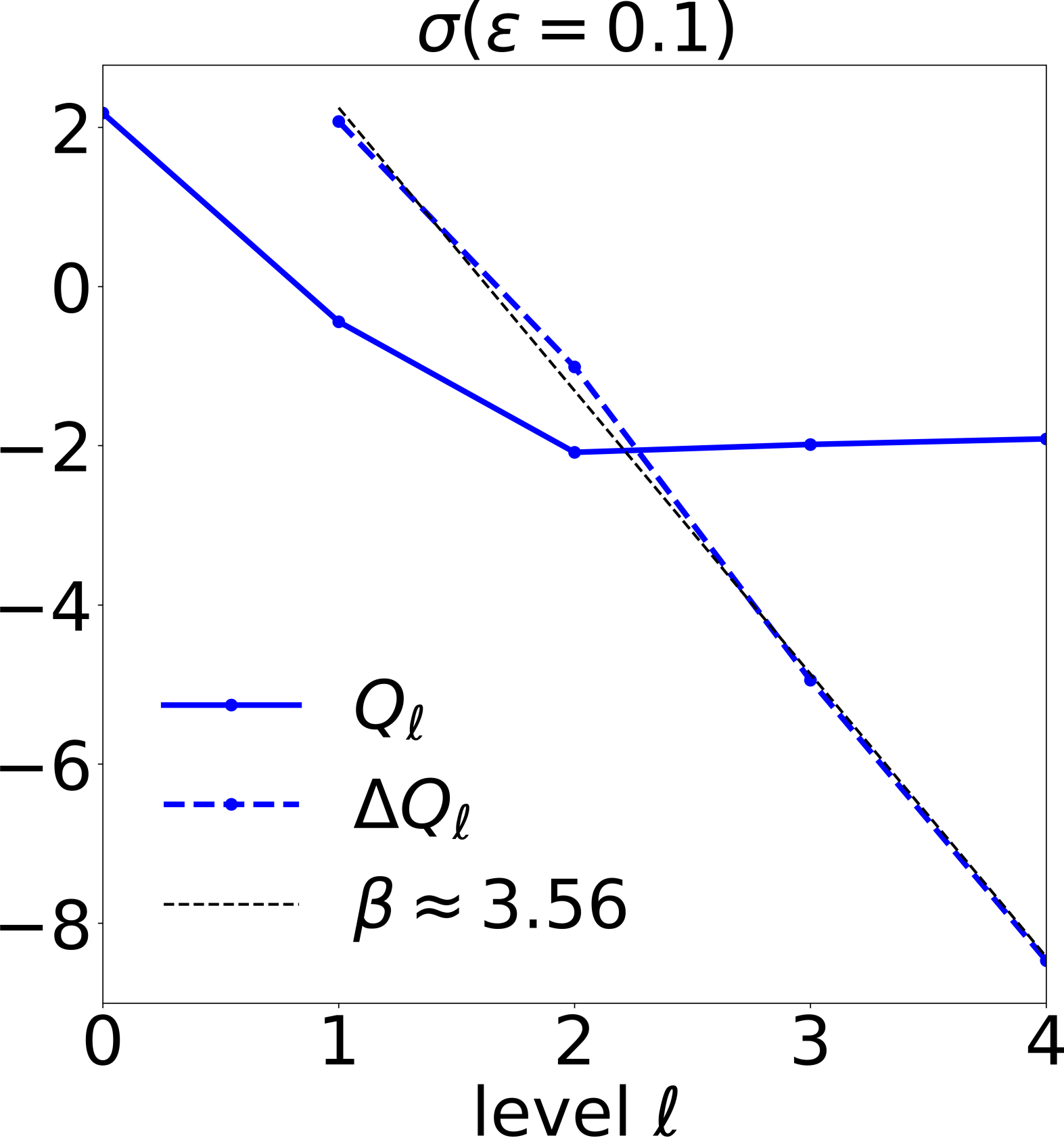}
\caption{$\varepsilon = 0.10$}
\end{subfigure}
\begin{subfigure}[b]{0.30\textwidth}
\includegraphics[width=\textwidth]{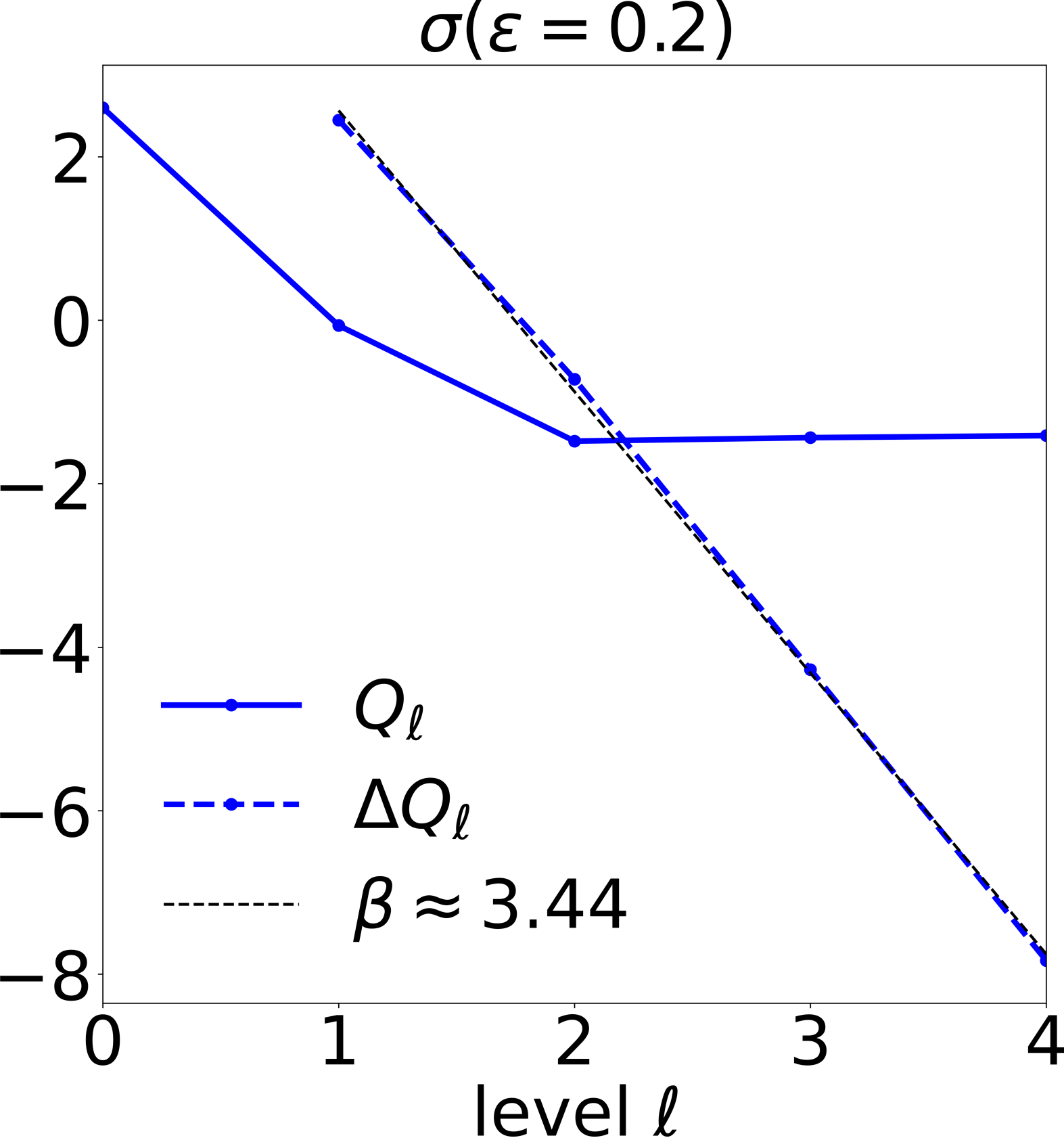}
\caption{$\varepsilon = 0.20$}
\end{subfigure}
\begin{subfigure}[b]{0.30\textwidth}
\includegraphics[width=\textwidth]{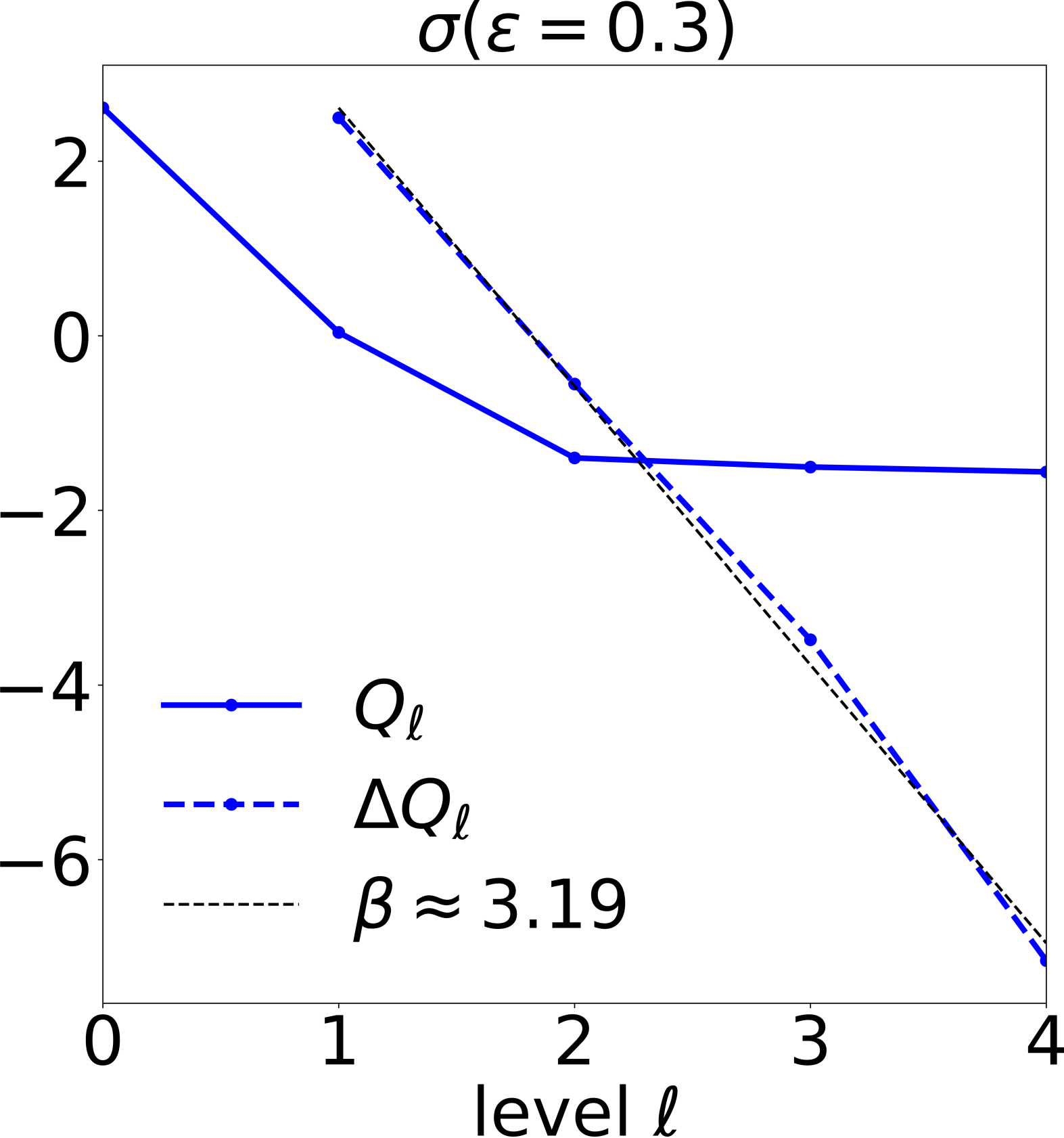}
\caption{$\varepsilon = 0.30$}
\end{subfigure}

\begin{subfigure}[b]{0.30\textwidth}
\includegraphics[width=\textwidth]{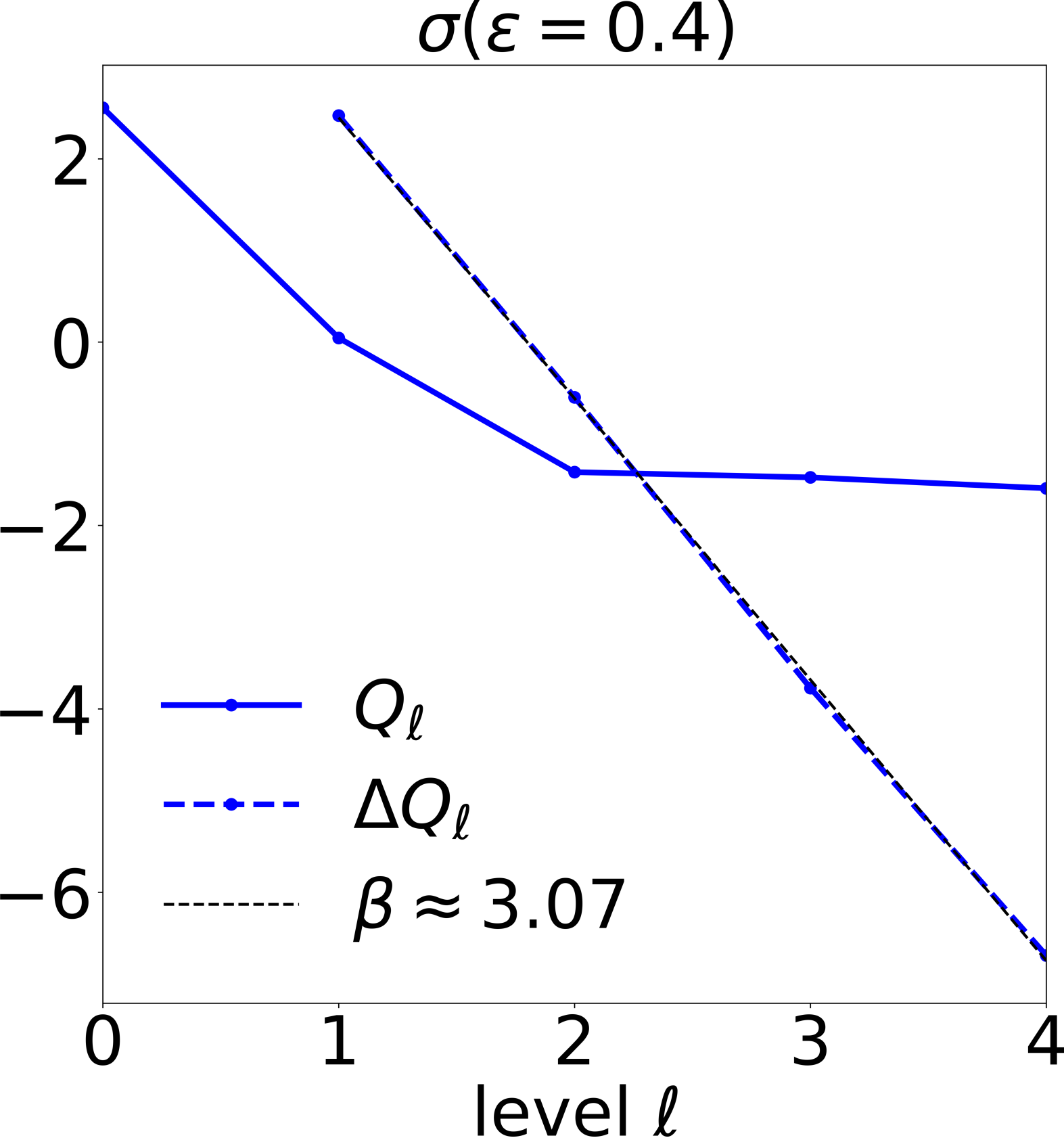}
\caption{$\varepsilon = 0.40$}
\end{subfigure}
\begin{subfigure}[b]{0.30\textwidth}
\includegraphics[width=\textwidth]{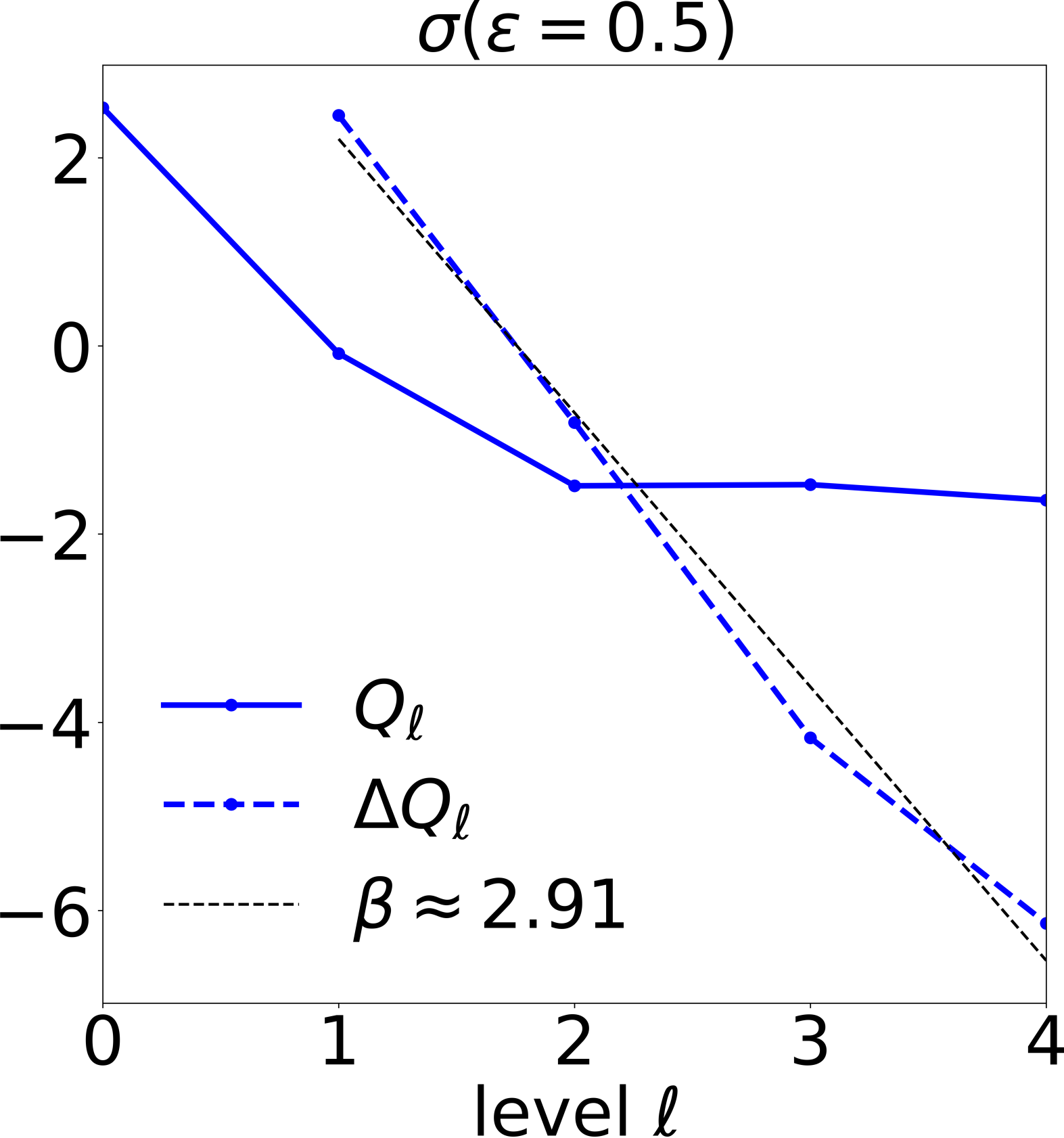}
\caption{$\varepsilon = 0.50$}
\end{subfigure}
\begin{subfigure}[b]{0.30\textwidth}
\includegraphics[width=\textwidth]{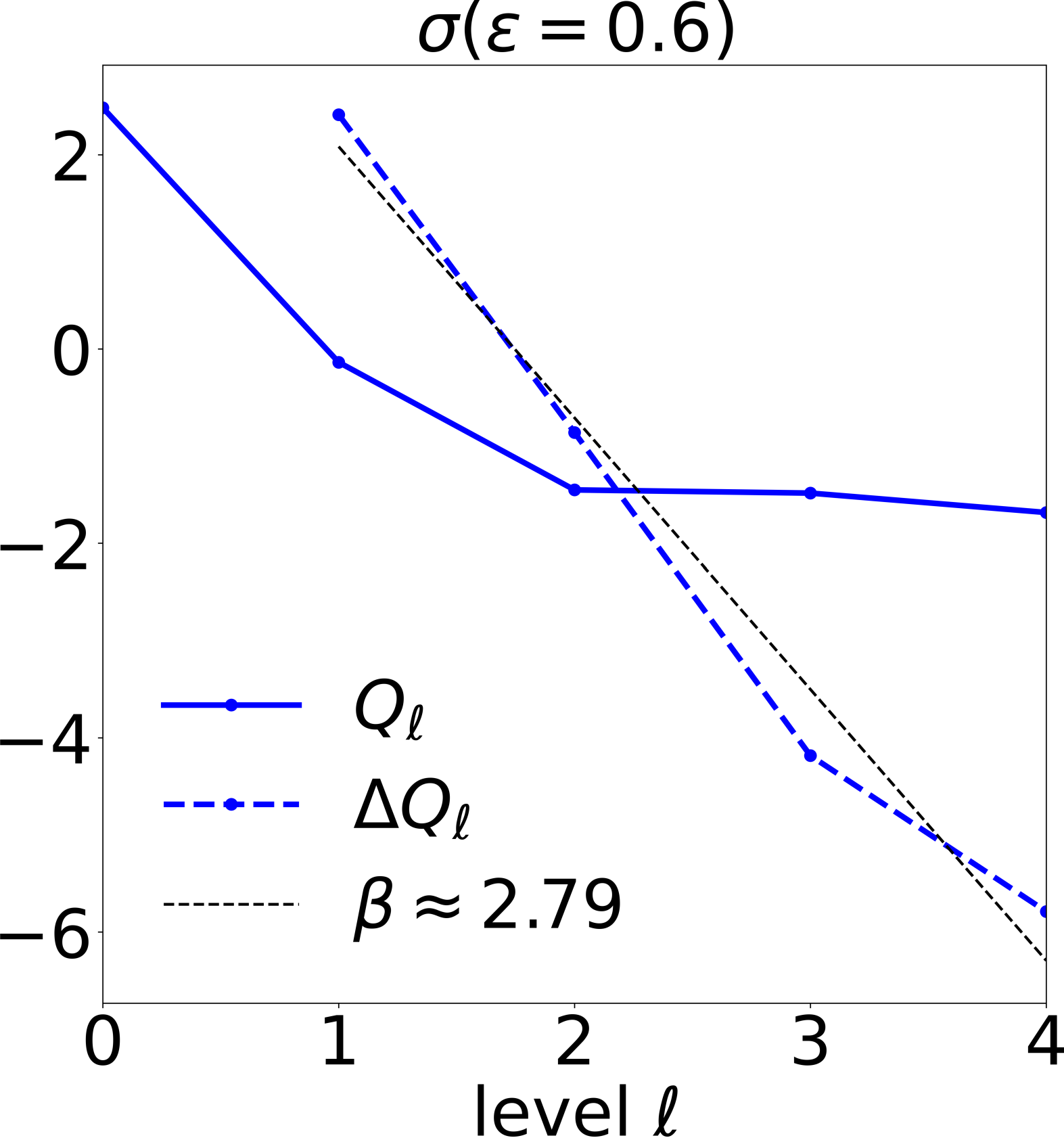}
\caption{$\varepsilon = 0.60$}
\end{subfigure}

\begin{subfigure}[b]{0.30\textwidth}
\includegraphics[width=\textwidth]{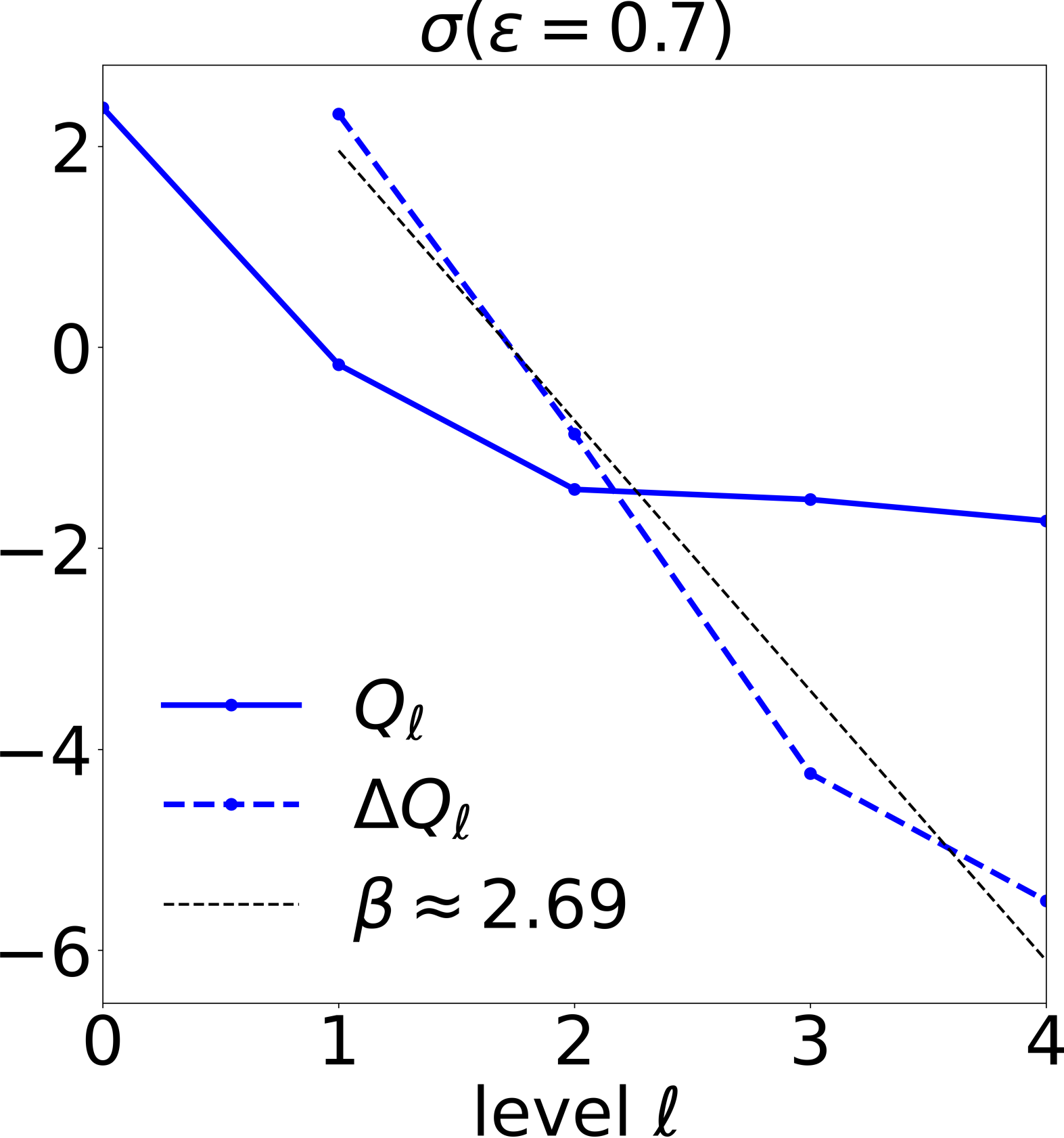}
\caption{$\varepsilon = 0.70$}
\end{subfigure}
\begin{subfigure}[b]{0.30\textwidth}
\includegraphics[width=\textwidth]{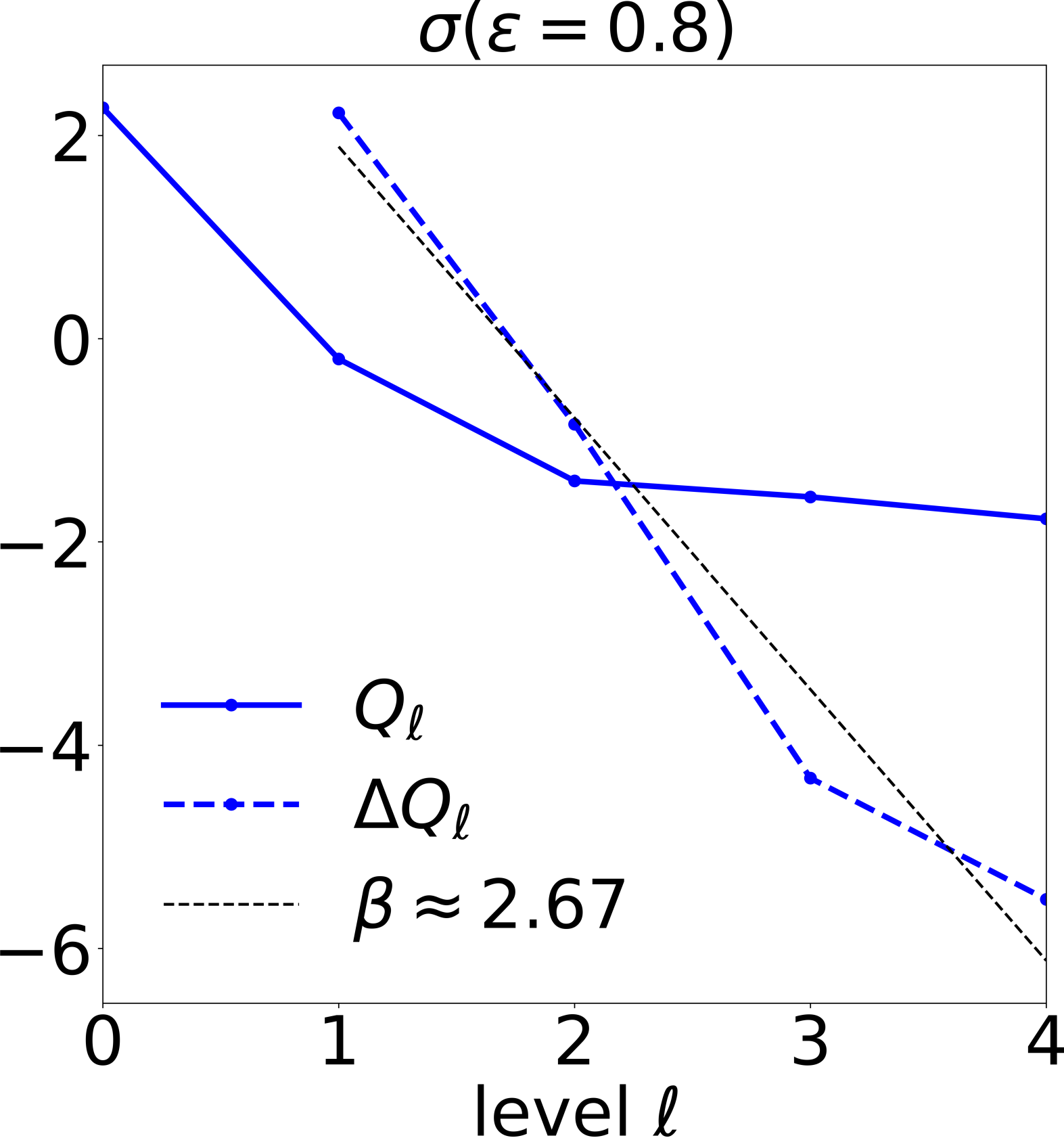}
\caption{$\varepsilon = 0.80$}
\end{subfigure}
\begin{subfigure}[b]{0.30\textwidth}
\includegraphics[width=\textwidth]{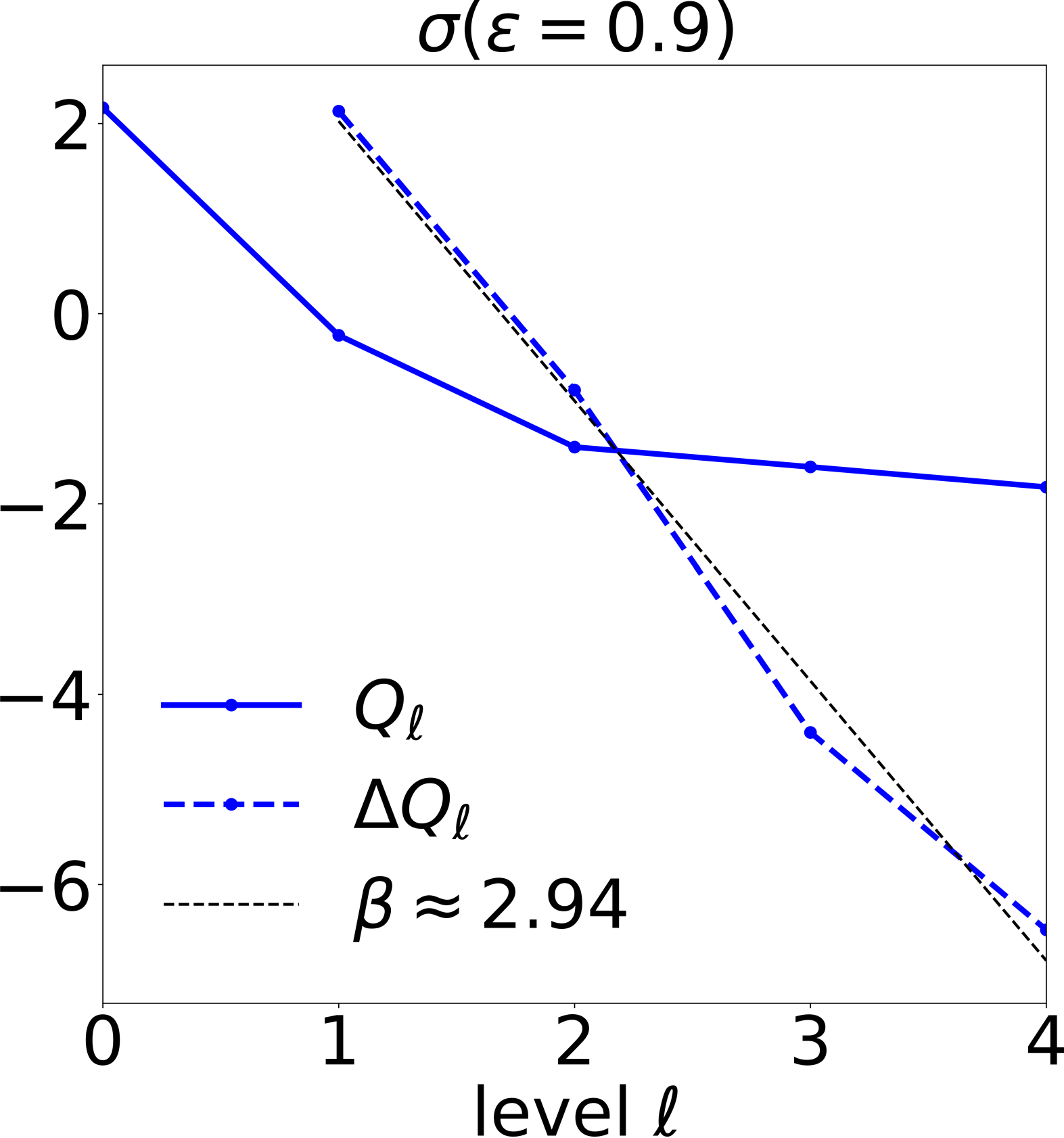}
\caption{$\varepsilon = 0.90$}
\end{subfigure}

\caption{Behavior of the variance $\log_2 (\mathbb{V} [ \left\vert \cdot \right\vert ])$ of the QoI $Q_\ell$ and the multi-level difference $\Delta Q_\ell$ as a function of the level $\ell$ using the mesh refinement as a fidelity parameter for each QoI. We numerically fitted the values $\bbV[\Delta Q_\ell] \;\propto\; 2^{-\beta \ell}$ with $ 2.62 \lesssim \beta \lesssim 3.53$, see Table~\ref{tab:alpha_beta}.}
\label{fig:mlmc_variance}
\end{figure}

The behavior of the expected value of the QoIs $Q_\ell$ and the multilevel difference $Q_\ell - Q_{\ell-1}$ is shown in Figure~\ref{fig:mlmc_expected_value}. The expected value $\E{Q_\ell}$ is stable across all levels $\ell=0, 1, \ldots, 4$, and the expected value of the multilevel difference decays as the mesh resolution level $\ell$ increases. Table~\ref{tab:alpha_beta} shows the numerically fitted decay rates $\alpha$ from condition~\eqref{eq:C1}.

The behavior of the variance of the QoIs $Q_\ell$ and the multilevel difference $Q_\ell - Q_{\ell-1}$ is shown in Figure~\ref{fig:mlmc_variance}. 
As shown in Figure~\ref{fig:stress_strain_compilation_magnesium}, the most materials variability comes from $2 \times 2 \times 2$, and the least materials variability comes from $32 \times 32 \times 32$, possibly due to the number of grains in SERVEs. This observation is consistent with the quantitative analysis of variance shown in Figure~\ref{fig:mlmc_variance}. 
For mesh resolution levels 1 and 2 the variance of the multilevel difference $\V{\Delta Q_\ell}$, is more than the variance of the quantity of interest $\V{Q_\ell}$ itself, i.e. $\V{\Delta Q_\ell} > \V{Q_\ell}$ with $\ell = 1, 2$. This means that the necessary condition for an efficient MLMC estimator required in~\eqref{eq:StrongCorr} is not satisfied (note the logarithmic axis for the variance). A lack of correlation between the QoIs derived from microstructures at levels $\ell =$ 0, 1, and 2 explains this constraint violation, see Figure~\ref{fig:mesh_refinements} for an illustration. However, the variance of the multilevel difference does satisfy condition~\eqref{eq:C2} with numerically fitted rates $\beta$ shown in Table~\ref{tab:alpha_beta}. In order to ensure that our MLMC estimator is efficient, we remove the first three levels in the multilevel hierarchy, leaving only mesh resolution levels $\ell = 3$ and $\ell = 4$. Since $2\alpha \geq \min(\beta, \gamma)$ for all QoIs, we expect an asymptotic cost complexity of $\mathcal{O}(\epsilon^{-2})$, see~\eqref{eq:mlmc_theorem}.

It is observed that the expectation across all fidelity levels $\ell$ is stable (Figure~\ref{fig:mlmc_expected_value}), whereas the variance across all fidelity levels $\ell$, while generally still decreasing, is only useful when $\ell \geq 3$. We postulate that the convergence of the expectation has more to do with the underlying constitutive models, whereas the variance has more to do with the number of grains in SERVEs. The converged expectation is consistent in the notion of \textit{unbiased estimator} (for both MLMC and MC methods, as any MC-based approach is naturally unbiased), but the variance depends on number of SERVEs $N$ (Equation~\eqref{eq:mc_variance_decays}), and implicitly the number of grains. In this sense, there may be a meaningful criteria to establish a bound for low-fidelity representation of SERVEs, so that the correlation constraint in Equation~\eqref{eq:StrongCorr} can be satisfied.

\begin{table}[!htbp]
\centering
\caption{Fitted decay rates $\alpha$ and $\beta$ in conditions~\eqref{eq:C1} and~\eqref{eq:C2}, respectively.}
\label{tab:alpha_beta}
\begin{tabular}{ccc} \hline
$\varepsilon$ & $\alpha$ & $\beta$ \\ \hline
0.1           & 2.93     & 3.53    \\
0.2           & 2.26     & 3.44    \\
0.3           & 2.64     & 3.15    \\
0.4           & 1.59     & 3.02    \\
0.5           & 1.30     & 2.88    \\
0.6           & 1.19     & 2.76    \\
0.7           & 1.16     & 2.64    \\
0.8           & 1.25     & 2.62    \\
0.9           & 1.56     & 2.93    \\ \hline
\end{tabular}
\end{table}

\subsection{Comparison between MC and MLMC}

% \hl{consider add $\mathrm{\mathbb{E}}[\cdot]$ and $\mathrm{\mathbb{V}}[\cdot]$ on $\sigma/\varepsilon$ curves.}

\begin{figure}[!htbp]
\centering
\includegraphics[height=300px,keepaspectratio]{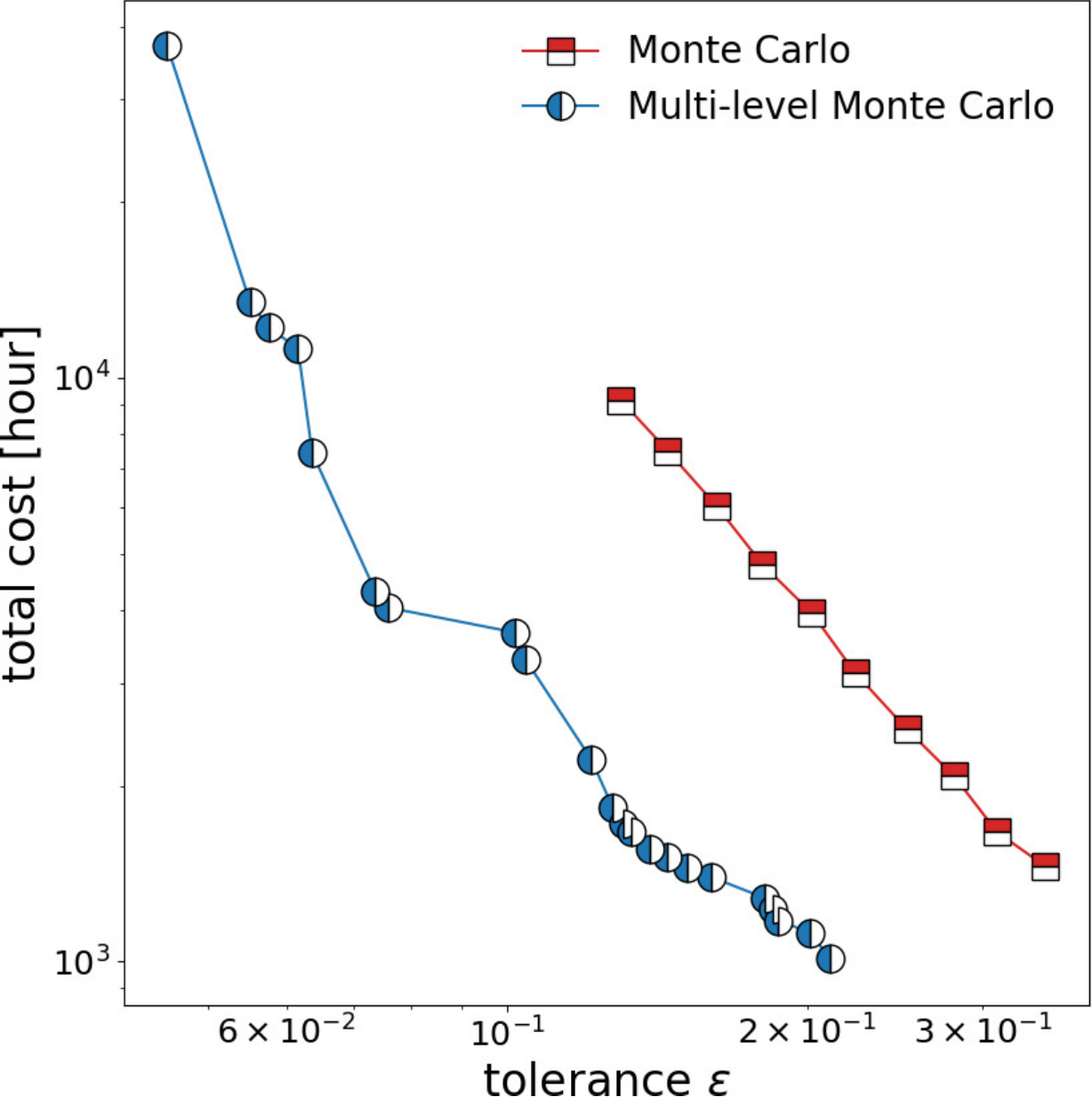}
\caption{Total cost comparison between MC and MLMC showing a factor of around 2.23$\times$ speedup.}
\label{fig:mc-vs-mlmc}
\end{figure}

In Figure~\ref{fig:mc-vs-mlmc}, we show a comparison of the computational cost, measured in wall clock time, of the MC and MLMC methods, for different values of the tolerance parameter $\epsilon$. The MLMC method outperforms the MC method by a factor of 2.23$\times$. For a tolerance $\epsilon = 1.30e-1$, the plain MC method takes 1283.38 hours, while our MLMC method takes approximately 68.14 hours. 
Table~\ref{tab:comparison_mlmc_mc} tabulates the comparison between MLMC and MC methods, where the number of SERVEs at level $\ell = 4$, $N_4$, in the MC method is approximated by $N_3$ in the MLMC method. The reasoning behind this approximation is that we assume the variance on the low-fidelity (coarse) level, i.e. $\ell = 3$ is approximately the same with the high-fidelity (fine) level, i.e. $\ell = 4$. The approximately same variance means approximately the same number of samples, and therefore, the number of SERVEs $N_4$ in the MC method can be approximated by $N_3$ in the MLMC method. Because MF UQ methods generally aim to argue that it is possible to conduct UQ efficiently by multiple levels of fidelity, the best way to highlight the efficiency is simply to compare the MF UQ method at hand with the brute-force high-fidelity counterpart. Therefore, we only consider an approximate $N_4$ and not $N_3$ for the MC method.

\begin{table}[!htbp]
\centering
\caption{Comparison of the distribution of the number of samples across the mesh resolution levels for MLMC and MC, for different tolerances $\epsilon$. 
% \pjr{(COMMENT: Please change this from seconds to hours, it's hard to read otherwise. In fact -- just delete the columns that show the cost, this is already shown in Figure 8.)}
}
\label{tab:comparison_mlmc_mc}
\begin{tabular}{|c|cc|ccc|c|} \hline
\multirow{2}{*}{$\epsilon$} & \multicolumn{2}{c|}{MC} & \multicolumn{3}{c|}{MLMC}           & \multirow{2}{*}{Speedup}   \\ 
                     & Approx. $N_4$    & Approx. Cost (hr)     & $N_3 $   & $N_4$        & Cost (hr)  &              \\ \hline

3.51000000e-01  &  \textcolor{blue}{7 }   &  \textcolor{blue}{ 24.28}   &    7     &   3          &    16.83 & \textcolor{blue}{1.44}$\times$   \\
3.14324193e-01  &  \textcolor{blue}{8 }   &  \textcolor{blue}{ 27.74}   &    8     &   3          &    17.75 & \textcolor{blue}{1.56}$\times$   \\
2.81480622e-01  &  \textcolor{blue}{10}   &  \textcolor{blue}{ 34.68}   &    10    &   3          &    19.58 & \textcolor{blue}{1.77}$\times$  \\
2.52068859e-01  &  \textcolor{blue}{12}   &  \textcolor{blue}{ 41.62}   &    12    &   3          &    21.42 & \textcolor{blue}{1.94}$\times$  \\
2.25730315e-01  &  \textcolor{blue}{15}   &  \textcolor{blue}{ 52.02}   &    15    &   3          &    24.17 & \textcolor{blue}{2.15}$\times$  \\
2.02143872e-01  &  \textcolor{blue}{19}   &  \textcolor{blue}{ 65.90}   &    19    &   3          &    27.84 & \textcolor{blue}{2.36}$\times$   \\
1.81021965e-01  &  \textcolor{blue}{23}   &  \textcolor{blue}{ 79.77}   &    23    &   3          &    31.52 & \textcolor{blue}{2.53}$\times$  \\
1.62107074e-01  &  \textcolor{blue}{29}   &  \textcolor{blue}{100.58}   &    29    &   3          &    37.02 & \textcolor{blue}{2.71}$\times$  \\
1.45168590e-01  &  \textcolor{blue}{36}   &  \textcolor{blue}{124.87}   &    36    &   4          &    46.92 & \textcolor{blue}{2.66}$\times$  \\
1.30000000e-01  &  \textcolor{blue}{44}   &  \textcolor{blue}{152.61}   &    44    &   8          &    68.14 & \textcolor{blue}{2.23}$\times$  \\ \hline

\end{tabular}
\end{table}

\begin{figure}[!htbp]
\centering
\includegraphics[width=1.0\textwidth,keepaspectratio]{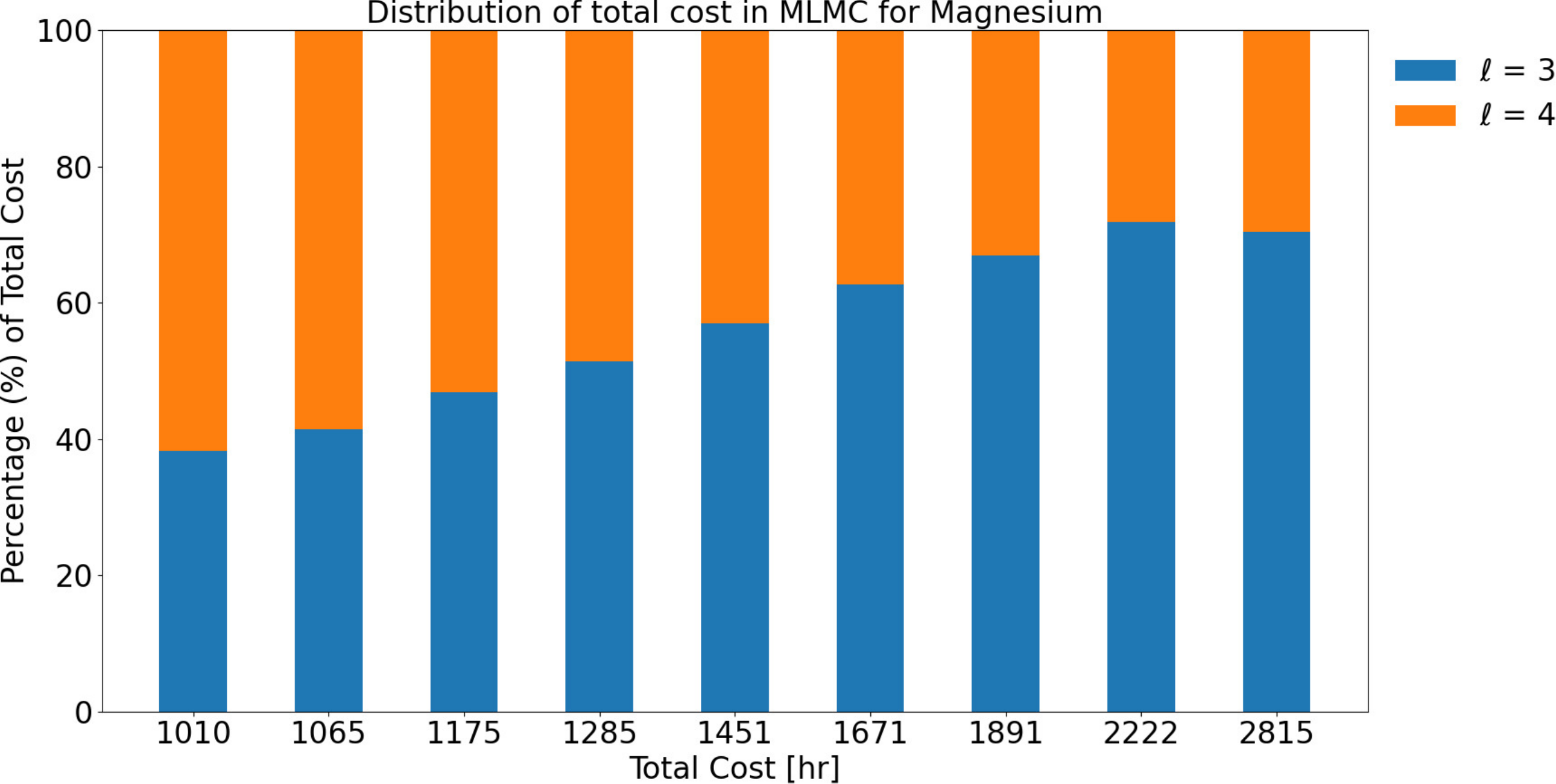}
\caption{Total cost distribution for MLMC run with only levels $\ell = 3, 4$ showing an increase of computational resource allocation toward the low-fidelity level as the total cost increases.}
\label{fig:frac_cost}
\end{figure}

In Figure~\ref{fig:frac_cost} we show the distribution of the total cost among the levels, measured as a percentage of the wall clock time, as a function of the total cost. As the total cost increases, more and more time is spent evaluating the QoIs using a coarse mesh resolution level. For a total cost over 2000 hours, only a quarter of the time is spent evaluating the QoIs at the finest resolution level. 
% \pjr{\sout{percentage distribution of total cost, measured in CPU$\times$hr, as a function of total cost in the MLMC approach, where level $\ell=3$ is plotted as \tikz \path [fill=orange] (0,0) rectangle (1,0.5); and $\ell=4$ are plotted as \tikz \path [fill=blue!55] (0,0) rectangle (1,0.5);. Readers are referred to color version online. As the total cost increases, we also observe an increasing resource allocation to the low-fidelity level $\ell = 3$ to reduce the overall total cost, thanks to a significant correlation across levels.}}

\section{Discussion}
\label{sec:Discussion}

The novelty of this paper lies in the integration of the methodology (i.e. MLMC) and the application (CPFEM) by approaching the classical problem in CPFEM at a different UQ perspective, but not the methodology or the application itself. Given that the interest in UQ for materials science has emerged in the last few years and the ensemble of SERVEs has been the main approach for the last two decades or so, it is non-trivial to formulate the UQ problem through the correct lens, where the right methodology is married to the right application.

Several techniques that exploit a hierarchy of model approximations have been introduced in the recent literature. The multi-index Monte Carlo (MIMC)~\cite{haji2016multi} method, for example, uses a multi-dimensional hierarchy to further increase the efficiency. This multi-dimensional hierarchy may include refinement levels such as constitutive models, integration time-step, or $hp$-FEM, where $h$ corresponds to mesh size and $p$ corresponds to polynomial degree of the finite element. In our previous work in~\cite{tran2023multi} we applied MIMC for a single-QoI CPFEM application. For non-hierarchical model approximations, we mention multi-fidelity MC (MFMC)~\cite{peherstorfer2016optimal,peherstorfer2018survey,qian2018multifidelity} and approximate control variates~\cite{gorodetsky2020generalized,schaden2020multilevel,schaden2021asymptotic,bomarito2022optimization}.

The analysis of the variance in Figure~\ref{fig:mlmc_variance} shows that there is a limitation to the amount of mesh coarsening we can perform in order to yield an efficient multilevel estimator. When the mesh resolution becomes too coarse, the QoIs extracted from the discretized microstructures are no longer correlated, a necessary requirement for a good control variate as per equation~\eqref{eq:StrongCorr}. Such behavior can only be observed through a preliminary analysis of the problem, such as the one performed in this paper. In order to improve the correlation between microstructure resolutions at intermediate mesh resolution levels, the element agglomeration technique from the algebraic multigrid literature may be useful.

There is a lower bound for low-fidelity levels (so that Equation~\eqref{eq:StrongCorr} is not violated), but there is no upper bound for high-fidelity levels. The MLMC method is \textit{at most} much better and \textit{at least} on par with the MC method. The equality in terms of performance occurs when there is only one level of fidelity, as the MLMC method is a natural MF extension of the MC method. 

One important limitation of the MLMC methods is that it requires a substantial correlation across multiple fidelity levels (i.e. $Q_\ell$ and $Q_{\ell-1}$ in \eqref{eq:StrongCorr}). Without sufficient correlation, the variance $\V{\Delta Q_\ell}$ will not be substantially reduced, and thus the gain of using MLMC methods would be marginal compared to MC methods.

% \pjr{\sout{In this work, we have extensively included a wide range from very coarse meshes to fine meshes. However, there is a limitation of mesh coarsening, to which extent the coarsest mesh is still reliable. By analyzing the variance decay with respect to each level $\ell$, Figure~\ref{fig:mlmc_variance} seems to suggest that the number of grains in the coarsest mesh and that in the finest mesh cannot vary significantly from each other. It remains an open question for future research. }}

The generalization of MLMC methods for CPFEM in various materials system is clear, where it can be applied to various materials systems of interests. The MLMC method is a generic framework for MF UQ problems, with the objective of quantifying the uncertainty by leveraging a multi-fidelity hierarchy. There are several mathematical restrictions, such as strong correlation \eqref{eq:StrongCorr} and exponential decay in expectation, variance, and computational cost, i.e. \eqref{eq:C1}, \eqref{eq:C2}, and \eqref{eq:C3}. The exponential decay assumptions, in general, are naturally satisfied as long as the CPFEM converges at all levels $\ell$'s. Perhaps, it is worthy to note that in this work, we focus solely on the case of 1-dimension in fidelity parameters (i.e. mesh resolutions) and multi objectives (i.e. stress observations at 9 collocated points), while in previous work~\cite{tran2023multi}, we focused on multi-dimension in fidelity parameters (i.e. mesh resolutions and constitutive models) with single objective (i.e. Young modulus or yield strength).

It is well-known that there is a close relationship between the materials texture characterized by the orientation distribution functions and the materials stress-strain response. In this work, we follow Mangal and Holm~\cite{mangal2018dataset} to study the texture of magnesium with $(\phi_1, \theta, \phi_2) = (90^\circ,0^\circ,0^\circ)$ as an exemplar for our proposed UQ problem. How material textures influence the material stress-strain responses is certainly a very interesting question, but it lies beyond the scope of this paper, and therefore, is an subject of interest for future study.

Finally, we point out that the CPFEM method, just as all other finite-element-based method, introduces several layers of approximation errors. A first error occurs because of the discretization of the geometry. A second source of error is the accuracy of the sparse linear solver used to compute the finite element solutions. In this work, we have focused on the former error, but it possible to include the solver accuracy as an additional dimension for refinement, possibly in a multi-index Monte Carlo framework.

% \pjr{\sout{For FEM simulations, it is worthy to observe that FEM, including CPFEM, includes two layers of approximations. The first approximation layer is the geometric approximation through a mesh that approximates the microstructure. The second approximation layer comes from the numerical approximation, typically associated with the numerical iterative solvers for solving a large-scale system. As a result, there are multiple sources of epistemic uncertainty in the form of approximation. In this work, we focus on the geometric approximation for SERVE realizations through microstructure reconstruction method, implemented in DREAM.3D}}~\cite{groeber2014dream}. 

\section{Conclusion}
\label{sec:Conclusion}

In this paper, we applied the MLMC method to study the effect of microstructure variation on the stress-strain curves for magnesium predicted by CPFEM. We imposed a geometric mesh hierarchy with five different mesh resolution levels, of which only the two finest levels were used in an efficient multilevel sampling strategy. Our MLMC method outperforms the standard MC method by a factor of 2.23$\times$.  
% \pjr{(COMMENT: change ZZ)}.

% \pjr{\sout{In this paper, we apply the multi-output MLMC method to study the materials variability of stress-strain curves for magnesium with multiple mesh resolutions. 
% Five geometric mesh resolutions are considered, $2 \times 2 \times 2$, $4 \times 4 \times 4$, $8 \times 8 \times 8$, $16 \times 16 \times 16$, and $32 \times 32 \times 32$ with exponentially increasing computational cost. 
% We exclude three coarsest meshes which corresponds to $2 \times 2 \times 2$, $4 \times 4 \times 4$, and $8 \times 8 \times 8$, due to a mathematical constraints in MLMC method. 
% By adaptively sampling between $\ell = 3$ (which corresponds to $16 \times 16 \times 16$ mesh) and $\ell = 4$ (which corresponds to $32 \times 32 \times 32$ mesh), we obtain a speedup factor of 18.83$\times$, compared to MC method with only two levels of fidelity. 
% The numerical results show that MLMC outperforms MC approach in quantifying materials variability. Furthermore, we also show that in this case study, the materials variability seems to increase, even though not always monotonically, as the strain increases into the plasticity domain. 
% This behavior is inconclusive, which begs for future study.  }}

\section*{Acknowledgment}
The views expressed in the article do not necessarily represent the views of the U.S. Department of Energy or the United States Government. Sandia National Laboratories is a multimission laboratory managed and operated by National Technology and Engineering Solutions of Sandia, LLC., a wholly owned subsidiary of Honeywell International, Inc., for the U.S. Department of Energy's National Nuclear Security Administration under contract DE-NA-0003525. 

The authors thank two anonymous reviewers for their critics, which has substantially improved the quality of the manuscript. On behalf of all authors, the corresponding author states that there is no conflict of interest.


% \bibliographystyle{spmpsci}      % mathematics and physical sciences
% \bibliographystyle{spphys}       % APS-like style for physics
% \bibliographystyle{apsrev}       % APS-like style for physics
\bibliography{lib}   % name your BibTeX data base

\begin{thebibliography}{10}
\providecommand{\url}[1]{{#1}}
\providecommand{\urlprefix}{URL }
\providecommand{\doi}[1]{\url{https://doi.org/#1}}
\bibcommenthead

\bibitem{olson1997computational}
G.B. Olson, Computational design of hierarchically structured materials.
\newblock Science \textbf{277}(5330), 1237--1242 (1997)

\bibitem{arroyave2019systems}
R.~Arr{\'o}yave, D.L. McDowell, Systems approaches to materials design: {Past},
  present, and future.
\newblock Annual Review of Materials Research \textbf{49}(1), 103--126 (2019)

\bibitem{hey2009fourth}
T.~Hey, S.~Tansley, K.M. Tolle, \emph{The fourth paradigm: data-intensive
  scientific discovery}, vol.~1 (Microsoft research Redmond, WA, 2009)

\bibitem{agrawal2016perspective}
A.~Agrawal, A.~Choudhary, Perspective: {M}aterials informatics and big data:
  {R}ealization of the ``fourth paradigm'' of science in materials science.
\newblock APL Materials \textbf{4}(5), 053,208 (2016)

\bibitem{national2008integrated}
U.N.R. Council, \emph{Integrated computational materials engineering: a
  transformational discipline for improved competitiveness and national
  security} (National Academies Press, 2008)

\bibitem{national2011materials}
{US NSTC}, \emph{{Materials Genome Initiative} for global competitiveness}
  (Executive Office of the President, National Science and Technology Council,
  2011)

\bibitem{de2019new}
J.J. de~Pablo, N.E. Jackson, M.A. Webb, L.Q. Chen, J.E. Moore, D.~Morgan,
  R.~Jacobs, T.~Pollock, D.G. Schlom, E.S. Toberer, J.~Analytis, I.~Dabo, D.M.
  DeLongchamp, G.A. Fiete, G.M. Grason, G.~Hautier, Y.~Mo, K.~Rajan, E.J. Reed,
  E.~Rodriguez, V.~Stevanovic, J.~Suntivich, K.~Thornton, J.C. Zhao, New
  frontiers for the materials genome initiative.
\newblock npj Computational Materials \textbf{5}(1), 1--23 (2019)

\bibitem{oden2010computer1}
T.~Oden, R.~Moser, O.~Ghattas, Computer predictions with quantified
  uncertainty, {P}art {I}.
\newblock SIAM News \textbf{43}(9), 1--3 (2010)

\bibitem{oden2010computer2}
T.~Oden, R.~Moser, O.~Ghattas, Computer predictions with quantified
  uncertainty, {P}art {II}.
\newblock SIAM News \textbf{43}(10), 1--4 (2010)

\bibitem{nguyen2021bayesian}
T.~Nguyen, D.C. Francom, D.J. Luscher, J.~Wilkerson, Bayesian calibration of a
  physics-based crystal plasticity and damage model.
\newblock Journal of the Mechanics and Physics of Solids \textbf{149}, 104,284
  (2021)

\bibitem{hasan2022microstructure}
M.~Hasan, P.~Acar, Microstructure-sensitive stochastic design of
  polycrystalline materials for quasi-isotropic properties.
\newblock AIAA Journal \textbf{60}(12), 6869--6880 (2022)

\bibitem{tran2019quantifying}
A.~Tran, D.~Liu, H.A. Tran, Y.~Wang, Quantifying uncertainty in the
  process-structure relationship for {Al-Cu} solidification.
\newblock Modelling and Simulation in Materials Science and Engineering
  \textbf{27}(6), 064,005 (2019)

\bibitem{tran2022microstructure}
A.~Tran, T.~Wildey, H.~Lim, Microstructure-sensitive uncertainty quantification
  for crystal plasticity finite element constitutive models using stochastic
  collocation method.
\newblock Frontiers in Materials \textbf{9}, 1--20 (2022)

\bibitem{venkatraman2023new}
A.~Venkatraman, S.~Mohan, R.~Joseph, D.L. McDowell, S.R. Kalidindi, A new
  framework for the assessment of model probabilities of the different crystal
  plasticity models for lamellar grains in $\alpha$+ $\beta$ {T}itanium alloys.
\newblock Modelling and Simulation in Materials Science and Engineering  (2023)

\bibitem{rixner2022self}
M.~Rixner, P.S. Koutsourelakis, Self-supervised optimization of random material
  microstructures in the small-data regime.
\newblock npj Computational Materials \textbf{8}(1), 46 (2022)

\bibitem{tran2020solving}
A.~Tran, T.~Wildey, Solving stochastic inverse problems for property-structure
  linkages using data-consistent inversion and machine learning.
\newblock JOM \textbf{73}, 72--89 (2020)

\bibitem{butler2018combining}
T.~Butler, J.~Jakeman, T.~Wildey, Combining push-forward measures and {Bayes}'
  rule to construct consistent solutions to stochastic inverse problems.
\newblock SIAM Journal on Scientific Computing \textbf{40}(2), A984--A1011
  (2018)

\bibitem{butler2018convergence}
T.~Butler, J.~Jakeman, T.~Wildey, Convergence of probability densities using
  approximate models for forward and inverse problems in uncertainty
  quantification.
\newblock SIAM Journal on Scientific Computing \textbf{40}(5), A3523--A3548
  (2018)

\bibitem{tran2023multi}
A.~Tran, P.~Robbe, H.~Lim, Multi-fidelity microstructure-induced uncertainty
  quantification by advanced monte carlo methods.
\newblock Materialia \textbf{27}, 101,705 (2023)

\bibitem{rodgers2020three}
T.M. Rodgers, H.~Lim, J.A. Brown, Three-dimensional additively manufactured
  microstructures and their mechanical properties.
\newblock JOM \textbf{72}(1), 75--82 (2020)

\bibitem{mitchell2023parallel}
J.A. Mitchell, F.~Abdeljawad, C.~Battaile, C.~Garcia-Cardona, E.A. Holm, E.R.
  Homer, J.~Madison, T.M. Rodgers, A.P. Thompson, V.~Tikare, E.~Webb, S.J.
  Plimpton, Parallel simulation via {SPPARKS} of on-lattice kinetic and
  {Metropolis Monte Carlo} models for materials processing.
\newblock Modelling and Simulation in Materials Science and Engineering
  \textbf{31}(5), 055,001 (2023)

\bibitem{tran2023high}
A.~Tran, H.~Lim, An asynchronous parallel high-throughput model calibration
  framework for crystal plasticity finite element constitutive models.
\newblock Computational Mechanics  (2023)

\bibitem{tran2022aphbo}
A.~Tran, M.~Eldred, T.~Wildey, S.~McCann, J.~Sun, R.J. Visintainer,
  {aphBO-2GP-3B}: a budgeted asynchronous parallel multi-acquisition functions
  for constrained {B}ayesian optimization on high-performing computing
  architecture.
\newblock Structural and Multidisciplinary Optimization \textbf{65}(4), 1--45
  (2022)

\bibitem{ricciardi2020uncertainty}
D.E. Ricciardi, O.A. Chkrebtii, S.R. Niezgoda, Uncertainty quantification
  accounting for model discrepancy within a random effects bayesian framework.
\newblock Integrating Materials and Manufacturing Innovation \textbf{9},
  181--198 (2020)

\bibitem{kennedy2001bayesian}
M.C. Kennedy, A.~O'Hagan, Bayesian calibration of computer models.
\newblock Journal of the Royal Statistical Society: Series B (Statistical
  Methodology) \textbf{63}(3), 425--464 (2001)

\bibitem{khalil2021modeling}
M.~Khalil, G.H. Teichert, C.~Alleman, N.~Heckman, R.E. Jones, K.~Garikipati,
  B.~Boyce, Modeling strength and failure variability due to porosity in
  additively manufactured metals.
\newblock Computer Methods in Applied Mechanics and Engineering \textbf{373},
  113,471 (2021)

\bibitem{ghoreishi2018multi}
S.F. Ghoreishi, A.~Molkeri, A.~Srivastava, R.~Arroyave, D.~Allaire,
  Multi-information source fusion and optimization to realize {ICME:
  Application} to dual-phase materials.
\newblock Journal of Mechanical Design \textbf{140}(11), 111,409 (2018)

\bibitem{groeber2008framework1}
M.~Groeber, S.~Ghosh, M.D. Uchic, D.M. Dimiduk, A framework for automated
  analysis and simulation of {3D} polycrystalline microstructures. {Part 1:
  Statistical characterization}.
\newblock Acta Materialia \textbf{56}(6), 1257--1273 (2008)

\bibitem{groeber2008framework2}
M.~Groeber, S.~Ghosh, M.D. Uchic, D.M. Dimiduk, A framework for automated
  analysis and simulation of {3D} polycrystalline microstructures. {Part 2:
  Synthetic structure generation}.
\newblock Acta Materialia \textbf{56}(6), 1274--1287 (2008)

\bibitem{torquato2002statistical}
S.~Torquato, Statistical description of microstructures.
\newblock Annual review of materials research \textbf{32}(1), 77--111 (2002)

\bibitem{bostanabad2018computational}
R.~Bostanabad, Y.~Zhang, X.~Li, T.~Kearney, L.C. Brinson, D.W. Apley, W.K. Liu,
  W.~Chen, Computational microstructure characterization and reconstruction:
  {R}eview of the state-of-the-art techniques.
\newblock Progress in Materials Science \textbf{95}, 1--41 (2018)

\bibitem{paulson2017reduced}
N.H. Paulson, M.W. Priddy, D.L. McDowell, S.R. Kalidindi, Reduced-order
  structure-property linkages for polycrystalline microstructures based on
  2-point statistics.
\newblock Acta Materialia \textbf{129}, 428--438 (2017)

\bibitem{teferra2018random}
K.~Teferra, L.~Graham-Brady, A random field-based method to estimate
  convergence of apparent properties in computational homogenization.
\newblock Computer Methods in Applied Mechanics and Engineering \textbf{330},
  253--270 (2018)

\bibitem{paulson2018data}
N.H. Paulson, M.W. Priddy, D.L. McDowell, S.R. Kalidindi, Data-driven
  reduced-order models for rank-ordering the high cycle fatigue performance of
  polycrystalline microstructures.
\newblock Materials \& Design \textbf{154}, 170--183 (2018)

\bibitem{robert2013monte}
C.~Robert, G.~Casella, \emph{Monte Carlo statistical methods} (Springer Science
  \& Business Media, 2013)

\bibitem{babuvska1990p}
I.~Babu{\v{s}}ka, M.~Suri, The $p-$ and $h-p$ versions of the finite element
  method, an overview.
\newblock Computer methods in applied mechanics and engineering
  \textbf{80}(1-3), 5--26 (1990)

\bibitem{giles2015multilevel}
M.B. Giles, Multilevel {Monte Carlo} methods.
\newblock Acta Numerica \textbf{24}, 259--328 (2015)

\bibitem{cliffe2011multilevel}
K.A. Cliffe, M.B. Giles, R.~Scheichl, A.L. Teckentrup, Multilevel monte carlo
  methods and applications to elliptic pdes with random coefficients.
\newblock Computing and Visualization in Science \textbf{14}(1), 3--15 (2011)

\bibitem{robbe2017multi}
P.~Robbe, D.~Nuyens, S.~Vandewalle, A multi-index quasi--monte carlo algorithm
  for lognormal diffusion problems.
\newblock SIAM Journal on Scientific Computing \textbf{39}(5), S851--S872
  (2017)

\bibitem{robbe2019recycling}
P.~Robbe, D.~Nuyens, S.~Vandewalle, Recycling samples in the multigrid
  multilevel (quasi-) {Monte Carlo} method.
\newblock SIAM Journal on Scientific Computing \textbf{41}(5), S37--S60 (2019)

\bibitem{groeber2014dream}
M.A. Groeber, M.A. Jackson, {DREAM.3D}: a digital representation environment
  for the analysis of microstructure in {3D}.
\newblock Integrating materials and manufacturing innovation \textbf{3}(1), 5
  (2014)

\bibitem{roters2019damask}
F.~Roters, M.~Diehl, P.~Shanthraj, P.~Eisenlohr, C.~Reuber, S.~Wong, T.~Maiti,
  A.~Ebrahimi, T.~Hochrainer, H.O. Fabritius, S.~Nikolov, M.~Fri\'{a}k,
  N.~Fujita, N.~Grilli, K.~Janssens, N.~Jia, P.~Kok, D.~Ma, F.~Meier,
  E.~Werner, M.~Stricker, D.~Weygand, D.~Raabe, {DAMASK--The D{\"u}sseldorf
  Advanced Material Simulation Kit for modeling multi-physics crystal
  plasticity, thermal, and damage phenomena from the single crystal up to the
  component scale}.
\newblock Computational Materials Science \textbf{158}, 420--478 (2019)

\bibitem{balay1997efficient}
S.~Balay, W.D. Gropp, L.C. McInnes, B.F. Smith, in \emph{Modern software tools
  for scientific computing} (Springer, 1997), pp. 163--202

\bibitem{petsc-web-page}
S.~Balay, S.~Abhyankar, M.F. Adams, S.~Benson, J.~Brown, P.~Brune,
  K.~Buschelman, E.M. Constantinescu, L.~Dalcin, A.~Dener, V.~Eijkhout,
  J.~Faibussowitsch, W.D. Gropp, V.~Hapla, T.~Isaac, P.~Jolivet, D.~Karpeev,
  D.~Kaushik, M.G. Knepley, F.~Kong, S.~Kruger, D.A. May, L.C. McInnes, R.T.
  Mills, L.~Mitchell, T.~Munson, J.E. Roman, K.~Rupp, P.~Sanan, J.~Sarich, B.F.
  Smith, S.~Zampini, H.~Zhang, H.~Zhang, J.~Zhang.
\newblock {PETS}c {W}eb page.
\newblock \url{https://petsc.org/} (2023).
\newblock \urlprefix\url{https://petsc.org/}

\bibitem{roters2010overview}
F.~Roters, P.~Eisenlohr, L.~Hantcherli, D.D. Tjahjanto, T.R. Bieler, D.~Raabe,
  Overview of constitutive laws, kinematics, homogenization and multiscale
  methods in crystal plasticity finite-element modeling: {T}heory, experiments,
  applications.
\newblock Acta Materialia \textbf{58}(4), 1152--1211 (2010)

\bibitem{roters2011crystal}
F.~Roters, P.~Eisenlohr, T.R. Bieler, D.~Raabe, \emph{Crystal plasticity finite
  element methods: in materials science and engineering} (John Wiley \& Sons,
  2011)

\bibitem{sedighiani2020efficient}
K.~Sedighiani, M.~Diehl, K.~Traka, F.~Roters, J.~Sietsma, D.~Raabe, An
  efficient and robust approach to determine material parameters of crystal
  plasticity constitutive laws from macro-scale stress--strain curves.
\newblock International Journal of Plasticity \textbf{134}, 102,779 (2020)

\bibitem{sedighiani2022determination}
K.~Sedighiani, K.~Traka, F.~Roters, D.~Raabe, J.~Sietsma, M.~Diehl,
  Determination and analysis of the constitutive parameters of
  temperature-dependent dislocation-density-based crystal plasticity models.
\newblock Mechanics of Materials \textbf{164}, 104,117 (2022)

\bibitem{wang2014situ}
F.~Wang, S.~Sandl{\"o}bes, M.~Diehl, L.~Sharma, F.~Roters, D.~Raabe, In situ
  observation of collective grain-scale mechanics in {Mg and Mg--rare} earth
  alloys.
\newblock Acta materialia \textbf{80}, 77--93 (2014)

\bibitem{tromans2011elastic}
D.~Tromans, Elastic anisotropy of {HCP} metal crystals and polycrystals.
\newblock Int. J. Res. Rev. Appl. Sci \textbf{6}(4), 462--483 (2011)

\bibitem{agnew2006validating}
S.~Agnew, D.~Brown, C.~Tom{\'e}, Validating a polycrystal model for the
  elastoplastic response of magnesium alloy {AZ31} using in situ neutron
  diffraction.
\newblock Acta materialia \textbf{54}(18), 4841--4852 (2006)

\bibitem{mangal2018dataset}
A.~Mangal, E.A. Holm, A dataset of synthetic hexagonal close packed 3d
  polycrystalline microstructures, grain-wise microstructural descriptors and
  grain averaged stress fields under uniaxial tensile deformation for two sets
  of constitutive parameters.
\newblock Data in brief \textbf{21}, 1833--1841 (2018)

\bibitem{virtanen2020scipy}
P.~Virtanen, R.~Gommers, T.E. Oliphant, M.~Haberland, T.~Reddy, D.~Cournapeau,
  E.~Burovski, P.~Peterson, W.~Weckesser, J.~Bright, S.J. van~der Walt,
  M.~Brett, J.~Wilson, K.J. Millman, N.~Mayorov, A.R.J. Nelson, E.~Jones,
  R.~Kern, E.~Larson, C.J. Carey, I.~Polat, Y.~Feng, E.W. Moore, J.~VanderPlas,
  D.~Laxalde, J.~Perktold, R.~Cimrman, I.~Henriksen, E.A. Quintero, C.R.
  Harris, A.M. Archibald, A.H. Ribeiro, F.~Pedregosa, P.~van Mulbregt,
  A.~Vijaykumar, A.P. Bardelli, A.~Rothberg, A.~Hilboll, A.~Kloeckner,
  A.~Scopatz, A.~Lee, A.~Rokem, C.N. Woods, C.~Fulton, C.~Masson,
  C.~H\"{a}ggstr\"{o}m, C.~Fitzgerald, D.A. Nicholson, D.R. Hagen, D.V.
  Pasechnik, E.~Olivetti, E.~Martin, E.~Wieser, F.~Silva, F.~Lenders,
  F.~Wilhelm, G.~Young, G.A. Price, G.L. Ingold, G.E. Allen, G.R. Lee,
  H.~Audren, I.~Probst, J.P. Dietrich, J.~Silterra, J.T. Webber, J.~Slavič,
  J.~Nothman, J.~Buchner, J.~Kulick, J.L. Schönberger, J.V.
  de~Miranda~Cardoso, J.~Reimer, J.~Harrington, J.L.C. Rodríguez,
  J.~Nunez-Iglesias, J.~Kuczynski, K.~Tritz, M.~Thoma, M.~Newville,
  M.~Kümmerer, M.~Bolingbroke, M.~Tartre, M.~Pak, N.J. Smith, N.~Nowaczyk,
  N.~Shebanov, O.~Pavlyk, P.A. Brodtkorb, P.~Lee, R.T. McGibbon, R.~Feldbauer,
  S.~Lewis, S.~Tygier, S.~Sievert, S.~Vigna, S.~Peterson, S.~More, T.~Pudlik,
  T.~Oshima, T.J. Pingel, T.P. Robitaille, T.~Spura, T.R. Jones, T.~Cera,
  T.~Leslie, T.~Zito, T.~Krauss, U.~Upadhyay, Y.O. Halchenko,
  Y.~Vázquez-Baeza, {SciPy 1.0 Contributors}, {SciPy 1.0}: fundamental
  algorithms for scientific computing in {P}ython.
\newblock Nature methods \textbf{17}(3), 261--272 (2020)

\bibitem{haji2016multi}
A.L. Haji-Ali, F.~Nobile, R.~Tempone, Multi-index {Monte Carlo}: when sparsity
  meets sampling.
\newblock Numerische Mathematik \textbf{132}(4), 767--806 (2016)

\bibitem{peherstorfer2016optimal}
B.~Peherstorfer, K.~Willcox, M.~Gunzburger, Optimal model management for
  multifidelity {Monte Carlo} estimation.
\newblock SIAM Journal on Scientific Computing \textbf{38}(5), A3163--A3194
  (2016)

\bibitem{peherstorfer2018survey}
B.~Peherstorfer, K.~Willcox, M.~Gunzburger, Survey of multifidelity methods in
  uncertainty propagation, inference, and optimization.
\newblock SIAM Review \textbf{60}(3), 550--591 (2018)

\bibitem{qian2018multifidelity}
E.~Qian, B.~Peherstorfer, D.~O'Malley, V.V. Vesselinov, K.~Willcox,
  Multifidelity {Monte Carlo} estimation of variance and sensitivity indices.
\newblock SIAM/ASA Journal on Uncertainty Quantification \textbf{6}(2),
  683--706 (2018)

\bibitem{gorodetsky2020generalized}
A.A. Gorodetsky, G.~Geraci, M.S. Eldred, J.D. Jakeman, A generalized
  approximate control variate framework for multifidelity uncertainty
  quantification.
\newblock Journal of Computational Physics \textbf{408}, 109,257 (2020)

\bibitem{schaden2020multilevel}
D.~Schaden, E.~Ullmann, On multilevel best linear unbiased estimators.
\newblock SIAM/ASA Journal on Uncertainty Quantification \textbf{8}(2),
  601--635 (2020)

\bibitem{schaden2021asymptotic}
D.~Schaden, E.~Ullmann, Asymptotic analysis of multilevel best linear unbiased
  estimators.
\newblock SIAM/ASA Journal on Uncertainty Quantification \textbf{9}(3),
  953--978 (2021)

\bibitem{bomarito2022optimization}
G.F. Bomarito, P.E. Leser, J.E. Warner, W.P. Leser, On the optimization of
  approximate control variates with parametrically defined estimators.
\newblock Journal of Computational Physics \textbf{451}, 110,882 (2022)

\end{thebibliography}

\end{document}